\documentclass[11pt]{article}

\usepackage{amsthm,amsmath,amssymb,bm}
\usepackage{graphicx,color,epsfig}
\usepackage[margin=1in]{geometry}
\usepackage{xr}
\newtheorem{Thm}{Theorem}[section]
\newtheorem{Lem}[Thm]{Lemma}
\newtheorem{Pro}[Thm]{Proposition}
\newtheorem{Claim}[Thm]{Claim}
\newtheorem{Stat}[Thm]{Statement}
\newtheorem{CondProb}[Thm]{Conditioning Problem}
\newtheorem{Cor}[Thm]{Corollary}

\newcommand{\N}{\mathbb{N}}
\renewcommand{\P}{\mathbb{P}}
\newcommand{\R}{\mathbb{R}}
\newcommand{\T}{\mathbb{T}}
\newcommand{\Z}{\mathbb{Z}}

\begin{document}

\title{Conditioning problems for invariant sets of expanding piecewise affine mappings: Application to loss of ergodicity in globally coupled maps}
\author{Bastien Fernandez and Fanni M.\ S\'elley}
\date{}   
\maketitle
\begin{center}
Laboratoire de Probabilit\'es, Statistique et Mod\'elisation\\
CNRS - Univ. Paris Cit\'e-  Sorbonne Univ.\\
Paris, France\\
fernandez@lpsm.paris and selley@lpsm.paris
\end{center}

\begin{abstract}
We propose a systematic approach to the construction of invariant union of polytopes (IUP) in expanding piecewise affine mappings whose linear components are isotropic scalings. The approach relies on using empirical information embedded in trajectories in order to infer, and then to solve, a so-called conditioning problem for some generating collection of polytopes. A conditioning problem consists of a series of requirements on the polytopes' localisation and on the dynamical transitions between these elements. The core element of the approach is a reformulation of the problem as a set of piecewise linear inequalities for some matrices which encapsulate geometric constraints. In that way, the original topological puzzle is converted into a standard problem in computational geometry. This transformation involves an optimization procedure that ensures that both problems are equivalent.

As a proof of concept, the approach is applied to the study of the loss of ergodicity in basic examples of globally coupled maps. The study explains, completes and substantially extends previous achievements about asymmetric IUP in these systems.
Comparison with the numerics reveals sharp existence conditions depending on the map parameters, and accurate fits of the empirical ergodic components. In addition, this application also reveals unanticipated features about conditioning problem solutions, especially as the dependence on the set of admissible face directions is concerned.
\end{abstract}

\leftline{\small\today.}

\section{Introduction}
The main motivation of this work resides in proving the loss of ergodicity in expanding piecewise affine systems of globally coupled maps. More precisely, this means the  emergence of several ergodic components of positive Lebesgue measure when the coupling strength increases, from an ergodic absolutely continuous invariant measure at weak coupling. (NB: An ergodic component of positive Lebesgue measure is called a {\bf Lebesgue ergodic component} for short.)

Systems of coupled maps are deterministic models of collective systems of interacting units. They have revealed a rich phenomenology, depending upon the coupling strength and the characteristics of the individual dynamics, while being amenable to some mathematical analysis  \cite{CF05,K93}. In particular, in the case of piecewise expanding dynamics, perturbative results have been obtained about topological and ergodic properties at (very) weak coupling, in the regime where ergodicity and mixing hold by continuation from the uncoupled limit, see e.g.\ \cite{AF00,BS88,KL06}. For larger coupling strengths, the analogy with particle systems in Statistical Mechanics suggests that ergodicity should be eventually lost via some analogue of a symmetry-breaking-induced phase transition \cite{BS88}. Indeed, inspired by Toom's cellular automaton, some ad-hoc examples of infinite coupled map lattices have been constructed, that exhibit such phase transitions  \cite{BK06,GM00}. Besides, from the dynamical systems viewpoint, examples of bifurcations of a symmetric attractor with positive volume \cite{CG88,TSU94} suggest that ergodicity may be also lost in systems with finitely many units.

In fact, transitions from ergodicity to the emergence of several Lebesgue ergodic components have been observed and proved for piecewise affine coupled maps with a small number $N$ of units \cite{F14,S18,SB16}. In addition, a computer-assisted proof has been developed in \cite{F20}, which in principle applies to any $N\in\N$. However, in practice, its implementation turns out to be computing resource-intensive. So far, the proposed construction could only be completed for $N$ up to 6. Therefore, more effective approaches remain to be provided in order to envisage addressing (very) large numbers of units.

While ergodicity is a central notion in dynamical systems \cite{KH95}, in particular as foundations of Statistical Mechanics are concerned \cite{G99}, and notwithstanding the various conceptual criteria for ergodicity that have been provided, such as the existence of a transitive orbit, no universal method exists to establish this property in an arbitrary system.
Nonetheless, examples have been given of proved ergodicity in non-trivial parameter families \cite{BY93}.

Similarly, no universal method exists to establish absence of ergodicity. However, in the case of expanding piecewise affine maps of $\R^d$ whose affine domains are polytopes, a simple way to establish the existence of several Lebesgue ergodic components is to build up dynamically invariant unions of (sub-)polytopes (IUP) that contain at least one but surely not all ergodic components. A particularly relevant subcase is when the maps commute with some $\Z_2$-symmetry. Then, it suffices to prove the existence of IUP that are disjoint from their symmetric image (AsIUP). (NB: Accurate definitions of IUP and AsIUP are given in Section \ref{S-CONPRO} below.)

Piecewise affine mappings of $\R^d$ with polytope domains must have finitely many Lebesgue ergodic components \cite{T01} (see also \cite{B97}), yet non-trivial\footnote{Non-trivial means distinct from the union of all domains.} IUP need not exist and no universal approach is at hand for their construction. However, the {\sl ad-hoc} construction of AsIUP in coupled maps in \cite{F14,S18,SB16} suggests that numerical simulations of orbits may contain enough relevant information to infer IUP. Indeed, those AsIUP have been intuited using knowledge about the location and dynamics of the empirical ergodic components.

Based on these insights, the purpose of this paper is to develop a systematic approach to the construction of IUP in piecewise affine mappings of $\R^d$ whose linear components are isotropic scalings, and whose atoms are open convex and bounded polytopes. The approach firstly consists in expressing any IUP as a collection of polytopes that satisfy certain topological conditions, a so-called {\bf conditioning problem}. Then, it aims to obtain solutions of this problem. The topological conditions are inspired from the empirical information contained in numerical trajectories. They specify the polytopes' location in the atomic partition and the location of the corresponding images under the dynamics (NB: naturally, the images are assumed to be contained in the IUP, in order to ensure dynamical invariance; see details in Section \ref{S-CONPRO}). 

To address conditioning problems implies manipulating polytopes in arbitrary dimension. To that goal, we find it convenient to regard polytopes as intersections of half-spaces. More precisely, polytopes will be represented using tables ({\bf constraint matrices}) that collect information about the direction and location of the constraining hyper-planes. In addition, constraint matrices will be equipped with an optimization procedure, which ensures a sharp description of the polytopes (namely all constraints are made active); and hence one-to-one correspondence between optimized constraint matrices and polytopes.

Furthermore, basic topological and geometric operations on polytopes will be expressed in terms of operations on (optimized) constraint matrices, so that the topological conditions of a conditioning problem will be converted into inequalities on matrix entries. In this way, any conditioning problem will be reformulated as the problem of finding a collection of constraint matrices whose entries satisfy certain multidimensional piecewise linear inequalities. In order to obtain the desired IUP, it will then remain to find solutions of these inequalities, via standard analysis of piecewise linear problems.

The most simple instance of application of this approach is in dimension one and involves either interval exchange transformations \cite{K75} (see also \cite{K77,KN76} for examples with multiple Lebesgue ergodic components) or piecewise expanding interval maps. IUP are invariant unions of intervals in this case.
The analytic formulation of a conditioning problem intends to determine the interval boundaries from imposing the location of the intervals and of their images in the IUP.

Instead, application in this paper focuses on establishing AsIUP in basic examples of expanding globally coupled maps. In a way, the resulting mathematical statements can be regarded as a reconsideration of the results in \cite{F14,S18,SB16}, which provides previously missing justifications and more thorough descriptions. In practice, the statements also consider the maps associated with more general population distributions (given that in the original coupled maps, the population distribution is uniform, see Section \ref{S-31} for the related definitions).  
For a detailed description of the contributions, we refer to the introductory paragraphs of the various subsections in Section \ref{S-COUPLEDMAP}.

This application interestingly reveals unanticipated features about conditioning problem solutions depending on the set of admissible face directions. In particular, it may happen that no solution exists when the face directions are limited to those given by the atoms of the map. Yet solutions exist when more directions are allowed, which are determined upon solving. 

The coupled map examples are families of mappings that are parametrized by a number $\epsilon\in (0,\tfrac12)$ which quantifies the coupling strength. In particular, ergodicity holds when $\epsilon$ is sufficiently close to 0 \cite{KK92}. Therefore, to determine conditions on $\epsilon$ for which ergodicity fails - or to say the least, for which a non-trivial IUP exists - is also part of the problem in this case. Another interesting feature in this setting is that, despite that the conditioning problems do not explicitly impose that the IUP should also be AsIUP, in each case, the resulting sets all turn out {\sl a posteriori} to be asymmetric, as expected from the numerics.\footnote{Asymmetry can be explicitly included in the conditioning problem, for instance by considering finer symbolic partitions. However, the necessity of such consideration, which generates computational complications, is questionable. Indeed, in tested examples of AsIUP for $G_{2,\epsilon}$ (see Section \ref{S-COUPLEDMAP}), the resulting set turns out to be identical to the original one. No AsIUP additional feature has resulted from such more elaborated consideration.} 

The rest of the paper is organised as follows. Section \ref{S-CONPRO} contains a detailed description of the principles of the approach and of its implementation, together with theoretical considerations about the polytope representation by constraint matrices and their manipulation. Section \ref{S-COUPLEDMAP} presents the results of a full analysis of the conditioning problems associated with coupled maps of low dimension (which include the dynamics of cluster configurations in arbitrary dimension). In order to illustrate the results, the computed AsIUP are plotted against numerical ergodic components and this reveals that sharp approximations of the numerics have been obtained in this way. The analysis itself and related proofs are reported in the Appendices and secondary details of the implementation are given in the Supplementary Material. Finally, some concluding remarks and suggestions of open questions are given in Section \ref{S-CONCL}.

\section{Conditioning problems for IUP: principles and implementation}\label{S-CONPRO}
This section introduces the approach to IUP in this paper, from the basic definitions of the dynamics and of the conditioning problem, to the polytopes' representation and their manipulation in terms of matrices of constraints, and to an algorithmic presentation of the implementation. 

\subsection{Dynamics, empirical information and conditioning problems}
\subsubsection{Definition of the dynamics and related notions}
Let $d\in\N$ and $M\subset \R^d$ be a bounded polytope, typically $M=(0,1)^d$. Let ${\cal A}$ be a finite set ({\bf alphabet}) and consider the following partition of $M$ ({\bf atomic partition})\footnote{Throughout, the notation mod 0 means modulo a set of zero Lebesgue measure.}
\[
M=\bigcup_{\omega\in {\cal A}}A_\omega\ \text{mod}\ 0
\]
where every set $A_\omega$ is an open convex and bounded polytope ({\bf atom}).

Let $F:\bigcup_{\omega\in {\cal A}}A_\omega\to M$ be the mapping defined by $F=a\textrm{Id}+B$ where $a>1$\footnote{The approach is mostly relevant in the case $a\geq 1$ (see in particular footnote 1). However, it equally applies to the case $a\in (0,1)$ where it can be used to prove the existence of stable periodic orbits. For simplicity, focus here will be made on the case $a>1$.} and $B:\bigcup_{\omega\in {\cal A}}A_\omega\to\R^d$ is some piecewise constant function, {\sl viz.}\ $B$ is constant on each $A_\omega$. Let $B_\omega:=B|_{A_\omega}$ be the corresponding constant vector.

A set $U=\bigcup_{k=1}^n P_k$, where the $P_k$ are polytopes, is called an {\bf invariant union of polytopes} (IUP) for $F$ if 
\[
F\left(U\cap \bigcup_{\omega\in {\cal A}} A_\omega \right)\subset U.
\]
If $U$ is an IUP of $F$, then $F$ must have an absolutely continuous invariant measure with support included in $U$ \cite{T01}. 

As mentioned before, the collection $\{A_\omega\}_{\omega\in {\cal A}}$ of all atoms is a trivial IUP. Here, we are interested in IUP whose Lebesgue measure is smaller than that of $M$. A mapping $F$ needs not have such non-trivial IUP (in particular if $F$ is transitive on $M$ \cite{KH95}). However, if it does, hints about such sets might be obtained from numerical simulations of trajectories. 

In the case where $F$ commutes with some $\Z^2$-symmetry $\Sigma$,\footnote{{\sl ie.}\ $\Sigma:M\circlearrowleft$ is an invertible affine transformation such that $\Sigma\neq\text{Id}$ and $\Sigma^2=\text{Id}$. For instance $\Sigma=1-\textrm{Id}$ when $M$ is the unit hypercube.} any IUP $U$ such that $\Sigma (U)\cap U=\emptyset$ is said to be an Asymmetric IUP ({\bf AsIUP}) with respect to $\Sigma$. As mentioned before, if an AsIUP exists, then $F$ must have two Lebesgue ergodic components, at least.

As described in Section \ref{S-214} below, we shall also be interested in the existence of IUP/AsIUP in the case where $F$ has a larger symmetry group than $\{\text{Id},\Sigma\}$.

\subsubsection{Empirical information from numerical orbits}\label{S-212}
Suppose that the simulation of some orbit $\{F^tx\}_{t\in\N}$ of a piecewise affine mapping $F$ as above reveals that the points $F^tx$ aggregate on a finite number of clusters in phase space which they visit perpetually. A set of aggregated points is called a {\bf cluster} if the minimal distance between its points is {\sl microscopic} and the distance between any two clusters is {\sl macroscopic}. Here microscopic means small and getting arbitrarily small as the number of orbit points increases. Macroscopic means bounded below by a positive number, independently of the number of orbit points, see Appendix \ref{A-SYSTEMATIC} for a quantitative criterion.

Assume that these features of the orbit under consideration convey the numerical evidence of a genuine invariant set with positive Lebesgue measure. Then, the standard construction of symbolic dynamics in piecewise expanding systems that possess a Markov partition \cite{KH95} suggests that this set might be characterized based upon certain topological features, in particular, its structure, localisation and dynamics.
Accordingly, the following details are to be extracted from the observation under consideration (Fig.\ \ref{ILLUSTR}, left panel):
\begin{figure}[ht]
\begin{center}
\includegraphics*[width=60mm]{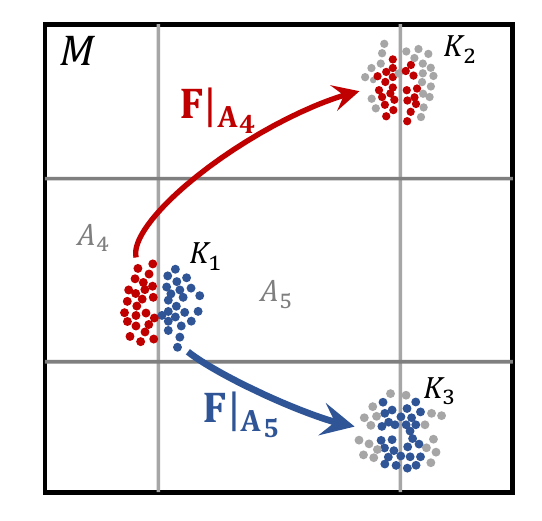}
\hspace{10mm}
\includegraphics*[width=60mm]{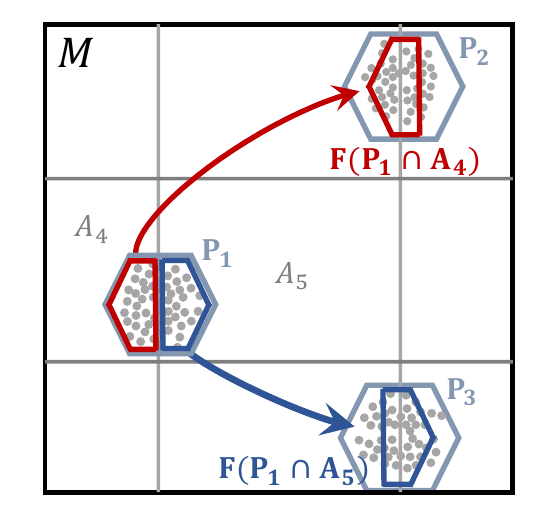}
\end{center}
\caption{Illustration of the approach. {\sl Left.} Assuming the observation that points $F^tx$ of an orbit $\{F^tx\}_{t\in\N}$ aggregate in a collection of clusters that they visit perpetually, the approach starts by collecting the following information: the number of clusters, their localisation in the atomic partition ({\sl ie.}\ $K_1\cap A_\omega\neq \emptyset$ iff $\omega\in\{4,5\}$ and likewise for $K_2$ and $K_3$) and the dynamical transitions between clusters (namely $F(K_1\cap A_4)\subset K_2$ and  $F(K_1\cap A_{5})\subset K_3$ and likewise for $K_2$ and $K_3$). {\sl Right.}\ This empirical information suggests to find an IUP $\bigcup_{k=1}^3P_k$ as a solution of a conditioning problem with similar dynamical characteristics.}  
\label{ILLUSTR}
\end{figure}
\begin{itemize}
\item the number of clusters, say $q$. Let then $\{K_k\}_{k=1}^q$ be an enumeration of the clusters,
\item their localisation in the atomic partition, {\sl ie.}\ for $k\in\{1,\cdots ,q\}$, let ${\cal A}_k\subset {\cal A}$ be such that $K_k\cap A_\omega\neq \emptyset$ iff $\omega\in {\cal A}_k$,
\item the transition between clusters, {\sl ie.}\ for each $K_k\cap A_\omega\neq \emptyset$, let $\ell(k,\omega)\in\{1,\cdots ,q\}$ be such that $F(K_k\cap A_\omega)\subset K_{\ell(k,\omega)}$.\footnote{That distinct clusters are separated by a macroscopic distance implies that $\ell(k,\omega)$ must be unique, {\sl ie.}\ no image $F(K_k\cap A_\omega)$ can intersect two clusters.}
\end{itemize}
To collect these details can be made largely automatized, or at least be made systematic, see Appendix \ref{A-SYSTEMATIC} for related instructions. Such an automatization is particularly useful when the dimension $d$ is large and direct visualisation in phase space is not practicable.

\subsubsection{Conditioning problems for IUP}
In order to mathematically confirm that an observation as above materializes an authentic invariant set of $F$, it suffices to prove the existence of an IUP that possesses the same characteristics as the observation \cite{T01}. This task can be formulated as a so-called conditioning problem, namely the question of determining a collection of polytopes that satisfy the topological conditions extracted from the observation (Fig.\ \ref{ILLUSTR}, right panel). In formal terms, a conditioning problem can be expressed as follows:

\noindent
{\bf Conditioning problem for an IUP.} {\sl Assume that a number $q\in \N$ of clusters is given together with a specification $\{{\cal A}_k\}_{k=1}^q$ (${\cal A}_k\subset {\cal A})$ of their localisation in the atomic partition and a specification $\{\{\ell(k,\omega)\}_{\omega\in {\cal A}_k}\}_{k=1}^q$ ($\ell(k,\omega)\in \{1,\cdots ,q\}$) of the inter-clusters transitions. Find a collection $\{P_k\}_{k=1}^q$ of polytopes that satisfy the following series of conditions:
\begin{itemize} 
\item For each $k\in\{1,\cdots ,q\}$, we have $P_k\cap A_\omega\neq \emptyset$ iff $\omega\in {\cal A}_k$.
\item For each $k\in\{1,\cdots ,q\}$ and $\omega\in {\cal A}_k$, we have $F(P_k\cap A_\omega)\subset P_{\ell(k,\omega)}$.
\end{itemize}}
\noindent
The following comments are in order:
\begin{itemize}
\item[] In practice, the conditions above may be imposed only for a sub-collection of polytopes, once such a restriction suffices to ensure that the full union set $U=\bigcup_{k=1}^qP_k$ is an IUP of $F$. This comment is particularly relevant in presence of symmetries, see Section \ref{S-214} below.
\item[] Some IUP might include dynamical iterates of some of their constituents, {\sl ie.}\ some $P_{k'}$ might be the image under $F$ of some atomic restriction of some $P_k$. In this case, we shall impose the stronger condition $F(P_k\cap A_\omega)= P_{\ell(k,\omega)}$ for some pairs $(k,\omega)$ of indices (see Conditioning Problem \ref{CONDPRO2} and \ref{CONDPRO3} in Section \ref{S-COUPLEDMAP} for examples in the coupled maps).
\item[] Unless we have $F(P_k\cap A_\omega)= P_{\ell(k,\omega)}$ for all admissible pairs $(k,\omega)$, the non-wandering set $\Omega$ of the restriction $F|_{U}$ is expected to have Lebesgue measure smaller than that of $U$, and $F|_\Omega$  need not be conjugated to a topological Markov chain. In this case, knowledge of the IUP $U$ is not enough to conclude about ergodic properties of the restriction $F|_\Omega$. However, most importantly for our purpose and as said before, every IUP must support an absolutely continuous invariant measure and the existence of two distinct IUP implies that the full map $F$ must have, at least, two Lebesgue ergodic components in $M$.
\end{itemize}

\subsubsection{Symmetric IUP and related reduced conditioning}\label{S-214}
Assume that $F$ has a {\bf symmetry group}, namely that there exists a group of invertible affine transformations $\{\sigma_i\}_{i=1}^s$ (called {\bf symmetries}) where each $\sigma_i:M\circlearrowleft$ commutes with $F$. Then, we may consider IUP that are composed by orbits under this group, {\sl ie.}\
\[
U=\bigcup_{i=1}^s\bigcup_{k=1}^{q'}\sigma_i(P_k).
\]
Of note, if $F$ also commutes with some $\Z_2$-symmetry $\Sigma$ and $\sigma_i\neq \Sigma$ for all $i\in \{1,\cdots ,s\}$, then such IUP $U$ might be a candidate for an AsIUP with respect to $\Sigma$.

In order to determine $U$, it obviously suffices to solve the conditioning problem for the polytopes $P_k$. The conditions on the remaining polytopes $\sigma_i(P_k)$ (when $\sigma_i\neq \text{Id}$) will automatically follow from the commutation property. This simplification ({\bf reduced conditioning}) will largely be employed in the coupled map examples below, where the symmetry groups are given by subgroups of the $N$-element permutation group (or more precisely, their representation in the space under consideration) and where $\Sigma=-\text{Id}$.

In order to specify a conditioning problem in presence of symmetries, one also needs to collect from the numerics, those labels $I(k,\omega)\in \{1,\cdots ,s\}$ of the symmetries involved in the dynamical conditions, {\sl viz.}\ such that $F(P_k\cap A_\omega)\subset \sigma_{I(k,\omega)}(P_{\ell(k,\omega)})$. In summary, a conditioning problem in presence of symmetries can be formulated as follows.

 \noindent
{\bf Conditioning problem for an IUP, in presence of symmetries.} {\sl In presence of a symmetry group $\{\sigma_i\}_{i=1}^s$, assume that a number\footnote{From now on, we use the symbol $q$ instead of $q'$, for simplicity.} $q\in \N$ of clusters is given together with a specification $\{{\cal A}_k\}_{k=1}^q$ of their localisation in the atomic partition, a specification $\{\{\ell(k,\omega)\}_{\omega\in {\cal A}_k}\}_{k=1}^q$ of the inter-clusters transitions and a specification $\{\{I(k,\omega)\}_{\omega\in {\cal A}_k}\}_{k=1}^q$ ($I(k,\omega)\in \{1,\cdots ,s\}$) of the transition-involved symmetries. Find a collection $\{P_k\}_{k=1}^{q}$ of polytopes that satisfy the following series of conditions:
\begin{itemize} 
\item For each $k\in\{1,\cdots ,q\}$, we have $P_k\cap A_\omega\neq \emptyset$ iff $\omega\in {\cal A}_k$.
\item For each $k\in\{1,\cdots ,q\}$ and $\omega\in {\cal A}_k$, we have $F(P_k\cap A_\omega)\subset \sigma_{I(k,\omega)}(P_{\ell(k,\omega)})$.
\end{itemize}}
\noindent
Naturally, this formulation of conditioning problem is an extension of the previous one, which can be recovered in the case where $\sigma_{I(k,\omega)}=\text{Id}$ for all $(k,\omega)$.

\subsection{Analytic representation of polytopes and their manipulation}
\subsubsection{Analytic representation of polytopes: definition and basic considerations on coefficient matrices}\label{S-221}
The analysis of conditioning problems requires dealing with arbitrary polytopes in arbitrary dimension. Since every (open) polytope in $\R^d$ can be regarded as the intersection of (open) half-spaces, see for instance \cite{G03}, polytopes can be characterized using interval constraints on linear combinations of the coordinates $\{x_i\}_{i=1}^d$. To that goal, we shall use and complete the formalism in \cite{F20} which provides a representation of polytopes in terms of matrices that collect the constraints. More precisely, \cite{F20} introduced the formalism of constraint matrices, their basic manipulation and the optimization procedure of their coordinates.\footnote{For completeness, we mention that the algorithm of the computer-assisted proof in \cite{F20} is based on an iterative procedure. It begins with an asymmetric cylinder of the atomic partition and then computes its iterates under the dynamics, until either the constructed union sets intersects its symmetric image - in which case the construction has to be restarted with another initial cylinder - or an invariant set results, which then must be an AsIUP.} Here, the topological operations and geometric manipulations are improved and complemented, in particular in order to include general linear transformations, see Section \ref{S-223}.

In a more formal way, given any open polytope $P\subset\R^d$, let a non-degenerate\footnote{A coefficient matrix is said to be {\bf non-degenerate} if all the vectors $\alpha_i=(\alpha_{ij})_{j=1}^d$ are distinct vectors in the real projective space $\R\P^d$.} {\bf coefficient matrix} $\alpha=\left(\alpha_{i1}\ \alpha_{i2}\ \cdots\ \alpha_{id}\right)_{i=1}^e\in \R^{e\times d}$
where $e\geq d$, and a {\bf constraint matrix} $m=(\underline{m}_{i}\ \overline{m}_{i})_{i=1}^{e}\in\R^{e\times 2}$ be so that 
\begin{equation}
P=P^\alpha_m:=\left\{x\in \R^d\ :\ \underline{m}_{i}<(\alpha x)_i<\overline{m}_{i}\ \text{for all}\ 1\leq i\leq e\right\}.
\label{REP-POLY}
\end{equation}
The matrix $m$ can also be regarded as the pair of vectors $m=(\underline{m},\overline{m})$ where $\underline{m}=(\underline{m}_i)_{i=1}^e$ and $\overline{m}=(\overline{m}_i)_{i=1}^e$.

A priori, both $\alpha$ and $m$ vary with the polytopes under consideration. However, given the nature of the conditioning problems, we shall be concerned by a single coefficient matrix $\alpha$ for the whole collection of polytopes. Only the constraint matrices $m$ will be specific to the elements of the IUP. 

A natural candidate for $\alpha$ is the {\bf canonical coefficient matrix} associated with the symbolic partition, {\sl viz.}\ every polytope face is parallel to some atom face.\footnote{In particular, all atoms, their images and all pre-images, can be captured by this canonical matrix, see Appendix \ref{A-PARTITION} for related notations of the constraint matrix entries.} For simplicity, we may drop the superscript $\alpha$ in $P^\alpha_m$ and use the notation $P_m$ in this case.

However, as some coupled map examples show, some conditioning problems turn out to have no solution in the canonical setting, {\sl ie.}\ such conditioned IUP cannot exist whose faces of its $P_k$ are (all) parallel to atom facets. 
More/different rows have to be added to $\alpha$ in order to include more directions and to obtain admissible solutions. We shall see in the examples that suitable coefficient matrices may be obtained (at the expense of cumbersome calculations) as part of the solution of the conditioning problem. In this case, this means that the information about face orientation is also implicitly embedded in the transcribed empirical knowledge.

Notice also that, given a coefficient matrix, multiple solutions may exist for the constraint matrices associated with IUP elements, as it is the case in some of the coupled map examples.

\subsubsection{Analytic representation of polytopes: optimization procedure}\label{S-222}
When the coefficient matrix has more rows than columns ({\sl ie.}\ $e>d$), for some polytopes $P$, some constraints in \eqref{REP-POLY} need not be active. An inequality in \eqref{REP-POLY} is said to be an {\bf active constraint} \cite{BV04} if there exists a point in the closure $\overline{P}$ for which it becomes an equality. Conversely, a constraint is said to be {\bf not active} or inactive, if it cannot be saturated by a point in  $\overline{P}$. 

When some constraints are not active, matrices $m$ such that $P_m^\alpha=P$ are not unique, because modifying (slightly) the entries associated with inactive constraints does not alter the polytope. Inactive constraints are problematic in the conditioning problem because they may yield stronger than necessary inequalities on some of the entries of the constraint matrices, which may in turn prevent one to obtain solutions.

Furthermore, to obtain necessary and sufficient conditions for the existence of solutions in the formalism of equation \eqref{REP-POLY} requires to be able to assert that a polytope $P_m^\alpha$ is not empty.\footnote{Naturally, the condition $\underline{m}_i< \overline{m}_i$ for all $i\in\{1,\cdots ,e\}$ is necessary for having $P^\alpha_m\neq\emptyset$.} This can be done by introducing an optimization scheme, namely a transformation $O:\R^{e\times 2}\circlearrowleft$ on constraint matrices, that aims to make all constraints active \cite{F20}.

Let $e>d$ and a coefficient matrix $\alpha$ be given. For each $i\in\{1,\cdots ,e\}$, consider the set $\Lambda_{i}$ of vectors $\lambda=(\lambda_{k})_{k=1}^{e}\in\R^e$ with at least $e-d$ vanishing entries,
which uniquely solve the system of equations
\begin{equation}
(\lambda^T\alpha)_{j}=\alpha_{ij},\ 1\leq j\leq d.
\label{LAGRANGE}
\end{equation}
More precisely, given any $s\in \{1,\cdots,d\}$ and $S\subset\{1,\cdots,e\}$ of cardinality $s$, if it exists, let $\lambda$ be the unique solution of the equations obtained from \eqref{LAGRANGE} by letting $\lambda_k=0$ for $k\in \{1,\cdots,e\}\setminus S$. The set $\Lambda_i$ is made of all such solutions $\lambda$ when $S$ and $s$ vary.\footnote{Each $\Lambda_{i}$ must be a non-empty finite set. In particular it contains the canonical vector $(\lambda_k)_{k=1}^{e}=(\delta_{k,i})_{k=1}^{e}$, where $\delta_{k,i}$ is the Kronecker symbol.}

Independently, given $a\in\R$ and a constraint matrix $m$, consider the vectors $\underline{e}(a,m):=(\underline{e}_k(a,m))_{k=1}^{e}$ and $\overline{e}(a,m):=(\overline{e}_k(a,m))_{k=1}^{e}$ defined by 
\[
\left\{\begin{array}{ccccl}
\underline{e}(a,m)=\underline{m}&\text{and}&\overline{e}(a,m)=\overline{m}&\text{if}&a\geq 0\\
\underline{e}(a,m)=\overline{m}&\text{and}&\overline{e}(a,m)=\underline{m}&\text{if}&a< 0
\end{array}\right.
\]
The {\bf optimized constraint matrix} $O(m)=(\underline{O(m)},\overline{O(m)})$ is defined by
\begin{equation*}
\underline{O(m)}_i=\max_{(\lambda_k)\in\Lambda_i}\sum_{k=1}^{e}\lambda_k\underline{e}_k(\lambda_k,m)\quad\text{and}\quad
\overline{O(m)}_i=\min_{(\lambda_k)\in\Lambda_i}\sum_{k=1}^{e}\lambda_k\overline{e}_k(\lambda_k,m),\ \forall i\in\{1,\cdots ,e\},
\end{equation*}
As announced, a crucial property of the optimization procedure is that it determines the existence of the corresponding polytope and also makes sure that all constraints are active. These properties are summarized in the following statement.
\begin{Lem}
{\rm \cite{F20}} (i) Given any constraint matrix $m$, the polytope $P_m^\alpha$ defined by \eqref{REP-POLY} is not empty iff $\underline{O(m)}_i<\overline{O(m)}_i$ for all $i\in\{1,\cdots ,e\}$.

\noindent
(ii) If $P_m^\alpha\neq \emptyset$, then $P_{O(m)}^\alpha=P_m^\alpha$ and all constraints in the definition of $P_{O(m)}^\alpha$ are active.

\noindent
(iii) The direction/plane $\sum_{j=1}^d\alpha_{ij}x_j=\underline{O(m)}_i$ (resp.\ $\sum_{j=1}^d\alpha_{ij}x_j=\overline{O(m)}_i$) defines an edge/a face of $P_m^\alpha$ iff 
\[
\max_{(\lambda_k)\in\Lambda_i\atop (\lambda_k)\neq (\delta_{k,i})}\sum_{k=1}^{e}\lambda_k\underline{e}_k(\lambda_k,O(m))<\underline{O(m)}_i\ 
\left(\text{resp.}\ \overline{O(m)}_i<\min_{(\lambda_k)\in\Lambda_i\atop (\lambda_k)\neq (\delta_{k,i})}\sum_{k=1}^{e}\lambda_k\overline{e}_k(\lambda_k,m)\right).
\]
\label{OPTIVEC}
\end{Lem}
\noindent
{\sl Remark:} In the coupled maps example, the optimization procedure can be used to determine the atoms of the partition, which are not know {\sl a priori}, see Appendix \ref{A-PARTITION}. Moreover, statement {\sl (iii)} of the Lemma will be employed in order to identify the faces/edges of the constructed polytopes.

\subsubsection{Analytic formulation of topological operations and geometric transformations on polytopes}\label{S-223}
To address conditioning problems also requires to manipulate and to execute certain topological operations on polytopes, especially to define the intersection between two such sets, to compute the image under affine transformations and to verify inclusion in a given polytope. The formalism of equation \eqref{REP-POLY} allows one to implement these operations on constraint matrices, as presented in this section. Below, we consider those operations that will be employed in the sequel, assuming that the non-degenerate coefficient matrix $\alpha$ is given, unless otherwise specifically mentioned.

The proofs of the claims in this section are all elementary and mostly left to the reader. 
\medskip

\noindent
$\bullet$ {\bf Intersection.} Given two constraint matrices $m$ and $m'$, the intersection constraint matrix $m\cap m'$ is defined by the following vectors
\[
\underline{(m\cap m')}:=\left(\max\{\underline{m}_i,\underline{m'}_i\}\right)_{i=1}^e\quad\text{and}\quad \overline{(m\cap m')}:=\left(\min\{\overline{m}_i,\overline{m'}_i\}\right)_{i=1}^e.
\]
An immediate consequence of statement {\sl (ii)} in Lemma \ref{OPTIVEC} is the following elementary characterization of the intersection of polytopes.
\begin{Claim}
$P_m^\alpha\cap P_{m'}^\alpha\neq\emptyset$ iff $P_{O(m\cap m')}^\alpha\neq\emptyset$. Moreover, if not empty, we have $P_m^\alpha\cap P_{m'}^\alpha=P_{O(m\cap m')}^\alpha$.
\end{Claim}
\noindent
Of note, the following observation is convenient in order to ensure empty intersection
\[
\text{If}\ \underline{(m\cap m')}_i\geq \overline{(m\cap m')}_i\ \text{for some $i\in \{1,\cdots,e\}$, then}\ P_m^\alpha\cap P_{m'}^\alpha=\emptyset.
\]
\medskip

\noindent
$\bullet$ {\bf Inclusion.} We say that $m\subset m'$ holds for the constraint matrices $m$ and $m'$ if $\underline{m'}_i\leq \underline{m}_i$ and $\overline{m}_i\leq \overline{m'}_i$ for all $i\in \{1,\cdots,e\}$.  Similarly as for the intersection, we have
\begin{Cor}
(i) $P^\alpha_{m}\subset P^\alpha_{m'}$ iff $O(m)\subset m'$.

\noindent
(ii) $P^\alpha_{m}= P^\alpha_{m'}$ iff $O(m)= O(m')$. 
\label{IFFINCLU}
\end{Cor}

\noindent
{\sl Proof:} {\sl (i)} $O(m)\subset m'$ implies $P^\alpha_{O(m)}=P^\alpha_{m}\subset P^\alpha_{m'}$. Conversely, assume that $O(m)\not\subset m'$. Then statement {\sl (ii)} in the Lemma implies that there is a point in $P^\alpha_m$ that does not satisfy all constraints that define $P^\alpha_{m'}$.

\noindent
{\sl (ii)} $O(m)\subset O(m')$ implies $P^\alpha_{m}=P^\alpha_{O(m)}\subset P^\alpha_{O(m')}=P^\alpha_{m'}$. Then, a similar argument obtained by exchanging the roles of $m$ and $m'$ easily yields the desired statement. \hfill $\Box$
\medskip

\noindent
$\bullet$ {\bf Basic affine transformations.}

\noindent
$\ast$ {\sl Translations.} Given a constraint matrix $m$ and $x\in\R^d$, let $m+\alpha x$ be the constraint matrix defined by $m+\alpha x :=(\underline{m}+\alpha x,\overline{m}+\alpha x)$. 
\begin{Claim}
$O(m+\alpha x)=O(m)+\alpha x$ and $P_m^\alpha+x=P_{O(m)+\alpha x}^\alpha$.
\end{Claim}

\noindent
$\ast$ {\sl Isotropic scalings.} Given a constraint matrix $m$ and $a>0$, let $am$ be the constraint matrix defined by $am:=(a\underline{m},a\overline{m})$. 
\begin{Claim}
$O(am)=aO(m)$ and $a\text{\rm Id}(P_m^\alpha)=P_{aO(m)}^\alpha$.
\end{Claim}
\noindent
As a consequence, we have for the restriction to atomic pieces, of a piecewise affine mapping $a\text{Id}+B$ as above 
\[
F|_{A_\omega}(P^\alpha_m)=aP_m^\alpha+B_{\omega}=P^\alpha_{aO(m)+\alpha B_{\omega}}.
\]

\noindent
$\ast$ {\sl Sign inversion.} Given a constraint matrix $m$, let $m'$ be the constraint matrix defined by $m':=(-\overline{m},-\underline{m})$.
\begin{Claim}
$O(m')=(-\overline{O(m)},-\underline{O(m)})$ and $-\text{\rm Id}(P_m^\alpha)=P_{O(m')}^\alpha$.
\end{Claim}
\medskip

\noindent
$\bullet$ {\bf More general affine transformations.} For conditioning problems that involve some symmetries, we also need to make sure that the action of such transformation on polytopes can be implemented on constraint matrices. To that goal, we introduce the following notion.

Given a non-degenerate coefficient matrix $\alpha$, an invertible affine transformation $\sigma$ of $\R^d$ is said to be {\bf $\alpha$-compatible} if the following relation holds 
\begin{equation}
\alpha\sigma x=A\pi\alpha x+B,\ \forall x\in\R^d
\label{SPETRANS}
\end{equation}
where $A=\text{diag}(a_i)$ is a diagonal matrix with $(a_i)\in\R^e\setminus \{0\}$, $\pi\in\Pi_e$ is the (representation in $\R^e$ of the) group of permutations of $\{1,\cdots ,e\}$ and $B\in\R^e$. If $\sigma$ is $\alpha$-compatible, then it induces a transformation $m\mapsto \sigma(m)$ on constraint matrices (NB: once again, we use the same symbol for simplicity) defined by 
\[
\left\{\begin{array}{l}
\underline{\sigma(m)}_i=a_i\pi\left(\underline{m}\right)_i+B_i\quad \text{and}\quad \overline{\sigma(m)}_i=a_i\pi\left(\overline{m}\right)_i+B_i\ \text{if}\ a_i>0,\\
\underline{\sigma(m)}_i=a_i\pi\left(\overline{m}\right)_i+B_i\quad \text{and}\quad \overline{\sigma(m)}_i=a_i\pi\left(\underline{m}\right)_i+B_i\ \text{if}\ a_i<0.
\end{array}\right.
\]
As before, the following statement readily follows from the definitions.
\begin{Claim}
Let a coefficient matrix $\alpha$ and an invertible affine transformation $\sigma$ be given.

\noindent
{\sl (i)} If $\sigma$ is $\alpha$-compatible, then we have $\sigma(P^\alpha_m)= P^\alpha_{\sigma(m)}$.

\noindent
{\sl (ii)} If $\sigma$ is $\alpha$-compatible and the diagonal matrix $A$ in \eqref{SPETRANS} writes $A=a\text{\rm Id}|_{\R^e}$ for some $a\in\R\setminus\{0\}$, then for every coefficient matrix $m$, we have $O(\sigma(m))=\sigma(O(m))$.
\label{CLAIMSPETRANS}
\end{Claim}
Any affine transformation of the form $a\text{Id}|_{\R^d}+B_\omega$ ($a\neq 0$), and in particular the inversion symmetry $\Sigma=1-\text{Id}$ trivially satisfies all assumptions of this statement, for any coefficient matrix $\alpha$. However, and as mentioned above, a real concern is to make sure that the symmetries involved in a conditioning problem are $\alpha$-compatible, for a suitable matrix $\alpha$.\footnote{Otherwise, when we only have $\sigma(P^\alpha_m)\subset P^\alpha_{\sigma(m)}$, this condition does not suffice to ensure that $m\subset \sigma(m')$ implies $P_m\subset \sigma(P_{m'})$.}  This issue is addressed by the following statement.
\begin{Claim}
Let $\{\sigma_i\}_{i=1}^s$ be a collection of invertible affine transformations of $\R^d$ which generate a finite group. Then there exist $e\geq d$ and an $e\times d$ non-degenerate coefficient matrix $\alpha$ such that every transformation in the group is $\alpha$-compatible.
\label{ADAPTMAT}
\end{Claim}
In particular, any collection of symmetry transformations in the coupled map examples below satisfies the assumptions of this claim.

\noindent
{\sl Proof.} Notice first that if a transformation $\sigma$ satisfies \eqref{SPETRANS} for some diagonal matrix $A$ and vector $B$, then the transformation $x\mapsto \sigma x+B'$ satisfies the same relation with the same $A$ and the vector $B+\alpha B'$. Hence, we may assume w.l.o.g.\ that the $\sigma_i$ are linear transformations of $\R^d$. 

Let $\{\sigma_i\}_{i=1}^g$ be an enumeration of the group of transformations, {\sl ie.}\ for $i>s$, $\sigma_i$ is a product of generators and assume that an $e\times d$ non-degenerate matrix $\alpha_0$ is given for some $e\geq d$ (for instance $e=d$ and $\alpha_0=\text{Id}$). We claim that the $(e\cdot g)\times d$ matrix $\alpha$ defined by 
\[
\alpha=\left(\begin{array}{c}
\alpha_0\sigma_1\\
\hline
\alpha_0\sigma_2\\
\hline
\vdots\\
\hline
\alpha_0\sigma_g
\end{array}\right)
\]
is as desired, up to non-degeneracy. Indeed, that $\{\sigma_i\}_{i=1}^g$ is a group implies that, for every $i\in\{1,\cdots ,s\}$, there exists a permutation $p_i$ of $\{1,\cdots, g\}$ such that 
\[
\sigma_j\circ \sigma_i=\sigma_{{p_i}_j},\ \forall j\in \{1,\cdots ,g\}.
\]
It easy to conclude that $\sigma_i$ satisfies \eqref{SPETRANS} with $A=\text{Id}|_{\R^{e\cdot g}}$, $B=0$, and $\pi$ the permutation in $\{1,\cdots , e\cdot g\}$ (or rather its representation in $\R^{e\cdot g}$) defined by the permutation of $e$-blocks of indices induced by $p_i$. 

It remains to prove that $\alpha$ can be chosen to be non-degenerate. Assume it is not, {\sl viz.}\ there are two rows that are multiple of each other, $\alpha_{\iota'}=c\alpha_\iota$ for some pair $\iota,\iota'$ and $c\in\R\setminus\{0\}$. In this case, the row $\iota'$ should be removed form $\alpha$. Moreover, equation \eqref{SPETRANS} implies that we must also have $\alpha_{\pi_{\iota'}}=c\alpha_{\pi_\iota}$ for the permutation $\pi$ associated with transformation $\sigma_i$; hence the row $\pi_{\iota'}$ should also be removed. Likewise, the rows $\pi^k_{\iota'}$ for every $k$ and every such permutation $\pi$ should be removed from $\alpha$, so that the matrix becomes non-degenerate. That $\pi=\text{Id}$ for $\sigma_i=\text{Id}$ and the assumption that $\alpha_0$ is non-degenerate imply that the resulting matrix cannot non-empty (and contains $\alpha_0$). \hfill $\Box$

\subsection{Analytic formulation of conditioning problems and implementation of the systematic approach to IUP construction}\label{S-ANALYSIS}
The results of the previous sections can be grouped into the following statements, which provide a formulation of the conditioning problems in terms of constraint matrices.  

Given a coefficient matrix $\alpha$, let $\{m_\omega\}_{\omega\in {\cal A}}$ be the collection of optimized constraint vectors associated with the atomic collection $\{A_\omega\}_{\omega\in {\cal A}}$, {\sl viz.}\ $m_\omega=O(m_\omega)$ and $A_\omega=P_{m_\omega}^\alpha$ for all $\omega\in {\cal A}$.
\begin{Pro}
{\bf (i) Analytic formulation of a standard conditioning problem.} Given a coefficient matrix $\alpha$, the collection of polytopes $\{P_{m_k}^\alpha\}_{k=1}^{q}$ associated with the constraint matrices $\{m_k\}_{k=1}^{q}$ satisfies the conditioning problem defined by the integer $q$, the localisation labels $\{{\cal A}_k\}_{k=1}^q$ and the transition indices $\{\{\ell(k,\omega)\}_{\omega\in {\cal A}_k}\}_{k=1}^q$ iff 
\begin{itemize} 
\item For each $k\in\{1,\cdots ,q\}$, we have $P_{O(m_k\cap m_\omega)}^\alpha\neq \emptyset$ iff $\omega\in {\cal A}_k$. 
\item For each $k\in\{1,\cdots ,q\}$ and $\omega\in {\cal A}_k$, we have $aO(m_k\cap m_\omega)+\alpha B_\omega\subset m_{\ell(k,\omega)}$ (or $aO(m_k\cap m_\omega)+\alpha B_\omega= O(m_{\ell(k,\omega)})$ if the original topological condition is an equality).
\end{itemize}

\noindent
{\bf (ii) Analytic formulation of a conditioning problem with symmetries.} In presence of a symmetry group $\{\sigma_i\}_{i=1}^s$, assume that a coefficient matrix $\alpha$ is given such that every transformation $\sigma_i$ is $\alpha$-compatible. Then, the collection of polytopes $\{P_{m_k}^\alpha\}_{k=1}^{q}$ associated with the constraint matrices $\{m_k\}_{k=1}^{q}$ satisfies the conditioning problem defined by the integer $q$, the localisation labels $\{{\cal A}_k\}_{k=1}^q$, the transition indices $\{\{\ell(k,\omega)\}_{\omega\in {\cal A}_k}\}_{k=1}^q$ and the symmetry indices $\{I(k,\omega)\}$ iff
\begin{itemize} 
\item For each $k\in\{1,\cdots ,q\}$, we have $P_{O(m_k\cap m_\omega)}^\alpha\neq \emptyset$ iff $\omega\in {\cal A}_k$. 
\item For each $k\in\{1,\cdots ,q\}$ and $\omega\in {\cal A}_k$, we have $aO(m_k\cap m_\omega)+\alpha B_\omega\subset \sigma_{I(k,\omega)} (m_{\ell(k,\omega)})$ (or $aO(m_k\cap m_\omega)+\alpha B_\omega= (O\circ \sigma_{I(k,\omega)})(m_{\ell(k,\omega)})$ if the original topological condition is an equality).
\end{itemize}
\label{SUMMARY}
\end{Pro}
This proposition is the core element of the following computational approach to the construction of IUP. 
In particular, the analysis of coupled map examples in the next section closely follows this sequence of instructions (in the case of a conditioning problem with symmetries).
\medskip

\noindent
{\bf Implementation of systematic approach to IUP construction.} For the sake of space, the following algorithm has been designed to simultaneously accommodate both standard conditioning problems and the conditioning problems in presence of symmetries. In particular, in order to apply it to the first case, it suffices to ignore the instructions related to the presence of symmetries.
\begin{itemize}
\item {\sl Preliminaries.} Given a polytope $M\subset \R^d$ and a piecewise affine mapping $F=a\textrm{Id}+B$ from $M$ into itself, suppose that the numerical simulation of some orbit provides evidence of a non-trivial invariant set with positive Lebesgue measure.
\item From the observation, get the topological information about the structure, localization and dynamics of the empirical orbit, {\sl viz.}\ the number of clusters $q$, the cluster localisation labels $\{{\cal A}_k\}_{k=1}^q$ and the inter-cluster transition indices $\{\{\ell(k,\omega)\}_{\omega\in {\cal A}_k}\}_{k=1}^q$ (and also, in presence of a symmetry group, the symmetry indices $\{I(k,\omega)\}$ ). See Section \ref{S-212} for details and definitions and Appendix \ref{A-SYSTEMATIC} for the corresponding systematic procedure.
\item  Consider some coefficient matrix $\alpha$ (Section \ref{S-221}), for instance the canonical matrix associated with the symbolic partition of $F$ or, ideally, an arbitrary coefficient matrix. Then, 
\begin{itemize}
\item[$\ast$] {\sl In presence of symmetries.} Check that all the symmetries $\{\sigma_i\}_{i=1}^s$ are $\alpha$-compatible (Section \ref{S-223}).
\item[$\ast$] Compute the solutions $\Lambda_i$ of \eqref{LAGRANGE} in order to specify the corresponding optimization function $O$ (Section \ref{S-222}). Then, compute the optimized constraints matrices $m_\omega=O(m_\omega)$ associated with the atoms $A_\omega$ for $\omega\in {\cal A}_k$ and $k\in \{1,\cdots, q\}$. 
\end{itemize}
\item Consider the collection $\{m_k\}_{k=1}^q$ of optimized constraint matrices ({\sl ie.}\ $m_k=O(m_k)$ for all $k$) which are assumed to satisfy the conditions of Proposition \ref{SUMMARY}, either {\sl (i)}, or {\sl (ii)} in the presence of symmetries. (NB: The $\{m_k\}_{k=1}^q$ are the unknown of the conditioning problem under consideration.)
\item Use the conditions $P_{O(m_k\cap m_\omega)}^\alpha\neq \emptyset$ iff $\omega\in {\cal A}_k$ and $P_{m_k}^\alpha\subset M$ to obtain restrictions on the range of the entries of $m_k$. Use these restrictions in order to simplify the formal expression of the constraint matrices $m_k\cap m_\omega$, and subsequently to simplify the expression of $O(m_k\cap m_\omega)$. 
\item Use these simplifications to solve the following piecewise linear inequalities/equations for the matrices $\{m_k\}_{k=1}^q$
\begin{itemize}
\item[$\ast$]  either $aO(m_k\cap m_\omega)+\alpha B_\omega\subset m_{\ell(k,\omega)}$/$aO(m_k\cap m_\omega)+\alpha B_\omega= m_{\ell(k,\omega)}$, in absence of symmetries,
\item[$\ast$]  or $aO(m_k\cap m_\omega)+\alpha B_\omega\subset \sigma_{I(k,\omega)} (m_{\ell(k,\omega)})$/$aO(m_k\cap m_\omega)+\alpha B_\omega= (O\circ \sigma_{I(k,\omega)})(m_{\ell(k,\omega)})$\footnote{If the symmetry $\sigma_{I(k,\omega)}$ satisfies the conditions in statement {\sl (ii)} of Claim \ref{CLAIMSPETRANS}, then this set of equations simplifies as $aO(m_k\cap m_\omega)+\alpha B_\omega= \sigma_{I(k,\omega)}(m_{\ell(k,\omega)})$} in presence of symmetries.
\end{itemize}
In the case where no solution exists for the current matrix $\alpha$, reconsider the problem with a larger, and to be determined, matrix $\alpha$.
\item Assuming solutions exist, check that all matrices $m_k$ are indeed optimized and that all conditions $P_{O(m_k\cap m_\omega)}^\alpha\neq \emptyset$ iff $\omega\in {\cal A}_k$ hold. 
\item In case of AsIUP, check that $\Sigma(U)\cap U=\emptyset$ holds for the IUP generated by the collection $\{P_{m_k}^\alpha\}_{k=1}^{q}$. 
\end{itemize}
\medskip

\noindent
{\sl Remarks:} (i) In low dimension, to complete the computations associated with this procedure remain accessible and one can determine all solutions (which sometimes are not unique). However, the computations become heavy when the number of rows of $\alpha$ ({\sl ie.} admissible number of faces) increases. Computation may require the use of formal tools for algebraic manipulation. In particular, computations related to the AsIUP that emerges at the second bifurcation of the three-dimensional map $G_{3,\epsilon}$ involved ten possible faces ({\sl viz.}\ the corresponding constraint matrices $m\in\R^{10\times 2}$) are executed using the software Mathematica (see Supplementary Material). 

(ii) For IUP/AsIUP that are given by other means ({\sl e.g.}\ {\sl ad-hoc} solutions), the same series of operations can be applied in order to prove that the proposed sets are indeed solutions of the conditioning problem. In particular, this is the adopted approach for the  AsIUP that emerges for $G_{3,\epsilon}$ at the second bifurcation. 

\section{Application to globally coupled maps and their cluster dynamics}\label{S-COUPLEDMAP}
This section presents results of the application of the approach above to the construction of AsIUP in expanding piecewise affine systems of coupled maps. We start by giving the definition of the dynamics for an arbitrary distribution and its connection to the cluster dynamics in basic globally coupled maps with $N$ units, which deal with rational distributions only. The definitions further include the restriction to some reduced mappings of the $(N-1)$-dimensional cube which capture any possible ergodicity failure. The results themselves describe the emergence of AsIUP in the lowest dimensions, namely dimension two and three, for certain  configurations (fairly general configurations in dimension two, uniform distribution in dimension three). The proofs, which proceed as indicated in the implementation scheme in the previous section, are all given in the Appendix.
 
\subsection{Systems of piecewise expanding globally coupled maps}\label{S-31}
{\bf Definition of the original system $F_{\rho,\epsilon}$:} An $N$-dimensional vector $\rho=(\rho_i)_{i=1}^N$ where all $\rho_i\geq 0$ and $\sum_{i=1}^N\rho_i=1$ is called a {\bf distribution}. Given a distribution and a number $\epsilon\in [0,\tfrac12)$, consider the mapping $F_{\rho,\epsilon}:\T^N\circlearrowleft$ defined by 
\[
(F_{\rho,\epsilon}(u))_i=2\left(u_i+\epsilon\sum_{j=1}^{N}\rho_i g(u_j-u_i)\right)\ \text{mod}\ 1,\ \forall i\in\{1,\cdots ,N\},\ u\in\T^N,
\] 
where $g$ represents pairwise elastic interactions on the circle \cite{KY10} and is defined by $g(u)=u-h(u)$ for all $u\in\T^1$ with
\[
h(u)=\left\{\begin{array}{ccl}
\lfloor u+\frac12\rfloor&\text{if}&u\not\in\frac12+\Z\\
0&\text{if}&u\in\frac12+\Z .
\end{array}\right.
\]
Here, $\lfloor \cdot \rfloor$ is the floor function. Hence $g$ is piecewise affine, with slope 1 and discontinuities at all points of $\frac12+\Z$. The symmetry $g(-u)=-g(u)\ \text{mod}\ 1$ implies that $F_{\rho,\epsilon}$ commutes with $-\text{\rm Id}$.

We shall denote by $F_{N,\epsilon}$ the restriction of $F_{\rho,\epsilon}$ to the {\bf uniform distribution} $\rho_i=\tfrac1{N}$ for all $i$. More precisely, the maps $F_{N,\epsilon}$, which are actually the  maps considered in \cite{F14}, are defined by
\[
F_{N,\epsilon}:=F_{(\tfrac1{N})_{i=1}^N,\epsilon}.
\]

For distributions $\rho$ with rational coordinates, the maps $F_{\rho,\epsilon}$ capture the so-called cluster dynamics in the maps $F_{N,\epsilon}$, namely the dynamics in some invariant subsets of $\T^N$. To see this, given any configuration $u=(u_i)_{i=1}^N\in\T^N$, let the distribution $(\tfrac{n_k}{N})_{k=1}^K$ be defined by the number $K\leq N$ of groups - called clusters\footnote{The meaning of the term {\sl cluster} here differs from the one in Section \ref{S-212}. However, no confusion should result since the context differs too.} - inside which the coordinates $u_i$ are equal, and by the number $n_k$ of coordinates in each group, see e.g.\ \cite{BKOVZ02}. The mean field coupling in $F_{N,\epsilon}$ implies that the set of configurations with given distribution $(\tfrac{n_k}{N})_{k=1}^K$ is invariant under the action of $F_{N,\epsilon}$ and the dynamics therein is governed by $F_{(\tfrac{n_k}{N})_{k=1}^K,\epsilon}$. 

Independently of the relationships with the dynamics of the $F_{N,\epsilon}$, upon selection of the distribution $\rho$, the mappings $F_{\rho,\epsilon}$ provide a natural instance where to investigate the emergence of AsIUP depending on additional symmetries of the dynamics. In particular, one may first consider the case of $F_{N,\epsilon}$ and then investigate AsIUP for non-uniform distributions $\rho$, either when any transformation in a sub-group of permutations commutes with the dynamics or when no other symmetry than $-\text{Id}$ prevails. The analysis of the examples below proceeds along these lines.
\medskip

\noindent
{\bf Definition of the reduced map $G_{\rho,\epsilon}$:} Ergodicity of the mappings $F_{\rho,\epsilon}$ can be examined in a more convenient family of piecewise affine mappings of the $(N-1)$-dimensional unit cube \cite{F20,SB16}, which turn out to be of the form considered in Section \ref{S-CONPRO} with $a=2(1-\epsilon)>1$. To see this, consider the change of variables $\pi_N$ of $\T^N$ defined by \cite{SB16}
\[
(\pi_Nu)_i=\left\{\begin{array}{ccl}
u_i-u_{i+1}\ \text{mod}\ 1&\text{if}&1\leq i\leq N-1\\
\sum_{j=1}^Nu_j\ \text{mod}\ 1&\text{if}&i=N
\end{array}\right.
\]
The mapping $\pi_N$ semi-conjugates $F_{\rho,\epsilon}$ to the direct product of a map $\tilde G_{\rho,\epsilon}:\T^{N-1}\circlearrowleft$ with the map $x\mapsto 2x\ \text{mod}\ 1$ of the unit circle, see \cite{F20,SB16} for more details in the case of $F_{N,\epsilon}$. Focusing on $\tilde G_{\rho,\epsilon}$, we consider its representation $G_{\rho,\epsilon}$ in the fundamental domain $M=(0,1)^{N-1}$ of the subset $\T_\ast^{N-1}$ of the torus whose configurations have no integer coordinates.

Letting $d=N-1$ , explicit computations yield that the expression of $G_{\rho,\epsilon}$ is as follows 
 \[
G_{\rho,\epsilon}=\text{\rm Frac}\left(2(1-\epsilon)\text{Id}+2\epsilon B_\rho \right)
\]
where for $x\in\R^{d}$, the notation $\text{\rm Frac}(x)=(x_i-\lfloor x_i\rfloor)_{i=1}^{d}$ stands for the vector with fractional part coordinates. The piecewise constant function $B_\rho$ reads
\[
(B_{\rho}(x))_i=(\rho_i+\rho_{i+1})h(x_i)+\sum_{j=1}^{i-1}\rho_j\left(h(\sum_{k=j}^{i}x_k)-h(\sum_{k=j}^{i-1}x_k)\right)+\sum_{j=i+1}^{d-1}\rho_{j+1}\left(h(\sum_{k=i}^{j}x_k)-h(\sum_{k=i+1}^{j}x_k)\right).
\]
for $i\in\{1,\cdots,d\}$. 

As before, for the uniform distribution $(\tfrac1{d+1})_{i=1}^{d+1}$, we use simplified notations, namely
\[
G_{d,\epsilon}:=G_{(\tfrac1{d+1})_{i=1}^{d+1},\epsilon}\quad \text{and}\quad B_d:=B_{(\tfrac1{d+1})_{i=1}^{d+1}}.
\]
Notice that $B_d(\R^d)\subset \tfrac{1}{d+1}\Z^d$.
\medskip

\noindent
{\bf Symmetries:} The change of variable $\pi_N$ implies that $G_{\rho,\epsilon}$ inherits the symmetries of $F_{\rho,\epsilon}$. In particular, the following properties hold.
\begin{itemize}
\item Let $\Sigma=1-\text{Id}$. We have $G_{\rho,\epsilon}\circ \Sigma|_{M\cap G_{\rho,\epsilon}^{-1}(M)}=\Sigma\circ G_{\rho,\epsilon}|_{M\cap G_{\rho,\epsilon}^{-1}(M)}$.\footnote{The restriction to $M\cap G_{\rho,\epsilon}^{-1}(M)$ incorporates the fact that $h(1-u)=1-h(u)$ iff $u\not\in\tfrac12+\Z$.} 
\item $G_{d,\epsilon}$ commutes with every element in a group of transformations that is isomorphic to $\Pi_{d+1}$. 
\end{itemize}

\subsection{AsIUP for two-dimensional mappings in the family of coupled maps}
In this paper, the conditioning problem approach to AsIUP construction in the coupled maps above addresses families of two and three-dimensional examples. This section deals with two-dimensional mappings. Following the lines above, we consider firstly maps $G_{(\rho_i)_{i=1}^3,\epsilon}$ with {\bf partly symmetric} distributions, namely those distributions for which $\rho_1=\rho_3$ (which forces $\rho_2=1-2\rho_1$). Then, we pass to maps $G_{(\rho_i)_{i=1}^3,\epsilon}$ with fully asymmetric distributions. In the three-dimensional case (next section), we shall only envisage the map $G_{3,\epsilon}$, for simplicity. 

\subsubsection{AsIUP for maps $G_{(\rho_i)_{i=1}^3,\epsilon}$ with distributions $\rho_1=\rho_3$}\label{S-RHO1-3}
Throughout this section we consider partly symmetric distributions $(\rho_i)_{i=1}^3$ so that they are described by a single parameter $\varrho:=\rho_1=\rho_3=\frac{1-\rho_2}2\in (0,\tfrac12)$. In particular, for $\varrho=\tfrac13$, one obviously obtains the uniform distribution for $N=d+1=3$.
\medskip

\noindent
{\bf Definitions:} A simple analysis concludes that the map $G_{\varrho,\epsilon}:=G_{(\rho_i)_{i=1}^3,\epsilon}$ explicitly writes
\[
(G_{\varrho,\epsilon}(x))_i=2(1-\epsilon)x_i+2\epsilon (B_{\varrho}(x))_i-h(x_i),\ i=1,2\quad \forall x=(x_i)_{i=1}^2\in (0,1)^2
\]
where $(B_{\varrho}(x))_i=(1-\varrho)h(x_i)+\varrho\left(h(x_1+x_2)-h(x_{3-i})\right)$. 
\medskip

\begin{figure}[ht]
\begin{center}
\includegraphics*[width=75mm]{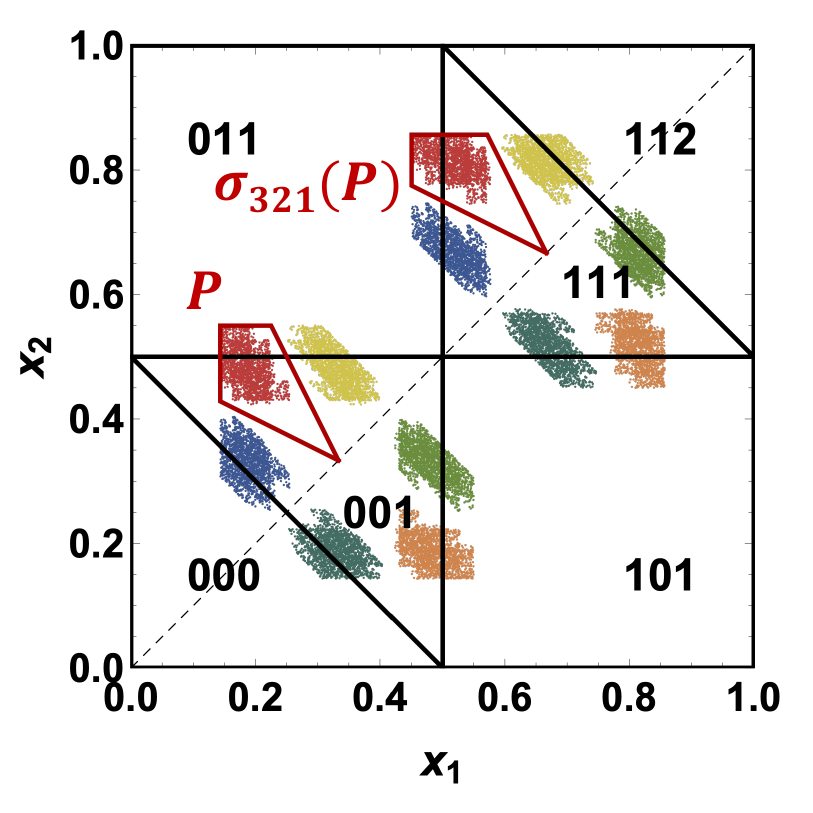}
\hspace{10mm}
\includegraphics*[width=75mm]{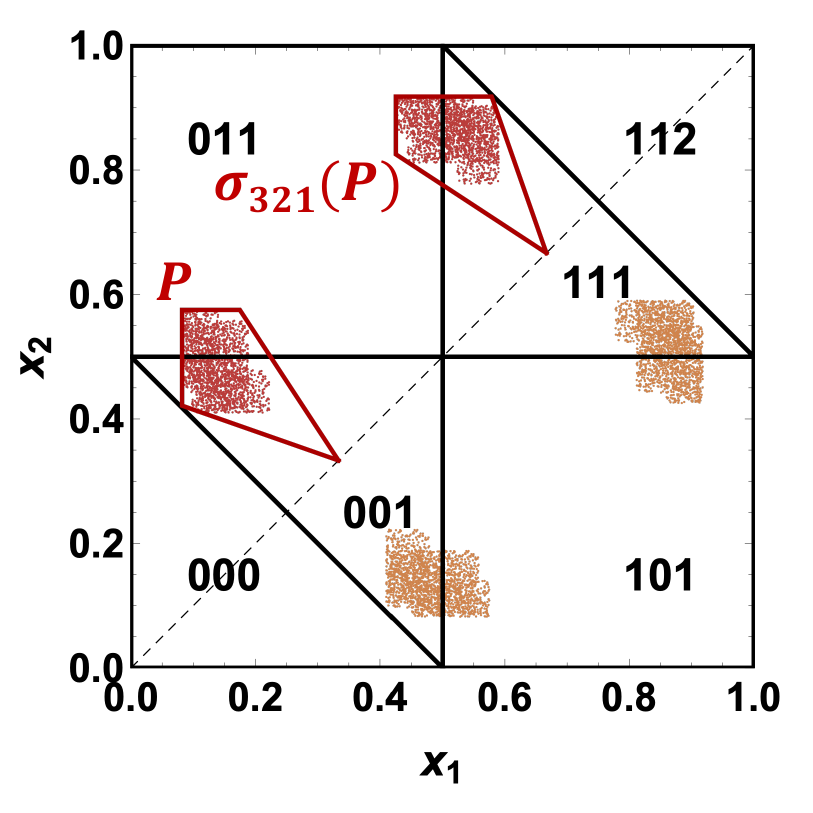}
\end{center}
\caption{Empirical ergodic components of $G_{\varrho,\epsilon}$ for $(\varrho,\epsilon)=(\tfrac13,0.43)$ (left, corresponding to the 3-unit globally coupled map $F_{3,\epsilon}$) and $(\varrho,\epsilon)=(\tfrac25,0.41)$ (right). Each component consists of $4\times 10^3$ consecutive orbit points (one color per orbit), after discarding transient behaviours (which are actually very short). One component is obtained by iterating $G_{\varrho,\epsilon}$ and the other ones follow from applying symmetries. 
The figure also displays a related maximal AsIUP $P\cup \sigma_{321}(P)$ (where $P=P^{\alpha_2}_{m_2}$ is the maximal solution of Conditioning Problem \ref{CONDPRO1}) which contains a single empirical ergodic component. By applying symmetries, more AsIUP can be obtained, in particular one for each of the five other ergodic components. 
Discontinuity lines $x_{1,2}=\frac12$ and $x_1+x_2=\frac12,\frac32$  and labels of the symbolic partition atoms, defined as concatenations of the values of $h(x_1),h(x_2)$ and $h(x_1+x_2)$, are also indicated, {\sl eg.}\ $000$, $001$, $011$ etc.}
\label{PhasePortraitG2}
\end{figure}

\noindent
{\bf Atomic partition:} For convenience, the atoms of the symbolic partition associated with $G_{\varrho,\epsilon}$ are labelled using concatenations of the values of $h(x_1),h(x_2)$ and $h(x_1+x_2)$. A simple geometric analysis (or use the systematic procedure described in Appendix \ref{A-PARTITION}) concludes that the square $(0,1)^2$ decomposes into 6 such atoms, namely (see Fig.\ \ref{PhasePortraitG2})
\[
{\cal A}=\left\{000, 001, 011, 101, 111, 112\right\}.
\]
\medskip

\noindent
{\bf Symmetries:} The map $\sigma_{321}$ acting in $(0,1)^2$, and induced, via the change of variables $\pi_3$, by the representation $(u_1,u_2,u_3)\mapsto (u_3,u_2,u_1)$ in $\T^3$ of the transposition $1\leftrightarrow 3$ (meaning that the units 1 and 3 are exchanged) in the original map $F_{(\rho_i)_{i=1}^3,\epsilon}$, writes
\[
\sigma_{321}(x)=(1-x_2,1-x_1).
\]
In addition to the commutation with $\Sigma$, for every $\varrho\in (0,\tfrac12)$, we have $G_{\varrho,\epsilon}\circ \sigma_{321}|_{M\cap G_{\rho,\epsilon}^{-1}(M)}=\sigma_{321}\circ G_{\varrho,\epsilon}|_{M\cap G_{\rho,\epsilon}^{-1}(M)}$ and $G_{\varrho,\epsilon}$ also commutes with the reflection $(x_1,x_2)\mapsto (x_2,x_1)$. 
\medskip

\noindent
{\bf State-of-the-art and contributions to the analysis of maps with partly symmetric distributions:} Numerical simulations of the dynamics of the uniform distribution map $G_{2,\epsilon}$ (obtained for $\varrho=\tfrac13$) have revealed the following features \cite{F14,SB16}. Ergodicity holds for $\epsilon < 0.417\ldots$. For larger values of $\epsilon$, the attractor decomposes into six ergodic components (see left panel in Figure \ref{PhasePortraitG2}). This phenomenology has been fairly well captured by mathematical statements. In particular, six distinct AsIUP have been identified, which proved to exist for $\epsilon\geq\tfrac{4-\sqrt{10}}2\sim 0.419$ \cite{SB16} (one of them is represented in the left panel of Fig.\ \ref{PhasePortraitG2}). 

While the exhibited AsIUP give a rather accurate description of the numerics, interrogations remain about their foundation and shape, in particular about those edges that are not aligned with any discontinuity lines of $G_{2,\epsilon}$, and about their persistence for more general distributions, when the symmetry group is smaller than $\Pi_3$, and in particular for $\varrho\neq \tfrac13$. (NB: The phenomenology for $\varrho\neq\tfrac13$ is similar to that of $G_{2,\epsilon}$, with the exception that only 2 two ergodic components may emerge. This is especially the case when $|\varrho-\tfrac13|$ is large, see right panel in Figure \ref{PhasePortraitG2}.)

As shown below, these questions will be answered for every value of $\varrho\in (0,\tfrac12)$, by solving the (reduced) IUP conditioning problem associated with the observations. Not only the AsIUP depicted in Fig.\ \ref{PhasePortraitG2} emerges as a solution, but a whole family of AsIUP is obtained, with distinct inclined lines (and whose existence domain in the coupling range depends on the inclination slope). Besides, it is also proved that no such IUP can exist, whose boundaries are only given by discontinuities' directions.
\medskip

\noindent
{\bf Conditioning problem for IUP of $G_{\varrho,\epsilon}$}: Given that the red trace in Fig.\ \ref{PhasePortraitG2} consists of two clusters and is globally invariant under the action of $\sigma_{321}$, given its location in the atomic partition and given the dynamics in the corresponding atoms, we consider the following reduced conditioning problem for the map $G_{\varrho,\epsilon}$:
\begin{CondProb}
Find a polytope $P\subset (0,1)^2$ that satisfies the following conditions: 
\begin{itemize}
\item $P\cap A_\omega\neq\emptyset$ iff $\omega\in\{001,011\}$.
\item $G_{\varrho,\epsilon}(P\cap A_{001})=\sigma_{321}(P)$ and $G_{\varrho,\epsilon}(P\cap A_{011})\subset P$.
\end{itemize}
\label{CONDPRO1}
\end{CondProb}
\noindent
If such a polytope $P$ exists, then commutation of $G_{\varrho,\epsilon}$ and $\sigma_{321}$ immediately implies that $P\cup \sigma_{321}(P)$ is an IUP for $G_{\varrho,\epsilon}$. The analysis of Conditioning Problem \ref{CONDPRO1} is given in Appendix \ref{A-P1} and yields the following conclusions.
\medskip

\noindent
{\bf Absence of solution associated with the canonical coefficient matrix:} 
\begin{Stat}
Let $\alpha$ be the canonical coefficient matrix associated with $G_{\varrho,\epsilon}$. Then no polytope of the type $P_m^\alpha$ can satisfy all conditions in Conditioning Problem \ref{CONDPRO1}.
\label{NOSOLUTIONP1}
\end{Stat}
The arguments in Appendix \ref{A-P1} in fact show that no solution of Problem \ref{CONDPRO1} can have edges along the direction of the anti-diagonal $x_1+x_2=\text{constant}$.
\medskip

\noindent
{\bf Families of solutions:} Given $a\in\R$, which we may assume $a>1$ w.l.o.g., let $\alpha_a$ be the matrix such that the polytopes $P^{\alpha_a}_m$ are defined as follows
\[
P^{\alpha_a}_m=\left\{x\in\R^2\ :\ \begin{array}{c}
\underline{m}_1<x_1<\overline{m}_1\\
\underline{m}_2<x_2<\overline{m}_2\\
\underline{m}_{a\cdot 1+2}<ax_1+x_2<\overline{m}_{a\cdot 1+2}\\
\underline{m}_{1+a\cdot 2}<x_1+ax_2<\overline{m}_{1+a\cdot 2}
\end{array}\right\}
\quad\text{where}\quad
m=\left(
\begin{array}{cc}
\underline{m}_1&\overline{m}_1\\
\underline{m}_2&\overline{m}_2\\
\underline{m}_{a\cdot 1+2}&\overline{m}_{a\cdot 1+2}\\
\underline{m}_{1+a\cdot 2}&\overline{m}_{1+a\cdot 2}
\end{array}\right)
\]
The matrix $\alpha_a$ is chosen in a way that the transformation acting on the corresponding constraint matrices, and induced by the symmetry $\sigma_{321}$ above, is $\alpha_a$-compatible.

In addition to the parameter $a$ which governs the polytopes $P^{\alpha_a}_m$ that have just been defined, recall that the parameter $\varrho=\rho_1=\rho_3=\frac{1-\rho_2}2$ characterizes the partly symmetric distribution $(\rho_i)_{i=1}^3$ and $\epsilon\in (0,\tfrac12)$ is the coupling strength.
\begin{Stat}
(i) Given $\varrho\in (0,\tfrac12),\ a>1$ and $\epsilon\in (0,\tfrac12)$, there exists at most one optimized constraint matrix $m_a$ such that $P^{\alpha_a}_{m_a}$ solves Conditioning Problem \ref{CONDPRO1}. The matrix $m_a$ is known explicitly.

\noindent
(ii) Given $\varrho\in \left(0,\tfrac12\right)$, the polytope $P^{\alpha_a}_{m_a}$ is not empty and solves Conditioning Problem \ref{CONDPRO1} iff $a\in \left(1,\min\left\{2,\tfrac{1+\varrho}{1-\varrho},\tfrac{\varrho}{3\varrho-1}\right\}\right]$ and $\epsilon\in [\epsilon_{\varrho,a},\tfrac12)$ for some uniquely defined $\epsilon_{\varrho,a}\in (0,\tfrac12)$. 

\noindent
(iii) When $P^{\alpha_a}_{m_a}$ is not empty, the union set $P^{\alpha_a}_{m_a}\cup \sigma_{321}(P^{\alpha_a}_{m_a})$ is an AsIUP with respect to $\Sigma=1-\text{Id}$ for $G_{\varrho,\epsilon}$.
\label{SOLP1}
\end{Stat}
\noindent
{\bf Solutions dependence on the parameters $\varrho$ and $\epsilon$:} For the uniform distribution $\varrho=\tfrac13$, we have $\epsilon_{\frac13,a}=\frac{3a+2-\sqrt{9a^2+4a-4}}4$ (see Appendix \ref{A-P1}), which decreases with $a\in \left(1,2\right]$\footnote{Notice that $\left(1,2\right]$ is the largest domain for the slope $a$, when $\varrho$ varies in $\left(0,\tfrac12\right)$. On the other hand, this domain shrinks to the singleton set $\{1\}$ when $\varrho\to 0$ and $\varrho\to\tfrac12$.} between the following values 
\[
\epsilon_{\frac13,2}=\tfrac{4-\sqrt{10}}2\quad \text{and}\quad \epsilon_{\frac13,1}=\tfrac12.
\]
For $\varrho$ close to $\tfrac13$, first order expansions (see Appendix \ref{A-P1}) yield
\[
\epsilon_{\varrho,a}=\epsilon_{\frac13,a}+\delta_a(\varrho-\tfrac13)+O\left((\varrho-\tfrac13)^2\right),
\]
(NB: here, the symbol $O$ stands for the big O notation) where $\delta_a$ is a negative constant for all $a\in (1,2]$, which decreases when $a$ increases. In particular, $\varrho\mapsto \epsilon_{\varrho,a}$ is decreasing for $\varrho$ near $\tfrac13$. In this neighbourhood, one can also show that $a\mapsto \epsilon_{\varrho,a}$ is decreasing; hence the coupling threshold for which we are sure that ergodicity breaking holds decreases when $\varrho$ increases. 

\begin{figure}[ht]
\begin{center}
\includegraphics*[width=70mm]{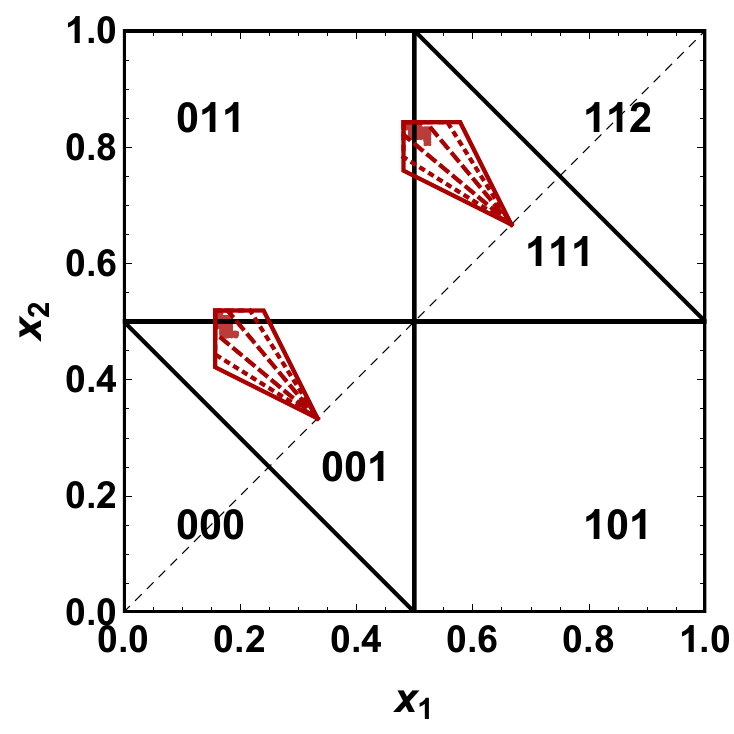}
\hspace{10mm}
\includegraphics*[width=70mm]{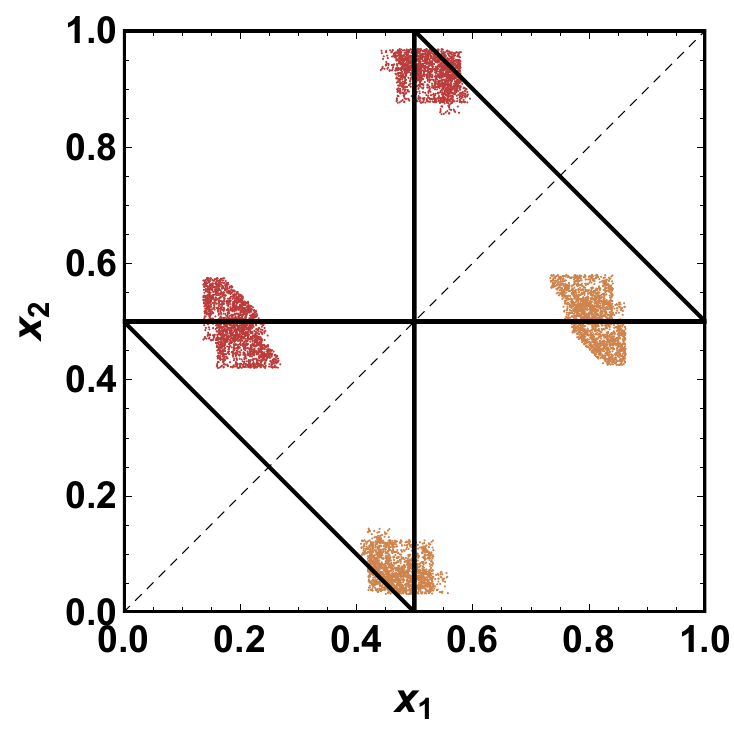}
\end{center}
\caption{{\sl Left.} Example of an empirical ergodic component of $G_{2,\epsilon}$ for $\epsilon=\epsilon_{a}\simeq 0.473$ for $a=1.2$, together with the AsIUP $P_{m_a}^{\alpha_a}$ for $a=1.2$ (dashed edges), $a=1.6$ (dotted edges) and $a=2$ (solid edges). {\sl Right.} Empirical ergodic components of $G_{\rho,\epsilon}$ for the distribution $\rho=\{0.463,0.2,0.337\}$ and $\epsilon=0.42$. Since the upper red cluster intersects the atom $A_{112}$, the conditioning problem of any IUP that contains this component must be distinct from Problem \ref{CONDPRO2}.}
\label{AsIUPD2Multi}
\end{figure}

The explicit expression of $m_a$ is given in Appendix \ref{A-P1} (see \eqref{ma}), which shows that the sets $P^{\alpha_a}_{m_a}$ are nested (see Fig.\ \ref{AsIUPD2Multi} left).\footnote{More precisely, we have $P_{m_a}^{\alpha_a}\subsetneq P_{m_{a'}}^{\alpha_{a'}}$ when $a<a'$, but the left edge $x_1=\epsilon(1-2\varrho)$, the top edge $x_2=1-2\epsilon\left(1-\varrho-\epsilon(1-2\varrho)\right)$ and the vertex $(\tfrac13,\tfrac13)$ of $P_{m_a}^{\alpha_a}$ do not depend on $a$.} 
In particular, for $\varrho=\tfrac13$, the previously obtained AsIUP for $G_{2,\epsilon}$, which is represented in Fig.\ \ref{PhasePortraitG2} right, is nothing but the largest one, namely $P^{\alpha_2}_{m_2}\cup \sigma_{321}(P^{\alpha_2}_{m_2})$. 

In addition, the results here imply that, for every $\epsilon\in [\epsilon_{\frac13,2},\tfrac12)$, there exists a (smallest) AsIUP $P^{\alpha_{a_\epsilon}}_{m_{a_\epsilon}}\cup \sigma_{321}(P^{\alpha_{a_\epsilon}}_{m_{a_\epsilon}})$ (where $a_\epsilon:=\min\{a\in \left(1,2\right]\ :\ \epsilon_{\frac13,a}\leq \epsilon\}$ is the reciprocal of $a\mapsto \epsilon_{\frac13,a}$), which is our best fit of the ergodic component (and appears to be a pretty accurate delimitation of it, see Fig.\ \ref{AsIUPD2Multi} left), and whose Lebesgue measure (area) vanishes as $\epsilon\to\tfrac12$.

\subsubsection{AsIUP for $G_{(\rho_i)_{i=1}^3,\epsilon}$ with asymmetric distributions} 
In this section, we consider the map $G_{(\rho_i)_{i=1}^3,\epsilon}$ for an (a priori) arbitrary distribution $(\rho_i)_{i=1}^3$, so that $\Sigma$ is the only symmetry. Notice that the atomic partition remains the same as the one of $G_{\varrho,\epsilon}$.
\medskip

\noindent
{\bf Conditioning problem for IUP:} The conditioning problem is obtained by duplicating the conditions in Conditioning Problem \ref{CONDPRO1}, taking into account that the symmetry $\sigma_{321}$ no longer applies. 
\begin{CondProb}
Find two polytopes $P_1,P_2\subset (0,1)^2$ which satisfy the following conditions: 
\begin{itemize}
\item $P_1\cap A_\omega\neq \emptyset$ iff $\omega\in\{001,011\}$ and $P_2\cap A_\omega\neq \emptyset$ iff $\omega\in\{111,011\}$,
\item $G_{(\rho_i)_{i=1}^3,\epsilon}(P_1\cap A_{001})=P_2$, $G_{(\rho_i)_{i=1}^3,\epsilon}(P_2\cap A_{111})=P_1$ and $G_{(\rho_i)_{i=1}^3,\epsilon}(P_j\cap A_{011})\subset P_j$ for $j=1,2$.
\end{itemize}
\label{CONDPRO2}
\end{CondProb}

\noindent
{\bf Families of solutions for weakly asymmetric distributions:} By using continuation arguments for distributions $(\rho_i)$ with $|\rho_1-\rho_3|$ small, namely distributions that are defined by $\rho_1=\varrho+\delta$, $\rho_3=\varrho-\delta$ and then $\rho_2=1-2\varrho$, for some arbitrary $\varrho\in (0,\tfrac12)$ and $|\delta|$ small, one obtains the following statement about the solutions of the Conditioning Problem \ref{CONDPRO2}. Recall that the parameter $a$ governs the polytopes $P^{\alpha_a}_m$ and $\epsilon$ is the coupling strength.
\begin{Stat}
{\sl (i)} Given $\varrho\in (0,\tfrac12),\ |\delta|<\varrho,\ a>1$ and $\epsilon\in (0,\tfrac12)$, there exists at most one pair $(m_1,m_2)$ such that $(P_{m_1}^{\alpha_a},P_{m_2}^{\alpha_a})$ solves the Conditioning Problem \ref{CONDPRO2} for the map $G_{(\rho_i)_{i=1}^3,\epsilon}$ with $\rho_1=\varrho+\delta$ and $\rho_3=\varrho-\delta$. The matrices $m_1$ and $m_2$ depend continuously on $\delta$.

\noindent
{\sl (ii)} For every $\varrho\in \left(0,\tfrac12\right)$, $a\in\left(1,\min\left\{2,\tfrac{1+\varrho}{1-\varrho},\tfrac{\varrho}{3\varrho-1}\right\}\right]$ and $\epsilon\in \left(\epsilon_{\varrho,a},\tfrac12\right)$, there exists $\Delta>0$ such that the polytopes $(P_{m_1}^{\alpha_a},P_{m_2}^{\alpha_a})$ are non-empty and solve the Conditioning Problem \ref{CONDPRO2} iff $|\delta|<\Delta$. For other values of $a$ and/or $\epsilon$, the conditioning problem has no solution $(P_{m_1}^{\alpha_a},P_{m_2}^{\alpha_a})$.

\noindent
{\sl (iii)} When $(P_{m_1}^{\alpha_a},P_{m_2}^{\alpha_a})$ are non-empty, the union set $P_{m_1}^{\alpha_a}\cup P_{m_2}^{\alpha_a}$ is an AsIUP  with respect to $\Sigma=1-\text{Id}$ for the map $G_{(\rho_i)_{i=1}^3,\epsilon}$ with $\rho_1=\varrho+\delta$ and $\rho_3=\varrho-\delta$.
\label{SOLP2}
\end{Stat}
The proof is given in Appendix \ref{A-P2}. Notice that while the restriction $\rho_1\sim\rho_3$ in this statement is a by-product of the continuation arguments used in the proof, some constraint should be imposed on $|\rho_1-\rho_3|$ for Conditioning Problem \ref{CONDPRO2} to have solution. Indeed, the right panel in Fig.\ \ref{AsIUPD2Multi} suggests that this problem has no solution when $|\rho_1-\rho_3|$ becomes large. To investigate a suitable conditioning problem in this case, and more generally to investigate all conditioning problems for AsIUP associated with $G_{\rho,\epsilon}$ will be the subject of future studies.

\subsection{AsIUP for the map $G_{3,\epsilon}$}
{\bf Definitions:} The map $G_{3,\epsilon}$, which corresponds to the 4-unit globally coupled map $F_{4,\epsilon}$, explicitly writes
\[
(G_{3,\epsilon}(x))_i=2(1-\epsilon)x_i+2\epsilon (B_3(x))_i-h(x_i),\ i=1,2,3,\quad \forall x=(x_i)_{i=1}^3\in (0,1)^3
\]
where
\[
\left\{\begin{array}{l}
(B_{3}(x))_1=\tfrac{2h(x_1)-h(x_2)+h(x_1+x_2)-h(x_2+x_3)+h(x_1+x_2+x_3)}4\\
(B_{3}(x))_2=\tfrac{2h(x_2)-h(x_1)-h(x_3)+h(x_1+x_2)+h(x_2+x_3)}4\\
(B_{3}(x))_3=\tfrac{2h(x_3)-h(x_2)-h(x_1+x_2)+h(x_2+x_3)+h(x_1+x_2+x_3)}4
\end{array}\right.
\]
\medskip

\noindent
{\bf Atomic partition:} The atoms of the symbolic partition associated with $G_{3,\epsilon}$ are labelled using concatenations of the values of $h(x_1),h(x_2),h(x_3),h(x_1+x_2),h(x_2+x_3)$ and $h(x_1+x_2+x_3)$.
A systematic analysis based on the optimization function (see Appendix \ref{A-PARTITION}) concludes that the cube $(0,1)^3$ decomposes into 26 atoms. For the sake of space, only those atoms involved in the conditioning problems will be considered here.
\medskip

\noindent
{\bf Symmetries:} The permutation group $\Pi_4$ can be generated by transpositions. Here, we shall be especially interested in the transformations $\sigma_{4231}$ and $\sigma_{1324}$ acting in $(0,1)^3$, and respectively induced by the transpositions $(u_1,u_2,u_3,u_4)\mapsto (u_4,u_2,u_3,u_1)$ and $(u_1,u_2,u_3,u_4)\mapsto (u_1,u_3,u_2,u_4)$ in $\T^4$. They write
\[
\sigma_{4231} (x)=\left(-x_2-x_3,x_2,-x_1-x_2\right)\quad\text{and}\quad \sigma_{1324} (x)=\left(x_1+x_2,-x_2,x_2+x_3\right).
\]
Let also
\[
\sigma_{4321}:=\sigma_{4231}\circ\sigma_{1324}=\sigma_{1324}\circ \sigma_{4231},\quad \text{ie.}\quad 
\sigma_{4321}(x)=\left(-x_3,-x_2,-x_1\right).
\]
In the analysis of the second bifurcation, we shall use the transformations $\sigma_{2134}$ induced by $(u_1,u_2,u_3,u_4)\mapsto (u_2,u_1,u_3,u_4)$ and $\sigma_{3124}:=\sigma_{2134}\circ \sigma_{1324}$ whose explicit expressions are
\[
\sigma_{2134}(x)=\left(-x_1,x_1+x_2,x_3\right)\quad\text{and}\quad \sigma_{3124}(x)=\left(-x_1-x_2,x_1,x_2+x_3\right).
\]
\medskip

\begin{figure}[ht]
\begin{center}
\includegraphics*[width=70mm]{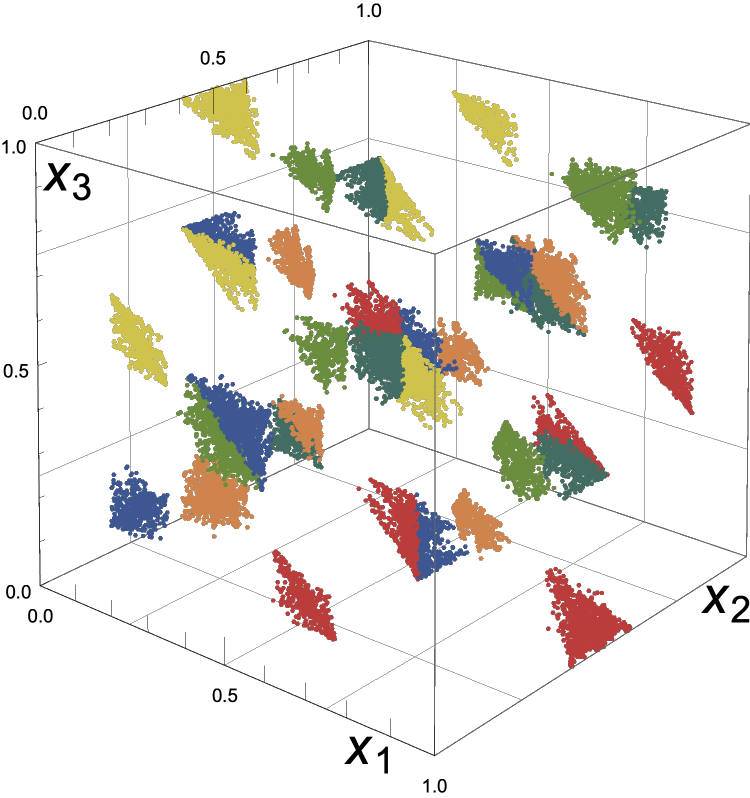}
\hspace{10mm}
\includegraphics*[width=82mm]{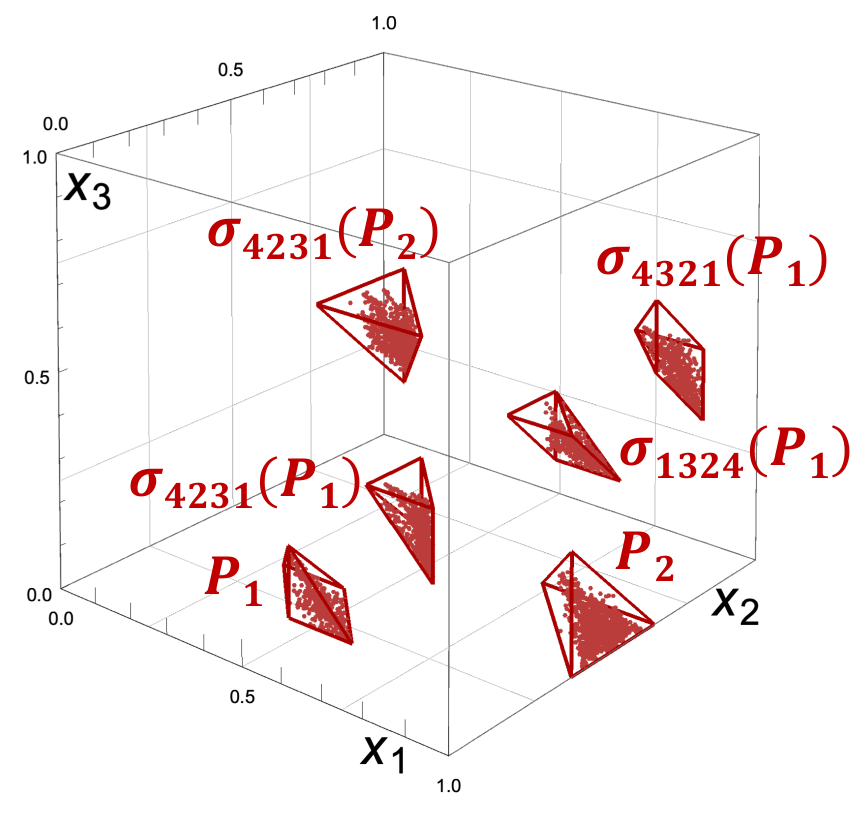}
\end{center}
\caption{{\sl Left.}\ Asymmetric empirical ergodic components of $G_{3,\epsilon}$ for $\epsilon=0.4$. Each component consists of $4\times 10^3$ consecutive orbit points (one color per orbit). One component is obtained by iterating $G_{3,\epsilon}$ and the other ones follow from applying symmetries. (NB: $G_{3,\epsilon}$ also has a symmetric ergodic component \cite{F20,S18}, which is not represented on this picture.) {\sl Right.}\ Related maximal AsIUP $\text{Orb}_{\langle \sigma_{4231},\sigma_{1324}\rangle}(P_1,P_2)$ (where $(P_1,P_2)=(P_{m_1(0)}^\alpha,P_{m_2(0)}^\alpha$) is the largest solution of the Conditioning Problem \ref{CONDPRO3}) which contains the red ergodic component (and that one only). Of note, the atomic partition of $G_{3,\epsilon}$ is not indicated here for the sake of clarity; see \cite{S18} for a representation of this partition. }
\label{PhasePortraitG3}
\end{figure}

\begin{figure}[ht]
\begin{center}
\includegraphics*[width=70mm]{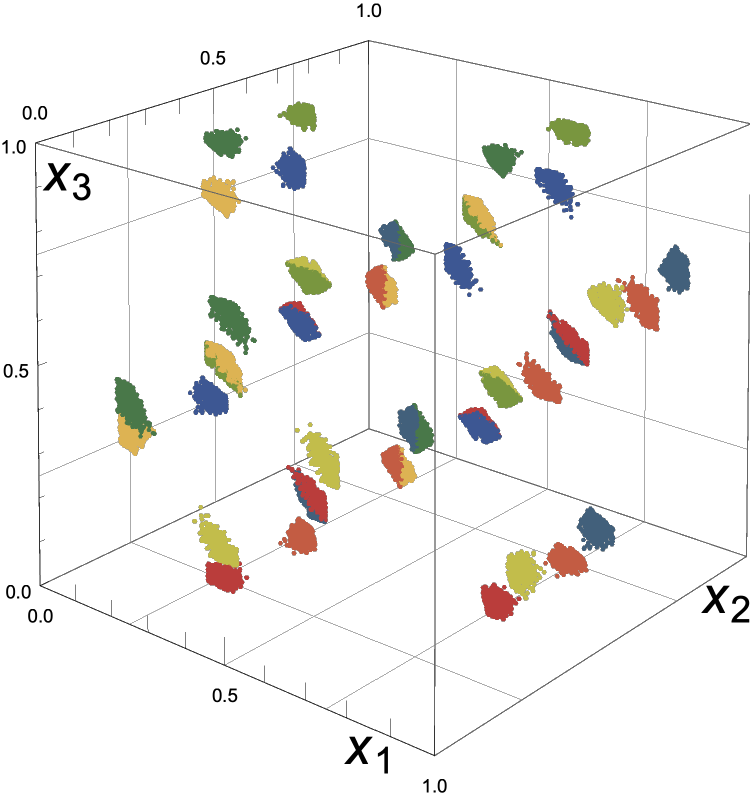}
\hspace{10mm}
\includegraphics*[width=77mm]{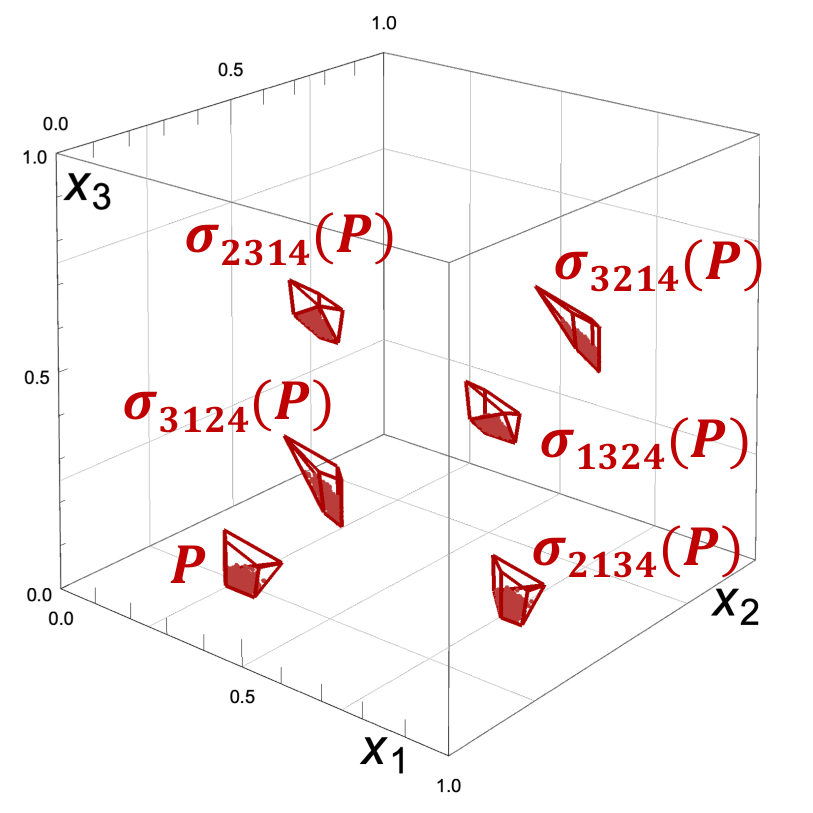}
\end{center}
\caption{{\sl Left.}\ Additional asymmetric empirical ergodic components of $G_{3,\epsilon}$ for $\epsilon=0.44$ (NB: For the sake of clarity, the symmetric ergodic component and the asymmetric components that emerge at $\epsilon\sim 0.397$ are not represented on this picture). 
{\sl Right.}\ Related AsIUP $\text{Orb}_{\langle \sigma_{2134},\sigma_{1324}\rangle}(P)$, where $P$ solves the Conditioning Problem \ref{CONDPRO4}, that contains a single one of these ergodic components.}
\label{PhaPorG3-2}
\end{figure}

\noindent
{\bf Basic phenomenology and state-of-the-art:} Numerical simulations of the dynamics of $G_{3,\epsilon}$ have revealed the following features \cite{F20,S18}. Unlike for $G_{2,\epsilon}$, a symmetric ergodic component exists for all $\epsilon\in [0,\tfrac12)$. For $\epsilon < 0.393\ldots$, this component is unique and ergodicity holds (see \cite{S18} for mathematical results related to the existence of a symmetric IUP). For larger values of $\epsilon$, six asymmetric components emerge, which persist for $\epsilon$ up to $\tfrac12$ (see left panel in Figure \ref{PhasePortraitG3}). At $\epsilon\sim 0.437$, an additional group of 8 asymmetric ergodic components emerges (see Fig.\ \ref{PhaPorG3-2}), which also persists for $\epsilon$ up to $\tfrac12$. 

\subsubsection{Conditioning problem for IUP of $G_{3,\epsilon}$ that emerges at the first bifurcation}  
The observed first bifurcation has been captured by mathematical statements. In particular, six distinct AsIUP have been identified, which proved to exist for $\epsilon\geq \epsilon_3$ \cite{S18} (where $\epsilon_3 \sim 0.397$ is the real root of some cubic polynomial, see expression before Statement \ref{SOLP3} below); one of them is represented in the right panel of Fig.\ \ref{PhasePortraitG3}. While the corresponding faces are all aligned with discontinuity planes, interrogations remain about the foundation and uniqueness of such invariant sets. As before, these questions can be answered by solving the following (reduced) IUP conditioning problem associated with the numerical observations. 
\begin{CondProb}
Find two polytopes $P_1,P_2\subset (0,1)^3$ which satisfy the following conditions: 
\begin{itemize}
\item $P_{1}\cap A_\omega\neq \emptyset$ iff $\omega\in \{000101,100101\}$
\item $\sigma_{4321}(P_2)=P_2$ and $P_{2}\cap A_\omega\neq\emptyset$ iff $\omega\in\{110111, 110112, 110212, 100101, 100111, 100112\}$.
\item $G_{3,\epsilon}(P_{1}\cap A_{000101})\subset P_{2}$ and $G_{3,\epsilon}(P_{1}\cap A_{100101})\subset P_{1}$.
\item $G_{3,\epsilon}(P_{2}\cap A_{110111})=P_{1}$, $G_{3,\epsilon}(P_{2}\cap A_{110112})\subset \sigma_{4231}(P_{1})$, $G_{3,\epsilon}(P_{2}\cap A_{110212})\subset P_{2}$.
\end{itemize}
\label{CONDPRO3}
\end{CondProb}
The transformations $\sigma_{4231}$ and $\sigma_{4321}$ are both $\alpha$-compatible for the canonical coefficient matrix $\alpha$ associated with $G_{3,\epsilon}$. From the expressions of the constraints matrices associated with the atoms, given at the beginning of Appendix \ref{A-P3}, it results that we have
\[
\sigma_{4321}(A_{110111})=A_{100112},\quad \sigma_{4321}(A_{110112})=A_{100111}\quad\text{and}\quad \sigma_{4321}(A_{110212})=A_{100101},
\]
which easily imply that, when such $P_1$ and $P_2$ exist, the corresponding orbit set under the sub-group $\langle \sigma_{4231},\sigma_{1324}\rangle$, namely\footnote{Notice that $\sigma_{4321}(P_2)=P_2$ implies $\sigma_{4231}(P_2)=\sigma_{1324}(P_2)$.}
\[
\text{Orb}_{\langle \sigma_{4231},\sigma_{1324}\rangle}(P_1,P_2)= P_{1}\cup \sigma_{4231}(P_{1})\cup \sigma_{1324}(P_{1})\cup \sigma_{4321}(P_{1})\cup P_{2}\cup \sigma_{4231}(P_{2}),
\]
is an IUP of $G_{3,\epsilon}$. 

The analysis reported in Appendix \ref{A-P3} reveals that the Conditioning Problem \ref{CONDPRO3} has indeed solutions of the form $(P^\alpha_{m_1},P^\alpha_{m_2})$ for the canonical matrix $\alpha$ of $G_{\epsilon,3}$.\footnote{We do not know whether or not this problem admits solutions for other coefficient matrices.} Let $\epsilon_3$ be the real root of the cubic polynomial
\[
4\epsilon^3-14\epsilon^2+15\epsilon-4.
\]
\begin{Stat}
Let $\alpha$ be the canonical coefficient matrix associated with $G_{\epsilon,3}$. A solution $(P_1,P_2)=(P^\alpha_{m_1},P^\alpha_{m_2})$ of Conditioning Problem \ref{CONDPRO3} exists iff $\epsilon\in \left[\epsilon_3,\tfrac12\right)$. Moreover, every such solution writes $(P^\alpha_{m_1(\delta)},P^\alpha_{m_2(\delta)})$ for an arbitrary $\delta:=(\delta_i)_{i=1}^5\in\Delta_\epsilon$, where the five-parameter families of constraint matrices $m_1(\cdot),m_2(\cdot)$ and the set $\Delta_\epsilon$ are respectively given by equations \eqref{m1} and \eqref{m2}, and by equation \eqref{DELTAEPS} in Appendix \ref{A-P3}. The corresponding orbit set 
\[
\text{Orb}_{\langle \sigma_{4231},\sigma_{1324}\rangle}(P^\alpha_{m_1(\delta)},P^\alpha_{m_2(\delta)})
\]
is an AsIUP  with respect to $\Sigma=1-\text{Id}$ for $G_{3,\epsilon}$. 
\label{SOLP3}
\end{Stat}
Moreover, we have $0\in \Delta_\epsilon$ for all $\epsilon \in \left[\epsilon_3,\tfrac12\right)$, $\Delta_{\epsilon_3}=\{0\}$ and equations \eqref{m1} and \eqref{m2} show that the polytopes $P^\alpha_{m_i(\delta)}$ are nested; in particular, we have
\[
m_i(\delta)\subset m_i(0)\ \text{for}\ i=1,2,\ \forall \delta\in \Delta_\epsilon.
\]
As for the solution in Statement \ref{SOLP1} above, best fit of the empirical ergodic component by the AsIUP here can accordingly be obtained by optimizing over $\delta\in\Delta_\epsilon$.\footnote{In practice, the smallest AsIUP might not be easily computed because the constraints on the coordinates of $\Delta$ are inter-dependent.} Besides, the AsIUP obtained in \cite{S18}, which is represented in Fig.\ \ref{PhasePortraitG3} right, is nothing but 
\[
\text{Orb}_{\langle \sigma_{4231},\sigma_{1324}\rangle}(P^\alpha_{m_1(0)},P^\alpha_{m_2(0)}).
\]

\subsubsection{Conditioning problem for IUP of $G_{3,\epsilon}$ that emerges at the second bifurcation}  
The additional asymmetric ergodic components that emerge in the numerics for $\epsilon> 0.437$ were clearly identified in \cite{S18}, with accurate localisation and dynamics. Nevertheless, to prove the existence of an IUP that would contain a single one of these components remained unsolved, especially because some of the directions of their faces remained evasive. Here, we first confirm that such IUP cannot be captured by using only the canonical coefficient matrix. Then we provide the missing faces and we check that a given candidate polytope indeed solves the following conditioning problem.

As usual, the following conditioning problem has been elaborated based on the information contained in numerical simulations.
\begin{CondProb}
Find a polytope $P\subset (0,1)^3$ which satisfies the following conditions
\begin{itemize}
\item $P\cap A_\omega\neq\emptyset$ iff $\omega\in \{000000, 000001,000101\}$,
\item $G_{3,\epsilon}(P\cap A_{000000})\subset P\cap A_{000101}$,  $G_{3,\epsilon}(P\cap A_{000001})\subset \sigma_{3124}(P)$ and $G_{3,\epsilon}(P\cap A_{000101})\subset \sigma_{2134}(P)$.
\end{itemize}
\label{CONDPRO4}
\end{CondProb}
When such a polytope exists, then the orbit set\footnote{Recall that $\sigma_{3124}=\sigma_{2134}\circ\sigma_{1324}$.} $\text{Orb}_{\langle \sigma_{2134},\sigma_{1324}\rangle}(P)$ (which consists of six elements) is clearly an IUP for $G_{3,\epsilon}$.
\begin{Stat}
Let $\alpha$ be the canonical coefficient matrix associated with $G_{3,\epsilon}$. Then no polytope of the type $P_m^\alpha$ can satisfy all conditions in Conditioning Problem \ref{CONDPRO4}.
\label{NOSOLUTIONP4}
\end{Stat}
The proof, given in Appendix \ref{A-NOSOLUTIONP4}, in fact shows that no solution of Conditioning Problem \ref{CONDPRO4} can have faces of the type $x_1=\text{constant}$.
\bigskip

For simplicity, we provide here a solution obtained from direct geometric considerations in phase space and use the constraint matrix formalism to compute its existence condition.
Let $p^\ast=\tfrac{2-\epsilon}{2(3-2\epsilon)}$\footnote{Notice that $1-3p^\ast<0<1-2p^\ast$ for $\epsilon\in (0,\tfrac12)$.} and 
\[
p_1^\ast=1-p^\ast,\quad p_2^\ast=1-2p^\ast,\quad\text{and}\quad p_3^\ast=1-3p^\ast,
\]
and let $P$ be the polytope resulting from the intersection of the half-spaces defined by the following inequalities
\begin{align}
&\tfrac{\epsilon}2\leq x_2,\quad 0\leq x_3,\quad \tfrac{\epsilon}2(3-2\epsilon)\leq x_1+x_2,\quad
x_1 + 2x_2 + 3 x_3\leq 1,\quad p_1^\ast x_1+p_2^\ast x_2+p_3^\ast x_3\leq p_2^\ast\nonumber\\
&\text{and}\quad -p^\ast x_1+p_2^\ast x_2+p_3^\ast x_3\leq 0.
\label{DEFP}
\end{align}
\begin{Stat}
The polytope $P$ is not empty and solves Conditioning Problem \ref{CONDPRO4} iff $\epsilon\in \left[\tfrac{5-\sqrt{17}}2,\tfrac12\right)$.\footnote{$\tfrac{5-\sqrt{17}}2\sim 0.438$} In this case, the orbit set $\text{Orb}_{\langle \sigma_{2134},\sigma_{1324}\rangle}(P)$ is an AsIUP of $G_{3,\epsilon}$.
\label{SOLP4}
\end{Stat}
The structure of the proof is presented in Appendix \ref{A-SOLP4}. Details of the calculations are provided in the Mathematica notebook "{\tt AsIUPG3SecondBif.nb"} (given in Supplementary Material).

To conclude this section, we mention that, when non-empty, the polytope $P$ actually satisfies the equality $G_{3,\epsilon}(P\cap A_{000101})= \sigma_{2134}(P)$, and not only the last inclusion in Conditioning Problem \ref{CONDPRO4}. A sketch of the proof of this claim is given at the end of Appendix \ref{A-SOLP4}.

\section{Concluding remarks and open questions}\label{S-CONCL}
Motivated by proofs of loss of ergodicity in systems of coupled maps, we have developed a systematic approach to the construction of IUP in piecewise affine and expanding systems, based on empirical information from numerical simulations. 

As a proof of concept, the approach has been applied to previously considered examples of globally coupled maps with $N=3,4$ units. It has provided justifications of the results in \cite{F14,S18,SB16} and has complemented them by addressing more general distributions and previously open questions. The approach has also revealed unexpected features of the dynamics, such as the impossibility, in some cases, of IUP with polytopes faces that are all aligned with some atomic faces. Altogether, the obtained accurate existence conditions and accurate fits of the numerics have demonstrated the feasibility and the efficiency of this approach. To a large extent, the numerical phenomenology has been confirmed by mathematical results.

Based on these results, future studies may envisage to address globally coupled maps for larger numbers of units, and especially for population sizes that have so far not be reached by the computed-assisted proof in \cite{F20}. Such a higher dimensional application would probably need to rely on formal computation tools, such as the Mathematica notebook employed in the proof of Statement \ref{SOLP4}.  The lists of instructions in Section \ref{S-ANALYSIS} and the systematic collection of empirical information in Appendix \ref{A-SYSTEMATIC} have been designed to anticipate such applications in higher/arbitrary dimension. 

Beyond systems of coupled maps, the approach can in principle be applied to any mapping in $\R^d$ of the form $a\text{Id}+B$ with $a>1$ with $B$ being piecewise constant on convex polytopes. Such a broad potential calls for a series of questions related to its theoretical foundations. 

Probably the most prominent theoretical challenge is to ensure the existence of (non-trivial) IUP for arbitrary mappings of the form $a\text{Id}+B$, and in particular to prove that every forward invariant set of such mappings is contained in some IUP. 

Equally relevant for our purpose would be to assert that every IUP can be obtained as a solution of a conditioning problem, {\sl ie.}\ to specify the number of polytopes, their location and the relative transitions between these sets suffices to determine an IUP. As the analysis of the two-dimensional map in Section \ref{S-RHO1-3} shows, this not only means to determine the constraint matrices $m$ but also the supporting coefficient matrix $\alpha$. 

In the examples, comparison of the IUP with the numerics shows that, given the limited number of faces that have been employed, the resulting fits appear to be somewhat accurate, especially as certain faces are concerned. However, the computed solutions do not always seem to exploit the full potential of polytopes associated with a given coefficient matrix. This raises the question of fit improvement. 

For instance, the coefficient matrix $\alpha_a$ in Section \ref{S-RHO1-3} allows for hexagons. Fig.\ \ref{PhasePortraitG2} suggests that such sets might provide better fits of the numerics of $G_{\varrho,\epsilon}$. Since the quadrilaterals described in Statement \ref{SOLP1} are the unique solutions of Conditioning Problem \ref{CONDPRO1}, given such matrix $\alpha_a$, one may investigate if, by imposing more topological/geometric conditions, hexagonal solutions coud result. More generally, given a coefficient matrix, one may intend to increase, if possible, the number of faces of the solution polytopes, in order to reduce the Lebesgue measure of the IUP. 

Alternatively, given a conditioning problem, fit improvement might result by 
considering coefficient matrices with larger number of rows. Can this process be iterated in order to obtain an arbitrary approximation of a Lebesgue ergodic component (or a union of such components)?

We conclude by a suggestion for an additional application. Claim \ref{ADAPTMAT} implies that the approach equally applies to piecewise affine mappings of the form $a\sigma+B$ provided that $\sigma$ generates a finite group. In particular, the approach could be applied to investigate the loss of ergodicity in certain piecewise isometries ({\sl viz.} $a=1$), the dynamics of which remains largely unknown beyond their topological entropy \cite{B01,G00}.
\medskip

\noindent
{\bf Acknowledgments}

\noindent
We are thankful to Anthony Quas for stimulating discussions and to P\'eter B\'alint for continuous support and incentives, especially as applications to coupled maps are concerned. We are also grateful to No\'e Cuneo, Stan Mintchev and Matteo Tanzi for their critical reading of the manuscript and suggestions for improvements. Several comments and open questions are direct consequences of their feedback. The research of FS is supported by the European Research Council (ERC) under the European Union's Horizon 2020 research and innovation programme (grant agreement No 787304).

\appendix

\section{Systematic procedure to define a conditioning problem and numerical considerations}\label{A-SYSTEMATIC}
The purpose of this Appendix is to provide a list of instructions for the setting up of a conditioning problem from the knowledge of a numerical orbit (which we assume to generate one empirical ergodic component - the procedure is similar for the union of several ergodic components). The instructions include the case of a reduced conditioning problem, to which considerations can be limited in the presence of symmetries. The instructions are provided in such a way that the procedure can be automatized.

Some instructions rely on a clustering algorithm. Such an algorithm typically involves a parameter, which quantifies the distance between any point and its closest neighbor. In order to identify the clusters, one starts from a large value of this parameter and decreases this value until the number of identified clusters ({\sl ie.}\ sets of points for which the minimal pairwise distance is not larger than the parameter) reaches a plateau. The value of this plateau is retained as the number of clusters (unless it is comparable to the number of orbit points.) 

Similar considerations apply to testing inclusion of (the image of) one cluster into another one. This implies considering inclusion of every point in the latter inside one sufficiently small ball around a point in the former. Once inclusion is asserted, equality can be also tested in a similar way, for those clusters that appear to have (hyper-)planar faces.
\medskip

Assuming that a clustering algorithm has been identified, the (macroscopic) clusters of an empirical orbit (whose number should be independent of the number of orbit points, provided that the latter is sufficiently large), the instructions can be enumerated as follows:
\begin{itemize}
\item {\sl in presence of symmetries:} 
\begin{itemize}
\item[$\ast$] Identify the symmetries of the cluster collection, namely those (sub-group of) symmetry transformations that leave invariant this collection\footnote{By invariance, we mean here invariance up to some accuracy in the location of individual points}
\item[$\ast$] Identify symmetric relationships between clusters, {\sl ie.}\ those clusters that are images of other, or simply invariant, under the previously identified symmetries. 
\end{itemize}
\item Identify those atoms each cluster intersects.
\item Compute the transition graph associated with clusters, {\sl ie.}\ identify localisation of the image of each atomic cluster component.
\item Suppress redundant dynamical information, {\sl ie.}\ those transitions that follow from other ones by applying the dynamics (or symmetries).
\end{itemize}

\section{Systematic procedure for the symbolic partition of $G_{d,\epsilon}$}\label{A-PARTITION}
Given the piecewise constant function $B_d$ in the map $G_{d,\epsilon}$ defined in Section \ref{S-COUPLEDMAP}, constraint matrices associated with the canonical coefficient matrix can be written as $m=(\underline{m}_{i+\cdots +j}\ \overline{m}_{i+\cdots +j})_{1\leq i\leq j\leq d}$ and the corresponding polytopes are defined by \footnote{The symbol $i+\cdots+j$ is a shortcut for the string $i+(i+1)+\cdots+(j-1)+j$, which is itself motivated by the expanded expression of the sum involved in the definition $P_m$,  ie.\ $\sum_{k=i}^jx_k=x_i+\cdots +x_j$.}
\[
P_m:=\left\{x\in\R^d:\underline{m}_{i+\cdots +j}<\sum_{k=i}^jx_k<\overline{m}_{i+\cdots +j}, 1\leq i\leq j\leq d\right\}.
\]
In this setting, the following procedure determines the symbolic partition of $G_{d,\epsilon}$ for arbitrary $d\in\N$ and, in particular the constraint matrices associated with each atom. It consists of the following operations:
\begin{itemize}
\item List all {\sl a priori} possible values of the vector $\left(h(\sum_{k=i}^jx_k)\right)_{1\leq i\leq j\leq d}$ for $x\in (0,1)^{d}$ by using that each $h(\sum_{k=i}^jx_k)$ may take any value in $\{0,\cdots, j-i+1\}$.
\item For each vector $\left(h(\sum_{k=i}^jx_k)\right)_{1\leq i\leq j\leq d}$,
\begin{itemize}
\item[$\ast$] Define the constraint matrix $m^\ast$ with entries 
\[
\underline{m}_{i+\cdots +j}^\ast=\max\{0,h(\sum_{k=i}^jx_k)-\tfrac12\}\quad \text{and}\quad \overline{m}_{i+\cdots +j}^\ast=\min\{h(\sum_{k=i}^jx_k)+\tfrac12,j-i+1\},\ 1\leq i\leq j\leq d
\]
\item[$\ast$] Compute the optimized matrix $O(m^\ast)$ and test $P_{m^\ast}\neq \empty 0$. If the test is positive, then retain $P_{O(m^\ast)}$ as the atom associated with the vector $\left(h(\sum_{k=i}^jx_k)\right)_{1\leq i\leq j\leq d}$.
\end{itemize}
\end{itemize}

\section{Analysis of Conditioning Problem \ref{CONDPRO1}}\label{A-P1}
This appendix presents the analysis of Conditioning Problem \ref{CONDPRO1} which in particular yields the proofs of Statements \ref{NOSOLUTIONP1} and \ref{SOLP1}. The analysis follows the plan described in Section \ref{S-ANALYSIS}. 

For the sake of space, we shall only consider the matrix $\alpha_a$ defined in Section \ref{S-RHO1-3}, assuming $a\geq 1$. Considerations about the canonical matrix associated with the symbolic partition (and the about corresponding constraint matrices in $\R^{3\times 2}$) will be assumed by letting $a=1$ and by using the following identifications
\[
\underline{m}_{1+2}=\underline{m}_{1\cdot 1+2}=\underline{m}_{1+1\cdot 2},
\]
and similarly for $\overline{m}_{1+2}$.
\medskip

\noindent
{\bf Expression of the optimization function $O$:} By computing, for instance using Mathematica, the elements of the sets $\Lambda_i$ of specific solutions of the equation \eqref{LAGRANGE} associated with $\alpha_a$, one obtains the following expression of the corresponding optimizing function $m\mapsto O(m)$
\begin{align*}
&\underline{O(m)}_1=\max\left\{\underline{m}_1, \tfrac{\underline{m}_{a\cdot 1+2}-\overline{m}_2}{a}, \underline{m}_{1+a\cdot 2}-a\overline{m}_2, \tfrac{a\underline{m}_{a\cdot 1+2}-\overline{m}_{1+a\cdot 2}}{a^2-1}\right\}\\
&\underline{O(m)}_2=\max\left\{\underline{m}_2, \tfrac{\underline{m}_{1+a\cdot 2}-\overline{m}_1}{a}, \underline{m}_{a\cdot 1+2}-a\overline{m}_1, \tfrac{a\underline{m}_{1+a\cdot 2}-\overline{m}_{a\cdot 1+2}}{a^2-1}\right\}\\
&\underline{O(m)}_{a\cdot 1+2}=\max\left\{\underline{m}_{a\cdot 1+2}, a\underline{m}_1+\underline{m}_2, \tfrac{(a^2-1)\underline{m}_1+\underline{m}_{1+a\cdot 2}}{a}, a\underline{m}_{1+a\cdot 2}-(a^2-1)\overline{m}_2\right\}\\
&\underline{O(m)}_{1+a\cdot 2}=\max\left\{\underline{m}_{1+a\cdot 2}, a\underline{m}_2+\underline{m}_1, \tfrac{(a^2-1)\underline{m}_2+\underline{m}_{a\cdot 1+2}}{a},a\underline{m}_{a\cdot 1+2}-(a^2-1)\overline{m}_1\right\}
\end{align*}
and the coordinates $\overline{O(m)}_i$ are obtained by replacing $\max$ by $\min$, $\underline{m}_{j}$ by $\overline{m}_{j}$ and {\sl vice-versa}.

In the case $a=1$ (canonical matrix), repetitions of identical terms in the resulting expressions, and terms of the form $\tfrac00$, both should be disregarded.
\medskip

\noindent
{\bf Optimized constraint matrices associated with atoms involved in Problem \ref{CONDPRO1}}: The atoms $A_{001}$ and $A_{011}$ (see Fig.\ \ref{PhasePortraitG2}) can be respectively written as $P_{m_{001}}^{\alpha_a}\cap \left\{x\in\R^2\ :\ \tfrac12<x_1+x_2\right\}$ and $P_{m_{011}}^{\alpha_a}$, where 
\[
m_{001}=\left(
\begin{array}{cc}
0&\tfrac12\\
0&\tfrac12\\
\tfrac12&\tfrac{1+a}2\\
\tfrac12&\tfrac{1+a}2
\end{array}\right)\quad\text{and}\quad
m_{011}=\left(
\begin{array}{cccc}
0&\tfrac12\\
\tfrac12&1\\
\tfrac{1}2&1+\tfrac{a}2\\
\tfrac{a}2&\tfrac12+a
\end{array}\right)
\]

\noindent
From now on, we assume that $m=O(m)$ is an optimized constraint matrix of a polytope $P_m^{\alpha_a}$ which satisfies Conditioning Problem \ref{CONDPRO1}.
\medskip

\noindent
{\bf Restricting the range of constraint matrices entries:} 
\begin{Claim}
The localisation condition in Problem \ref{CONDPRO1} implies that we must have 
\[
0\leq \underline{m}_1,\quad\overline{m}_1,\underline{m}_2\leq \tfrac12\leq \overline{m}_2\leq 1\quad\text{and}\quad \tfrac12\leq \underline{m}_{a\cdot 1+2}, \underline{m}_{1+a\cdot 2}.
\]
\label{BASICm21}
\end{Claim}
\noindent
{\sl Proof.} The inequalities $\underline{m}_2\leq \tfrac12\leq \overline{m}_2$ follow from the fact that $m \cap m_{001}\neq \emptyset$ and $m \cap m_{011}\neq\emptyset$. The inequality $\overline{m}_1\leq \tfrac12$ is proved by contradiction. If, otherwise, we had $\overline{m}_1>\tfrac12$, then the opimization assumption and statement {\sl (ii)} of Lemma \ref{OPTIVEC} would imply that we would have $P_m^{\alpha_a}\cap (A_{101}\cup A_{111})\neq\emptyset$ which is incompatible with the localisation condition in Problem \ref{CONDPRO1}. 

Similar arguments based on the condition $P_m^{\alpha_a}\subset (0,1)^2$ imply the inequalities
\[
0\leq \underline{m}_i\ \text{and}\ \overline{m}_i\leq 1,\ \text{for}\ i=1,2.
\]
The inequalities $\tfrac12\leq \underline{m}_{a\cdot 1+2}, \underline{m}_{1+a\cdot 2}$ are obtained in the same way. If otherwise, we had $\underline{m}_{a\cdot 1+2}<\tfrac12$ (resp.\ $\underline{m}_{1+a\cdot 2}<\tfrac12$) then, by optimality and $a\geq 1$, we would have $x_1+x_2\leq ax_1+x_2<\tfrac12$ (resp.\ $x_1+x_2\leq x_1+ax_2<\tfrac12$) for some $(x_1,x_2)\in P_m^{\alpha_a}$; hence $P_m^{\alpha_a}\cap A_{000}\neq \emptyset$ which is also incompatible with the localisation condition. \hfill $\Box$

In addition to these localisation-induced restrictions, anticipated considerations about restrictions resulting from the constraints on dynamics in Problem \ref{CONDPRO1} substantially shorten and simplify the analysis to follow. The following statement is particularly useful.
\begin{Claim}
We must have $\overline{m}_{a\cdot 1+2}\leq a\overline{m}_1+\tfrac12$.
\label{ma12}
\end{Claim} 
\noindent
{\sl Proof.} By contradiction, assume $\overline{m}_{a\cdot 1+2}> a\overline{m}_1+\tfrac12$ and suppose that we have shown that this implies the inequality 
\begin{equation}
\overline{m}_{1+a \cdot 2} > \overline{m}_1+\tfrac{a}{2}.
\label{m1a2}
\end{equation}
Then, a simplification of the expression $\overline{O(m\cap m_{011})}_1$, based on the assumption $m=O(m)$ (see Section \ref{SM-011} in Supplementary Material), yields the following conclusion
\[
\overline{O(m\cap m_{011})}_1=\min\left\{\overline{m}_1, \tfrac{\overline{m}_{a\cdot 1+2}-\tfrac12}{a}, \overline{m}_{1+a\cdot 2}-\tfrac{a}2\right\}=\overline{m}_1.
\]
However, when applied to this entry, the dynamics-induced condition $2(1-\epsilon)O(m\cap m_{011})+2\epsilon B_{\varrho,011}- {0\choose 1}\subset m$ writes
\[
2(1-\epsilon)\overline{O(m\cap m_{011})}_1=2(1-\epsilon)\overline{m}_1\leq \overline{m}_1,
\]
and the last inequality is impossible given that $2(1-\epsilon)>1$ and $\overline{m}_1>0$.

In order to complete the proof, it remains to show that $\overline{m}_{a\cdot 1+2}> a\overline{m}_1+\tfrac12$ implies \eqref{m1a2}. By contradiction again, suppose that we have $\overline{m}_{1+a \cdot 2} \leq \overline{m}_1+\tfrac{a}{2}$. Then the optimization assumption $m=O(m)$ implies the first inequality below
\[
\overline{m}_{a \cdot 1+2} \leq \frac{\overline{m}_{1+a \cdot 2}+(a^2-1)\overline{m}_1}{a} \leq \frac{\overline{m}_1+\frac{a}{2}+(a^2-1)\overline{m}_1}{a}=a\overline{m}_1+\tfrac{1}{2},
\] 
and the last inequality contradicts the assumption $\overline{m}_{a\cdot 1+2}> a\overline{m}_1+\tfrac12$. \hfill $\Box$
\medskip

\noindent
{\bf Constraint matrices associated with atomic restrictions:} Together with the expression of the optimization function, the inequalities in Claim \ref{BASICm21} and \ref{ma12} imply the following restrictions on the other entries of $m$
\[
\overline{m}_{a\cdot 1+2}\leq \tfrac{1+a}2\quad\text{and}\quad
\overline{m}_{1+a\cdot 2}\leq \tfrac12+a.
\]
Using all those inequalities here and above, the expressions of the constraint matrices $m\cap m_{001}$ and $m\cap m_{011}$ associated with atomic restrictions simplify as follows
\[
m\cap m_{001}=\left(
\begin{array}{cccccc}
\overline{m}_1&\tfrac12&\overline{m}_{a\cdot 1+2}&\min\{\overline{m}_{1+a\cdot 2},\tfrac{1+a}2\}\\
\underline{m}_1&\underline{m}_2&\underline{m}_{a\cdot 1+2}&\underline{m}_{1+a\cdot 2}
\end{array}\right)^T,
\]
and
\[
m\cap m_{011}=\left(
\begin{array}{cccccc}
\overline{m}_1&\overline{m}_2&\overline{m}_{a\cdot 1+2}&\overline{m}_{1+a\cdot 2}\\
\underline{m}_1&\tfrac12&\underline{m}_{a\cdot 1+2}&\max\{\underline{m}_{1+a\cdot 2},\tfrac{a}2\}
\end{array}\right)^T.
\]
Notice that, in addition to $\underline{m}_i<\overline{m}_i$ for all $i$, necessary conditions for the corresponding polytopes to be non-empty yield the following additional requirements 
\[
\underline{m}_{1+a\cdot 2}<\tfrac{1+a}2\quad\text{and}\quad \tfrac{a}2< \overline{m}_{1+a\cdot 2}.
\]
\medskip

\noindent
{\bf Expressions of the optimized constraint matrix $O(m\cap m_{001})$:} The inequalities above on the entries of $m$ can be used to simplify the expression of the optimized matrix $O(m\cap m_{001})$, see Section \ref{SM-001} of Supplementary Material. This yields
\begin{flalign*}
&\underline{O(m\cap m_{001})}_1=\max\left\{\underline{m}_1, \tfrac{\underline{m}_{a\cdot 1+2}-\tfrac12}{a}, \underline{m}_{1+a\cdot 2}-\tfrac{a}2\right\}\quad\text{and}\quad \overline{O(m\cap m_{001})}_1=\overline{m}_1\\
&\underline{O(m\cap m_{001})}_2=\underline{m}_2\quad\text{and}\quad\overline{O(m\cap m_{001})}_2=\tfrac12\\
&\underline{O(m\cap m_{001})}_{a\cdot 1+2}=\max\left\{\underline{m}_{a\cdot 1+2}, a\underline{m}_{1+a\cdot 2}-\tfrac{a^2-1}2\right\}\quad\text{and}\quad\overline{O(m\cap m_{001})}_{a\cdot 1+2}=\overline{m}_{a\cdot 1+2}\\
&\underline{O(m\cap m_{001})}_{1+a\cdot 2}=\underline{m}_{1+a\cdot 2}\quad\text{and}\quad\overline{O(m\cap m_{001})}_{1+a\cdot 2}=\min\left\{\overline{m}_{1+a\cdot 2}, \tfrac{a}2+\overline{m}_1, \tfrac{\tfrac{a^2-1}2+\overline{m}_{a\cdot 1+2}}{a}\right\}
\end{flalign*}
\medskip

\noindent
{\bf (In)equalities induced by the conditions on dynamics:} The transposition-induced transformation $\sigma_{321}$ on constraint matrices is $\alpha_a$-compatible and writes (denoted by the same symbol)
\[
\sigma_{321}(m)=\left(
\begin{array}{cccccc}
1-\underline{m}_2&1-\underline{m}_1&1+a-\underline{m}_{1+a\cdot 2}&1+a-\underline{m}_{a\cdot 1+2}\\
1-\overline{m}_2&1-\overline{m}_1&1+a-\overline{m}_{1+a\cdot 2}&1+a-\overline{m}_{a\cdot 1+2}
\end{array}\right)^T
\] 
Using also the expression of $G_{\varrho,\epsilon}$, the values $B_{\varrho,001}=\varrho{1\choose 1}$ and $B_{\varrho,011}={0\choose 1}$, the dynamics conditions in Problem \ref{CONDPRO1} express as follows on constraint matrices (especially as the left equation is concerned)
\[
2(1-\epsilon)O(m\cap m_{001})+2\varrho\epsilon{1\choose 1}=\sigma_{321}(m)\quad\text{and}\quad 2(1-\epsilon)O(m\cap m_{011})+(2\epsilon-1){0\choose 1}\subset m.
\]
For those entries of $O(m\cap m_{001})$ that are trivial, the left equation here immediately gives some of the entries of $m$. In particular, that $\overline{O(m\cap m_{001})}_2=\tfrac12$ yields 
\[
\underline{m}_1=\epsilon(1-2\varrho),
\]
and $\overline{O(m\cap m_{001})}_1=\overline{m}_1$ together with $\underline{O(m\cap m_{001})}_2=\underline{m}_2$ yield
\[
\overline{m}_1=\underline{m}_2=r_{\varrho}\quad\text{where}\quad r_{\varrho}:=\tfrac{1-2\varrho\epsilon}{3-2\epsilon}.
\]
Likewise, considerations that imply $\overline{O(m\cap m_{001})}_{a\cdot 1+2}$ and $\underline{O(m\cap m_{001})}_{1+a\cdot 2}$ result in 
\[
\overline{m}_{a\cdot 1+2}=\underline{m}_{1+a\cdot 2}=(1+a)r_{\varrho}
\]
\medskip

\noindent
{\bf Solving the remaining equations:} Altogether, the remaining entries of $O(m\cap m_{001})$ simplify as follows
\begin{flalign*}
&\underline{O(m\cap m_{001})}_1=\max\left\{\epsilon(1-2\varrho), \tfrac{\underline{m}_{a\cdot 1+2}-\tfrac12}{a},(1+a)r_{\varrho}-\tfrac{a}2\right\}\\
&\underline{O(m\cap m_{001})}_{a\cdot 1+2}=\max\left\{\underline{m}_{a\cdot 1+2}, a(1+a) r_{\varrho}-\tfrac{a^2-1}2\right\}\\
&\overline{O(m\cap m_{001})}_{1+a\cdot 2}=\min\left\{\overline{m}_{1+a\cdot 2}, \tfrac{1+a}{a}r_{\varrho}+\tfrac{a^2-1}{2a}\right\}
\end{flalign*}
and the corresponding equations write
\begin{align*}
&2(1-\epsilon)\max\left\{\epsilon(1-2\varrho), \tfrac{\underline{m}_{a\cdot 1+2}-\tfrac12}{a},(1+a)r_{\varrho}-\tfrac{a}2\right\}+\overline{m}_2=1-2\varrho\epsilon\\
&2(1-\epsilon)\max\left\{\underline{m}_{a\cdot 1+2}, a(1+a) r_{\varrho}-\tfrac{a^2-1}2\right\}+\overline{m}_{1+a\cdot 2}=(1+a)\left(1-2\varrho\epsilon\right)\\
&2(1-\epsilon)\min\left\{\overline{m}_{1+a\cdot 2}, \tfrac{1+a}{a}r_{\varrho}+\tfrac{a^2-1}{2a}\right\}+\underline{m}_{a\cdot 1+2}=(1+a)\left(1-2\varrho\epsilon\right)
\end{align*}
The following statement is proved independently, see Section \ref{SM-ma12Bis} of Supplementary Material.
\begin{Claim}
We must have $\underline{m}_{a\cdot 1+2}>a(1+a) r_{\varrho}-\tfrac{a^2-1}2$.
\label{ma12Bis}
\end{Claim}
Accordingly, when solving the last two equations, focus can be made on the alternative in the last equation. However, we must have $\overline{m}_{1+a\cdot 2}> \tfrac{1+a}{a}r_{\varrho}+\tfrac{a^2-1}{2a}$ because otherwise the last two equations would be solved as 
\[
\underline{m}_{a\cdot 1+2}=\overline{m}_{1+a\cdot 2}=(1+a)r_{\varrho}
\]
and the same values of $\overline{m}_{a\cdot 1+2}$ and $\underline{m}_{1+a\cdot 2}$ above would imply that $P_m^{\alpha_a}$ would be empty. Using the inequality $\overline{m}_{1+a\cdot 2}> \tfrac{1+a}{a}r_{\varrho}+\tfrac{a^2-1}{2a}$ and Claim \ref{ma12Bis}, one can solve the last two equations as
\[
\underline{m}_{a\cdot 1+2}=\tfrac{1+a}{a}\left(a\epsilon(1-2\varrho)+(1-\epsilon)(1-2r_\varrho)\right)=\tfrac{1+a}{a}\left(r_\varrho+(a-1)\epsilon(1-2\varrho)\right)
\]
and
\[
\overline{m}_{1+a\cdot 2}=\tfrac{1+a}{a}\left(a r_\varrho+2(a-1)(1-\epsilon)^2(1-2r_\varrho)\right).
\]
This expression of $\underline{m}_{a\cdot 1+2}$ implies that Claim \ref{ma12Bis} is equivalent to $\epsilon> 1-\tfrac{a}2$.

It remains to compute the value of $\overline{m}_2$ using the first equation. Notice that the expression of $\underline{m}_{a\cdot 1+2}$ shows that 
\[
\max\left\{\epsilon(1-2\varrho), \tfrac{\underline{m}_{a\cdot 1+2}-\tfrac12}{a}, (1+a)r_{\varrho}-\tfrac{a}2\right\}=\epsilon(1-2\varrho)\ \text{iff}\ \epsilon\geq 1-\tfrac{a}2.
\]
Therefore, we have $\overline{m}_2:=\overline{m_a}_2=1-2\epsilon\left(1-\varrho-\epsilon(1-2\varrho)\right)$. Altogether, we have shown that Problem \ref{CONDPRO1} has at most a unique solution $m_a$ (which also depends on $\varrho$ and $\epsilon$) given by 
{\small\begin{equation}
m_a:=\left(
\begin{array}{cccc}
r_{\varrho}&1-2\epsilon\left(1-\varrho-\epsilon(1-2\varrho)\right)&(1+a)r_{\varrho}&\tfrac{1+a}{a}\left(a r_\varrho+2(a-1)(1-\epsilon)^2(1-2r_\varrho)\right)\\
\epsilon(1-2\varrho)&r_{\varrho}&\tfrac{1+a}{a}\left(r_\varrho+(a-1)\epsilon(1-2\varrho)\right)&(1+a)r_{\varrho}
\end{array}\right)^T
\label{ma}
\end{equation}}
\medskip

\noindent
{\bf Proof of Statement \ref{NOSOLUTIONP1}:} This is immediate, because for $a=1$, expression \eqref{ma} yields
\[
\underline{m_a}_{a\cdot 1+2}=\overline{m_a}_{a\cdot 1+2}\quad\left(\text{and}\quad \underline{m_a}_{1+a\cdot 2}=\overline{m_a}_{1+a\cdot 2}\right).
\]
which both imply $P_{m_a}^{\alpha_a}=\emptyset$. \hfill $\Box$
\medskip

\noindent
{\bf Additional restrictions on the parameters that result from Problem \ref{CONDPRO1}:} In order to ensure the existence claim in Statement \ref{SOLP1}, one needs to verify that, in addition to the optimality condition $m_a=O(m_a)$, the entries of $m_a$ satisfy both the conditions that have been exhibited above in this section and the remaining ones in the conditioning problem. To that goal, we assume that $\rho\in (0,\tfrac12)$, $a>1$ and $\epsilon\in (\max\{0,1-\tfrac{a}2\},\tfrac12)$ and we first establish those extra conditions on these parameters that result from this verification. The conditions in Statement \ref{SOLP1} will follow from an analysis of all conditions combined together. 

Checking the conditions related to $P_{m_a}^{\alpha_a}\neq \emptyset$ and to localisation yields the following conclusions.
\begin{itemize}
\item The only additional restriction resulting from the inequalities in Claim \ref{BASICm21} and \ref{ma12} and the inequalities $\underline{m_a}_i<\overline{m_a}_i$ comes from $\tfrac12\leq \underline{m_a}_{a\cdot 1+2}$ and writes\footnote{We always have $\tfrac14<r_\varrho<\tfrac12$.}
\begin{equation*}
\tfrac{a}{2(1+a)}-r_\varrho\leq (a-1)\epsilon(1-2\varrho)
\end{equation*}
This inequality is actually weaker than the inequality \eqref{COND3} below. Moreover, the inequalities $\overline{m_a}_{a\cdot 1+2}\leq \tfrac{1+a}2$ and $\overline{m_a}_{1+a\cdot 2}\leq \tfrac12+a$ will follow from the ones here, once optimality is ensured.
\item The only additional restriction resulting from the localisation conditions $\underline{m_a}_{1+a\cdot 2}<\tfrac{1+a}2$ and $\tfrac{a}2< \overline{m_a}_{1+a\cdot 2}$ comes from the second inequality and writes
\begin{align}
\tfrac{a}{2(1+a)}- r_\varrho<2\tfrac{a-1}{a}(1-\epsilon)^2(1-2r_\varrho)
\label{COND2}
\end{align}
\item It remains to make sure that $P_{m_a}^{\alpha_a}\cap A_{000}= \emptyset$, {\sl viz.}\ $P_{m_a}^{\alpha_a}\cap\left\{x\in\R^2\ :\ x_1+x_2<\tfrac12\right\}=\emptyset$. This condition can be regarded as $\tfrac12\leq\underline{m_a}_{1+2}=\underline{O(m_a)}_{1+2}$ when adding one coefficient to the matrix $\alpha$ so that constraints on the sum $x_1+x_2$ are also imposed. Explicit computation of the corresponding optimization function then yields the inequality
\begin{align*}
\tfrac12\leq &\max\left\{\underline{m_a}_1+\underline{m_a}_2,\underline{m_a}_{a\cdot 1+2}-(a-1)\overline{m_a}_1,\underline{m_a}_{1+a\cdot 2}-(a-1)\overline{m_a}_2,\tfrac{(a-1)\underline{m_a}_1+\underline{m_a}_{1+a\cdot 2}}{a},\right.\\
&\left.\tfrac{(a-1)\underline{m_a}_2+\underline{m_a}_{a\cdot 1+2}}{a}, \tfrac{\underline{m_a}_{a\cdot 1+2}+\underline{m_a}_{1+a\cdot 2}}{a+1}\right\}
\end{align*}
The expression of the entries of $m_a$ imply that this inequality simplifies as follows
\[
\tfrac12\leq \max\left\{(1+a)r_{\varrho}-(a-1)(1-2\epsilon\left(1-\varrho-\epsilon(1-2\varrho)\right)),\tfrac{1+a}{a}r_\varrho+\tfrac{a-1}{a}\epsilon(1-2\varrho)\right\}.
\]
Elementary algebra then shows that this is equivalent to 
\begin{equation}
\tfrac{a}{2(1+a)}-r_\varrho\leq \tfrac{a-1}{1+a}\epsilon(1-2\varrho)
\label{COND3}
\end{equation}
\end{itemize}

In order to check those inequalities that are associated with the dynamical conditioning $2(1-\epsilon)O(m_a\cap m_{011})+(2\epsilon-1){0\choose 1}\subset m_a$, it is convenient to proceed similarly as before to simplify the entries of $O(m\cap m_{011})$ (Section \ref{SM-011} in Supplementary Material), prior to evaluate the resulting expressions for $O(m_a\cap m_{011})$. Using the restriction $\epsilon\in (\max\{0,1-\tfrac{a}2\},\tfrac12)$, one then obtains the following expressions
\begin{flalign*}
&\underline{O(m_a\cap m_{011})}_1=\underline{m_a}_1&\text{and}&& \overline{O(m_a\cap m_{011})}_1=\tfrac{\overline{m_a}_{a\cdot 1+2}-\tfrac12}{a}&\\
&\underline{O(m_a\cap m_{011})}_2=\tfrac12&\text{and}&& \overline{O(m_a\cap m_{011})}_2=\overline{m_a}_2&\\
&\underline{O(m_a\cap m_{011})}_{a\cdot 1+2}=a\underline{m_a}_1+\tfrac12 &\text{and}&& \overline{O(m_a\cap m_{011})}_{a\cdot 1+2}=\overline{m_a}_{a\cdot 1+2}&\\
&\underline{O(m_a\cap m_{011})}_{1+a\cdot 2}=\tfrac{a}2+\underline{m_a}_1&\text{and}&&\overline{O(m_a\cap m_{011})}_{1+a\cdot 2}=\overline{m_a}_{1+a\cdot 2}&
\end{flalign*}
Accordingly, for the entries $\underline{O(m_a\cap m_{011})}_1$ and $\overline{O(m_a\cap m_{011})}_2$, the condition $2(1-\epsilon)O(m_a\cap m_{011})+(2\epsilon-1){0\choose 1}\subset m_a$ respectively implies the redundant inequalities
\[
0\leq \underline{m_a}_1\quad \text{and}\quad \overline{m_ a}_2\leq 1.
\]
Combining the conditions that result for the entries $\underline{O(m_a\cap m_{011})}_2$ and $\overline{O(m_a\cap m_{011})}_{a\cdot 1+2}$ yield
\begin{equation}
r_\varrho\leq \min\left\{\epsilon,\tfrac1{1+a}\right\}.
\label{COND4}
\end{equation}
For the entry $\overline{O(m_a\cap m_{011})}_{1+a\cdot 2}$, one gets $\overline{m_a}_{1+a\cdot 2}\leq a$ which rewrites as 
\begin{equation}
2\tfrac{a-1}{a}(1-\epsilon)^2(1-2r_\varrho)\leq \tfrac{a}{1+a}-r_\varrho
\label{COND5}
\end{equation}
Now evaluating the condition for the entry $\overline{O(m_a\cap m_{011})}_1$ yields
\begin{align}
&2(1-\epsilon)(1+a)r_\varrho\leq a r_\varrho+1-\epsilon
\label{COND6}
\end{align}
Finally, for the entries $\underline{O(m_a\cap m_{011})}_{a\cdot 1+2}$ and $\underline{O(m_a\cap m_{011})}_{1+a\cdot 2}$, we respectively get the inequalities
\begin{align}
&(1+a)r_\varrho \leq (1+(1-2\epsilon)a^2)\epsilon(1-2\varrho)+a\epsilon\quad\text{and}\quad (1+a)r_\varrho\leq 2(1-\epsilon)\epsilon(1-2\varrho)+a\epsilon
\label{COND7}
\end{align}
The condition $a>1$ implies that the right inequality is stronger than the left one.
\medskip

\noindent
{\bf Proof of the existence conditions in Statement \ref{SOLP1}:} In order to determine the domains of parameters for which all conditions \eqref{COND2} - \eqref{COND7} hold, we consider separately the cases $\rho\in \left(0,\tfrac13\right]$ and $\rho\in (\tfrac13,\tfrac12)$ because $\epsilon\mapsto r_\varrho(\epsilon)$ is non-decreasing (resp.\ decreasing) in the first (resp.\ second) case. Notice that 
\[
r_\varrho(0^+)=\tfrac13\quad\text{and}\quad r_\varrho(\tfrac12^-)=\tfrac{1-\varrho}2.
\]
{\sl Analysis for $\rho\in \left(0,\tfrac13\right]$:}
\begin{itemize}
\item Equation \eqref{COND2}: Both sides are non-increasing with $\epsilon$; hence it suffices to verify that the value of the LHS for $\epsilon=0^+$ is not larger than the value of the RHS for $\epsilon=\tfrac12^-$, which requires $a\geq \tfrac12$.
\item Equation \eqref{COND3}: Monotonicity implies that it suffices to check this condition for $\epsilon=0^+$, which imposes $a\leq 2$.
\item Equation \eqref{COND4}: The map $\epsilon\mapsto r_\varrho$ is concave; hence if \eqref{COND4} holds, then it does for $\epsilon$ in a left neighborhood of $\tfrac12$. Moreover the condition $a>1$ imposes that it holds iff $r_\varrho(\tfrac12^-) \leq\tfrac1{1+a}$, {\sl viz.} $a\leq \tfrac{1+\varrho}{1-\varrho}$ (and we have $\tfrac{1+\varrho}{1-\varrho}\leq 2$).
\item Equation \eqref{COND5}: By continuity, \eqref{COND5} holds for $\epsilon$ in a left neighborhood of $\tfrac12$ iff $\tfrac{a-1}{2a}\varrho<\tfrac{a}{1+a}-\tfrac{1-\varrho}2$, which holds under our conditions on $a$ and $\varrho$.
\item Similar considerations apply for equations \eqref{COND6} and \eqref{COND7}.
\end{itemize}
Altogether, it follows that, for every $\rho\in \left(0,\tfrac13\right]$ and every $a\in \left(1,\tfrac{1+\varrho}{1-\varrho}\right]$, all conditions \eqref{COND2} - \eqref{COND7} hold iff $\epsilon\in \left[\epsilon_{\varrho,a},\tfrac12\right)$ where $\epsilon_{\varrho,a}$ is the maximum of the roots in $(0,\tfrac12)$ of the polynomials associated with the inequalities \eqref{COND5} - \eqref{COND7}. 

For the uniform distribution $\varrho=\tfrac13$, one checks that $\epsilon_{\frac13,a}$ is actually given by the right inequality in \eqref{COND7}; hence the value $\epsilon_{\frac13,a}=\frac{3a+2-\sqrt{9a^2+4a-4}}4$. By continuity (and anticipating the conclusion of the analysis below), the property holds for $\varrho$ close to $\tfrac13$. Assuming that $\epsilon_{\varrho,a}=\epsilon_{\frac13,a}+\delta_a(\varrho-\frac13)+O\left((\varrho-\tfrac13)^2\right)$ and solving the right inequality in \eqref{COND7} to first order in $\varrho-\tfrac13$ yields a cumbersome value of $\delta_a$, which is negative for $a=1$ and which decreases when $a$ increases.
\medskip

\noindent
{\sl Analysis for $\rho\in \left(\tfrac13,\tfrac12\right)$:}
\begin{itemize}
\item Equation \eqref{COND2}: The LHS increases with $\epsilon$ and the RHS decreases; hence we must have $\tfrac{a}{2(1+a)}-\tfrac{1-\varrho}2\leq \tfrac{a-1}{2a}\varrho$, which holds under our conditions on $a$ and $\varrho$.
\item Equation \eqref{COND3}: By continuity, \eqref{COND3} holds for $\epsilon$ in a left neighbourhood of $\tfrac12$ if $a-(1+a)(1-\varrho)<(a-1)(1-2\varrho)$, which is equivalent to $a<\tfrac{\varrho}{3\varrho-1}$. Moreover, anticipating the condition $a\leq 2$ below, notice that $\tfrac{\varrho}{3\varrho-1}\leq 2$ iff $\varrho\geq\tfrac25$. When $\varrho$ satisfies this condition, one checks that the derivative wrt $\epsilon$ of the LHS in \eqref{COND3} at $\epsilon=\tfrac12$ is larger than the derivative of the RHS; hence \eqref{COND3} also holds for $\epsilon$ in a left neighbourhood of $\tfrac12$ when $a=\tfrac{\varrho}{3\varrho-1}$.
\item Equation \eqref{COND4}: That the LHS decreases with $\epsilon$ and the condition $a>1$ impose that \eqref{COND4} holds for $\epsilon$ in a left neighbourhood of $\tfrac12$ iff $r_\varrho(\tfrac12^-) <\tfrac1{1+a}$, {\sl viz.} $a< \tfrac{1+\varrho}{1-\varrho}$.
\item Equation \eqref{COND5}: The LHS decreases with $\epsilon$ and the RHS increases; hence it suffices to check this condition for $\epsilon=0^+$, which imposes $a\leq 2$ (and we have $2<\tfrac{1+\varrho}{1-\varrho}$).
\item Strict inequalities at $\epsilon=\tfrac12$ (which holds under our conditions on $a$ and $\varrho$) and continuity arguments imply that equations \eqref{COND6} and \eqref{COND7} hold for $\epsilon$ in a left neighbourhood of $\tfrac12$. 
\end{itemize}
Similarly as before, it follows that, for every $\rho\in \left(\tfrac13,\tfrac12\right)$ and every $a\in \left(1,\min\left\{2,\tfrac{\varrho}{3\varrho-1}\right\}\right]$, all conditions \eqref{COND2} - \eqref{COND7} hold iff $\epsilon\in \left[\epsilon_{\varrho,a},\tfrac12\right)$.

Finally, it remains to verify that $m_a$ is indeed an optimized constraint matrix, {\sl ie.} we have $m_a=O(m_a)$. This is done in Section 1.3 of Supplementary Material, where the edges of $P^{\alpha_a}_{m_a}$ are also identified. \hfill $\Box$
\medskip

\noindent
{\bf Proof that $P^{\alpha_a}_{m_a}\cup \sigma_{321}(P^{\alpha_a}_{m_a})$ is an AsIUP with respect to $\Sigma=1-\text{Id}$.} That $P^{\alpha_a}_{m_a}$ satisfies the conditions of Problem \ref{CONDPRO1} makes $P^{\alpha_a}_{m_a}\cup \sigma_{321}(P^{\alpha_a}_{m_a})$ an IUP. Moreover we have
\begin{itemize}
\item $\overline{m_a}_1<\frac12<\underline{\Sigma(m_a)}_1=1-\overline{m_a}_1$; hence $P^{\alpha_a}_{m_a}\cap P^{\alpha_a}_{\Sigma(m_a)}=\emptyset$ where $P^{\alpha_a}_{\Sigma(m_a)}:=\Sigma(P^{\alpha_a}_{m_a})$,
\item $\overline{\sigma_{321}(m_a)}_1=1-\underline{m_a}_2=\underline{\Sigma(m_a)}_1$; hence $\sigma_{321}(P^{\alpha_a}_{m_a})\cap \Sigma(P^{\alpha_a}_{m_a})=\emptyset$,
\end{itemize}
which, together with the commutation of $\sigma_{321}$ and $\Sigma$ imply that the IUP does not intersect its image under $\Sigma$, and hence it is actually an AsIUP. Statement \ref{SOLP1} is proved. \hfill $\Box$

\section{Analysis of Conditioning Problem \ref{CONDPRO2}}\label{A-P2}
In principle, Conditioning Problem \ref{CONDPRO2} could be analysed following the plan in Section \ref{S-ANALYSIS}, as in the previous Appendix for Conditioning Problem \ref{CONDPRO1}. However, for simplicity, we apply instead a continuation argument to the solution of Problem \ref{CONDPRO1}. Indeed, as shown in Section \ref{SM-OPTI} of Supplementary Material, the solution $P_{m_a}^{\alpha_a}$ is a quadrilateral with edges 
\[
x_1=\underline{m_a}_1,\ x_1+ax_2=\underline{m_a}_{1+a\cdot 2},\ x_2=\overline{m_a}_2,\ \text{and}\quad ax_1+x_2=\overline{m_a}_{a\cdot 1+2},
\]
which are determined by the condition $G_{\varrho,\epsilon}(P\cap A_{001})=\sigma_{321}(P)$. More precisely, we have
\begin{itemize}
\item the edge $x_1=\underline{m_a}_1$ is obtained as the image of the line $x_1=\tfrac12$, viewed as belonging to $A_{111}$,
\item the edge $x_2=\overline{m_a}_2$ is obtained as the image, of the image (located in $A_{111}$), of the line $x_2=\tfrac12$ viewed as belonging to $A_{001}$,
\item the pair of edges $x_1+ax_2=\underline{m_a}_{1+a\cdot 2}$ and $ax_1+x_2=\overline{m_a}_{a\cdot 1+2}$ are obtained as solving a pair of equations induced by the condition $G_{\varrho,\epsilon}(P\cap A_{001})=\sigma_{321}(P)$.
\end{itemize}

For Conditioning Problem \ref{CONDPRO2}, these conditions can be duplicated as follows. For $j=1,2$, let $P_{m_j}^{\alpha_a}$ be a quadrilateral with edges
\[
x_1=\underline{m_j}_1,\ x_1+ax_2=\underline{m_j}_{1+a\cdot 2},\ x_2=\overline{m_j}_2,\ \text{and}\quad ax_1+x_2=\overline{m_j}_{a\cdot 1+2},
\]
which are determined by the conditions $G_{(\rho_i)_{i=1}^3,\epsilon}(P_1\cap A_{001})=P_2$, $G_{(\rho_i)_{i=1}^3,\epsilon}(P_2\cap A_{111})=P_1$. More precisely, we have
\begin{itemize}
\item the edge $x_1=\underline{m_1}_1$ is obtained as the image of the line $x_1=\tfrac12$, viewed as belonging to $A_{111}$, The edge $x_1=\underline{m_2}_1$ is obtained as the image of $x_1=\underline{m_1}_1$ (located in $A_{001}$).
\item the edge $x_2=\overline{m_2}_2$ is obtained as the image of the line $x_2=\tfrac12$, viewed as belonging to $A_{001}$. The edge $x_2=\overline{m_1}_2$ is obtained as the image of $x_2=\overline{m_2}_2$ (located in $A_{111}$) .
\item the pair $x_1+ax_2=\underline{m_1}_{1+a\cdot 2}$ and $x_1+ax_2=\underline{m_2}_{1+a\cdot 2}$ (resp.\ $ax_1+x_2=\overline{m_1}_{a\cdot 1+2}$ and $ax_1+x_2=\overline{m_2}_{a\cdot 1+2}$) is obtained by imposing that these edges must be exchanged under the dynamics, and that the first one belongs to $A_{001}$ and the second one belongs to $A_{111}$. 
\end{itemize}
\medskip

\noindent
{\bf Proof of Statement \ref{SOLP2}:} The entry pairs $(\underline{m_j}_1)_{j=1}^2$ and $(\overline{m_j}_2)_{j=1}^2$ follow immediately from the conditions above.  Moreover, both pairs $(\underline{m_j}_{1+a\cdot 2})_{j=1}^2$ and $(\overline{m_j}_{a\cdot 1+2})_{j=1}^2$ are solutions of pairs of linear equations whose discriminant is $(1-2\epsilon)(3+2\epsilon)$. Hence they are also uniquely determined. The remaining entries of $m_1$ and $m_2$ are obtained using linear combinations in the definition of the optimization function $O$, so that both matrices are optimized matrices. 

The map $G_{(\rho_i)_{i=1}^3,\epsilon}$ continuously depends on $(\rho_i)_{i=1}^3$; hence so do the entries of $m_1$ and $m_2$. For $\delta=0$, we have $m_1(0)=m_a$ and $m_2(0)=\sigma_{321}(m_a)$ by optimality; hence for $\epsilon\in (\epsilon_{\varrho,a},\tfrac12)$, the entries of $m_1(0)$ satisfy the inequalities \eqref{COND2} - \eqref{COND7} above associated with existence in Problem \ref{CONDPRO1} and all these inequalities are strict. Similarly, the entries of $m_2(0)$ satisfy all strict inequalities associated with $\sigma_{321}(m_a)$ and induced by the conditioning \ref{CONDPRO1}. 

By continuity, the same inequalities must hold for the matrices $m_1(\delta)$ and $m_2(\delta)$ when $\delta\neq 0$ is small enough. As shown above for the Problem \ref{CONDPRO1}, it is easy to check that under these conditions, the solution pair $(m_1,m_2)$ is unique. Existence and uniqueness of the solution of the Conditioning Problem \ref{CONDPRO2} immediately follow (and by optimality, the continued polytopes $P_{m_1(\delta)}^{\alpha_a}$ and $P_{m_2(\delta)}^{\alpha_a}$ must be quadrilaterals defined by the same types of edges as $P_{m_a}^{\alpha_a}$ and $\sigma_{321}(P_{m_a}^{\alpha_a})$). 
\medskip

\noindent
{\bf Proof that $P_{m_1(\delta)}^{\alpha_a}\cup P_{m_2(\delta)}^{\alpha_a}$ is an AsIUP with respect to $\Sigma=1-\text{Id}$.} Similarly as in the end of Appendix \ref{A-P1}, it suffices to show that we have 
\[
\overline{m_1(\delta)}_1<1-\overline{m_1(\delta)}_1\quad\text{and}\quad \overline{m_2(\delta)}_1=1-\overline{m_1(\delta)}_1,
\]
because these conditions imply $P_{m_1(\delta)}^{\alpha_a}\cap\Sigma(P_{m_1(\delta)}^{\alpha_a})=\emptyset$ and  $P_{m_2(\delta)}^{\alpha_a}\cap\Sigma(P_{m_1(\delta)}^{\alpha_a})=\emptyset$. 

The inequality above holds for $\delta=0$ and hence for $\delta>0$ small by continuity. In order to show the equality, we observe that the knowledge of the polytope edges imply that we must have (the dependence on $\delta$ is not mentioned for clarity)
\[
\overline{m_1}_1=\tfrac{a \overline{m_1}_{a\cdot 1+2}-\underline{m_1}_{1+a\cdot 2}}{a^2-1}
\quad\text{and}\quad 
\overline{m_2}_1=\tfrac{a \overline{m_2}_{a\cdot 1+2}-\underline{m_2}_{1+a\cdot 2}}{a^2-1}
\]
which, together with the equations of the edges
\[
\left\{\begin{array}{l}
2(1-\epsilon)\underline{m_1}_{1+a\cdot 2}+2\epsilon(\varrho-\delta+a(\varrho+\delta))=\underline{m_2}_{1+a\cdot 2}\\
2(1-\epsilon)\underline{m_2}_{1+a\cdot 2}+2\epsilon(1-\varrho+\delta+a(1-\varrho-\delta))-(1+a)=\underline{m_1}_{1+a\cdot 2}
\end{array}\right.
\]
\[
\left\{\begin{array}{l}
2(1-\epsilon)\overline{m_1}_{a\cdot 1+2}+2\epsilon(a(\varrho-\delta)+\varrho+\delta)=\overline{m_2}_{a\cdot 1+2}\\
2(1-\epsilon)\overline{m_2}_{a\cdot 1+2}+2\epsilon(a(1-\varrho+\delta)+1-\varrho-\delta)-(1+a)=\overline{m_2}_{a\cdot 1+2}
\end{array}\right.
\]
implies that $\overline{m_2(\delta)}_1+\overline{m_1(\delta)}_1$ does not depend on $\delta$ (nor on $\varrho$), and hence it must be equal to 1 for all $\delta\geq 0$. Statement \ref{SOLP2} is proved. \hfill $\Box$

\section{Analysis of Conditioning Problem \ref{CONDPRO3}}\label{A-P3}
As for Conditioning Problem \ref{CONDPRO1}, the analysis in this Appendix follows the plan described in Section \ref{S-ANALYSIS}. 
For Conditioning Problem \ref{CONDPRO3}, it suffices to consider the canonical coefficient matrix $\alpha$ associated with $G_{3,\epsilon}$. Therefore, we may use the adapted notations introduced in Appendix \ref{A-PARTITION}. For convenience, we write the corresponding constraint vector in $\R^{6\times 2}$ using the following representation
\[
m=\left(
\begin{array}{cccccc}
\overline{m}_1&\overline{m}_2&\overline{m}_3&\overline{m}_{1+2}&\overline{m}_{2+3}&\overline{m}_{1+2+3}\\
\underline{m}_1&\underline{m}_2&\underline{m}_3&\underline{m}_{1+2}&\underline{m}_{2+3}&\underline{m}_{1+2+3}
\end{array}\right)^T
\]

\subsection{Preliminaries}
{\bf Expression of the optimization function $O$:} Proceeding as for the optimizing function in Problem \ref{CONDPRO1}, one obtains the following expression of the corresponding optimizing function $m\mapsto O(m)$
\begin{align*}
&\underline{O(m)}_1=\max\left\{\underline{m}_1,\ \underline{m}_{1+2}-\overline{m}_2,\ \underline{m}_{1+2+3}-\overline{m}_{2+3},\ \underline{m}_{1+2+3}-\overline{m}_2-\overline{m}_3,\ \underline{m}_3+\underline{m}_{1+2}-\overline{m}_{2+3}\right\}\\
&\underline{O(m)}_2=\max\left\{\underline{m}_2,\ \underline{m}_{1+2}-\overline{m}_1,\ \underline{m}_{2+3}-\overline{m}_3,\ \underline{m}_{1+2+3}-\overline{m}_1-\overline{m}_3,\ \underline{m}_{1+2}+\underline{m}_{2+3}-\overline{m}_{1+2+3}\right\}\\
&\underline{O(m)}_3=\max\left\{\underline{m}_3,\ \underline{m}_{2+3}-\overline{m}_2,\ \underline{m}_{1+2+3}-\overline{m}_{1+2},\ \underline{m}_{1+2+3}-\overline{m}_1-\overline{m}_2,\ \underline{m}_1+\underline{m}_{2+3}-\overline{m}_{1+2}\right\}\\
&\underline{O(m)}_{1+2}=\max\left\{\underline{m}_{1+2},\ \underline{m}_1+\underline{m}_2,\ \underline{m}_{1+2+3}-\overline{m}_3,\ \underline{m}_1+\underline{m}_{2+3}-\overline{m}_3,\ \underline{m}_2+\underline{m}_{1+2+3}-\overline{m}_{2+3}\right\}\\
&\underline{O(m)}_{2+3}=\max\left\{\underline{m}_{2+3},\ \underline{m}_2+\underline{m}_3,\ \underline{m}_{1+2+3}-\overline{m}_1,\ \underline{m}_3+\underline{m}_{1+2}-\overline{m}_1,\ \underline{m}_2+\underline{m}_{1+2+3}-\overline{m}_{1+2}\right\}\\
&\underline{O(m)}_{1+2+3}=\max\left\{\underline{m}_{1+2+3},\ \underline{m}_1+\underline{m}_{2+3},\ \underline{m}_3+\underline{m}_{1+2},\ \underline{m}_1+\underline{m}_2+\underline{m}_3,\ \underline{m}_{1+2}+\underline{m}_{2+3}-\overline{m}_2\right\}
\end{align*}
As before, the coordinates $\overline{O(m)}_{i+\cdots +j}$ can be obtained by replacing $\max$ by $\min$, $\underline{m}_{i'+\cdots +j'}$ by $\overline{m}_{i'+\cdots +j'}$ and {\sl vice-versa}.

\noindent
{\bf Optimized constraint vectors associated with atoms involved in Conditioning Problem} \ref{CONDPRO3}: The optimized constraint vectors $m_\omega$ associated with the involved atoms $P_{m_{\omega}}$ are
\[
m_{000101}=\left(
\begin{array}{cccccc}
\tfrac12&\tfrac12&\tfrac12&1&\tfrac12&1\\
0&0&0&\tfrac12&0&\tfrac12
\end{array}\right)^T,
\quad
m_{100101}=\left(
\begin{array}{cccccc}
1&\tfrac12&\tfrac12&\tfrac32&\tfrac12&\tfrac32\\
\tfrac12&0&0&\tfrac12&0&\tfrac12
\end{array}\right)^T,
\]
\[
m_{100111}=\left(
\begin{array}{cccccc}
1&\tfrac12&\tfrac12&\tfrac32&1&\tfrac32\\
\tfrac12&0&0&\tfrac12&\tfrac12&1
\end{array}\right)^T,
\quad
m_{100112}=\left(
\begin{array}{cccccc}
1&\tfrac12&\tfrac12&\tfrac32&1&2\\
\tfrac12&0&0&1&\tfrac12&\tfrac32
\end{array}\right)^T,
\]
\[
m_{110111}=\left(
\begin{array}{cccccc}
1&1&\tfrac12&\tfrac32&1&\tfrac32\\
\tfrac12&\tfrac12&0&1&\tfrac12&1
\end{array}\right)^T,
\quad
m_{110112}=\left(
\begin{array}{cccccc}
1&1&\tfrac12&\tfrac32&\tfrac32&2\\
\tfrac12&\tfrac12&0&1&\tfrac12&\tfrac32
\end{array}\right)^T,
\]
and
\[
m_{110212}=\left(
\begin{array}{cccccc}
1&1&\tfrac12&2&\tfrac32&\tfrac52\\
\tfrac12&\tfrac12&0&\tfrac32&\tfrac12&\tfrac32
\end{array}\right)^T.
\]
\medskip

\noindent
{\bf Consequence of the symmetry $\sigma_{4321}$ on  constraint matrices entries:} The symmetry transformation $\sigma_{4321}(x)=\left(-x_3,-x_2,-x_1\right)$ induces the following transformation on constraint matrices
\[
\sigma_{4321}(m)=\left(
\begin{array}{cccccc}
1-\underline{m}_3&1-\underline{m}_2&1-\underline{m}_1&2-\underline{m}_{2+3}&2-\underline{m}_{1+2}&3-\underline{m}_{1+2+3}\\
1-\overline{m}_3&1-\overline{m}_2&1-\overline{m}_1&2-\overline{m}_{2+3}&2-\overline{m}_{1+2}&3-\overline{m}_{1+2+3}
\end{array}\right)^T
\]
which is obviously $\alpha$-compatible. This immediately yields the following property. 
\begin{Claim}
Let $m_2=O(m_2)$ be the optimized constraint matrix of polytope $P_{m_2}$ which satisfies the symmetry $\sigma_{4321}(P_{m_2})=P_{m_2}$. Then we have
\begin{align*}
&\underline{m_2}_1+\overline{m_2}_3=\overline{m_2}_1+\underline{m_2}_3=\underline{m_2}_2+\overline{m_2}_2=1\\
&\underline{m_2}_{1+2}+\overline{m_2}_{2+3}=\overline{m_2}_{1+2}+\underline{m_2}_{2+3}=2\\
&\underline{m_2}_{1+2+3}+\overline{m_2}_{1+2+3}=3
\end{align*}
\label{SYM-M2}
\end{Claim}
\medskip

\noindent
From now on, we assume that $m_1=O(m_1)$ and $m_2=O(m_2)$ are optimized constraint matrices of polytopes $P_1:=P_{m_1}$ and $P_2:=P_{m_2}$ which satisfy Conditioning Problem \ref{CONDPRO3}.
\medskip

\noindent
{\bf Restricting the range of constraint matrices entries:} Similar arguments to those in the proof of Claim \ref{BASICm21} can be used to show the following statement.
\begin{Claim}
The localization conditions in Problem \ref{CONDPRO3} imply that we must have
\begin{equation*}
\underline{m_1}_1,\overline{m_1}_2,\overline{m_1}_3,\overline{m_1}_{2+3}\leq \tfrac12\leq \overline{m_1}_1,\underline{m_1}_{1+2},\underline{m_1}_{1+2+3}
\end{equation*}
and 
\begin{equation*}
\underline{m_2}_2,\overline{m_2}_3,\underline{m_2}_{2+3}\leq \tfrac12\leq \underline{m_2}_1,\overline{m_2}_2,\overline{m_2}_{2+3}\ \text{and}\ 
\underline{m_2}_{1+2},\underline{m_2}_{1+2+3}\leq \tfrac32\leq \overline{m_2}_{1+2},\overline{m_2}_{1+2+3}.
\end{equation*}
\label{BASICm1m2}
\end{Claim}
The previous restrictions are supplemented by the following ones (which are consequences of the assumption $P_k\subset (0,1)^3$)
\[
0\leq \underline{m_k}_{i+\cdots +j}< \overline{m_k}_{i+\cdots +j}\leq j-i+1.
\]

\subsection{Expressions of the optimized constraint matrices associated with atomic restrictions}
{\bf Expression of $O(m_1\cap m_{000101})$:} Anticipating that some of the entries of $m_1$ can be immediately computed using those entries of $O(m_2\cap m_{110111})$ that are numerical values and the condition $G_{3,\epsilon}(P_2\cap A_{110111})=P_{1}$ (see next Subsection) and using also the inequalities in Claim \ref{BASICm1m2}, we have
\[
m_1\cap m_{000101}=\left(
\begin{array}{cccccc}
\tfrac12&\overline{m_1}_2&\overline{m_1}_3&1-\tfrac{\epsilon}2&\overline{m_1}_{2+3}&1-\tfrac{\epsilon}2\\
\underline{m_1}_1&\tfrac{\epsilon}2&0&\underline{m_1}_{1+2}&\tfrac{\epsilon}2&\underline{m_1}_{1+2+3}
\end{array}\right)^T.
\]
The expression of $O(m_1\cap m_{000101})$ then simplifies as follows (see Section 2.1 of Supplementary Material for details of the computation)
\begin{flalign*}
&\underline{O(m_1\cap m_{000101})}_1=\underline{m_1}_1\quad\text{and}\quad \overline{O(m_1\cap m_{000101})}_1=\tfrac12&\\
&\underline{O(m_1\cap m_{000101})}_2=\max\left\{\tfrac{\epsilon}2,\underline{m_1}_{1+2}-\tfrac12\right\}\quad\text{and}\quad\overline{O(m_1\cap m_{000101})}_2=\overline{m_1}_2&\\
&\underline{O(m_1\cap m_{000101})}_3=\max\left\{0,\underline{m_1}_{1+2+3}-\tfrac12-\overline{m_1}_2\right\}\quad\text{and}\quad&\\
&\overline{O(m_1\cap m_{000101})}_3=\min\left\{\overline{m_1}_3,\tfrac12+\overline{m_1}_{2+3}-\underline{m_1}_{1+2}\right\}&\\
&\underline{O(m_1\cap m_{000101})}_{1+2}=\underline{m_1}_{1+2}\quad\text{and}\quad\overline{O(m_1\cap m_{000101})}_{1+2}=\min\left\{1-\tfrac{\epsilon}2,\tfrac12+\overline{m_1}_2\right\}&\\
&\underline{O(m_1\cap m_{000101})}_{2+3}=\max\left\{\tfrac{\epsilon}2,\underline{m_1}_{1+2+3}-\tfrac12\right\}\quad\text{and}\quad\overline{O(m_1\cap m_{000101})}_{2+3}=\overline{m_1}_{2+3}&\\&\underline{O(m_1\cap m_{000101})}_{1+2+3}=\underline{m_1}_{1+2+3}\quad\text{and}\quad\overline{O(m_1\cap m_{000101})}_{1+2+3}=\min\left\{1-\tfrac{\epsilon}2,\tfrac12+\overline{m_1}_{2+3}\right\}&
\end{flalign*} 
\medskip

\noindent
{\bf Expressions of $O(m_2\cap m_{110111})$:} Using the inequalities in Claim \ref{BASICm1m2} and the trivially computed values below of $\overline{m_2}_1$ and $\underline{m_2}_3$, we have
\[
m_2\cap m_{110111}=\left(
\begin{array}{cccccc}
1&\overline{m_2}_2&\overline{m_2}_3&\tfrac32&\min\{\overline{m_2}_{2+3},1\}&\tfrac32\\
\underline{m_2}_1&\tfrac12&0&\max\{\underline{m_2}_{1+2},1\}&\tfrac12&\max\{\underline{m_2}_{1+2+3},1\}
\end{array}\right)^T,
\]
The expression of the corresponding optimized matrix simplifies as follows (see Section 2.1 of Supplementary Material)
\begin{flalign*}
&\underline{O(m_2\cap m_{110111})}_1=\underline{m_2}_1\quad\text{and}\quad \overline{O(m_2\cap m_{110111})}_1=1&\\
&\underline{O(m_2\cap m_{110111})}_2=\tfrac12\quad\text{and}\quad \overline{O(m_2\cap m_{110111})}_2=\min\left\{\overline{m_2}_2, \tfrac32-\underline{m_2}_{1}\right\}&\\
&\underline{O(m_2\cap m_{110111})}_3=0\quad\text{and}\quad \overline{O(m_2\cap m_{110111})}_3=\min\left\{\overline{m_2}_3, \overline{m_2}_{2+3}-\tfrac12\right\}&\\
&\underline{O(m_2\cap m_{110111})}_{1+2}=\max\left\{\underline{m_2}_{1+2}, \underline{m_2}_1+\tfrac12\right\}\quad\text{and}\quad \overline{O(m_2\cap m_{110111})}_{1+2}=\tfrac32&\\
&\underline{O(m_2\cap m_{110111})}_{2+3}=\tfrac12\quad\text{and}\quad \overline{O(m_2\cap m_{110111})}_{2+3}=\min\left\{\overline{m_2}_{2+3}, \tfrac32-\underline{m_2}_1\right\}&\\
&\underline{O(m_2\cap m_{110111})}_{1+2+3}=\max\left\{\underline{m_2}_{1+2+3}, \underline{m_2}_{1}+\tfrac12\right\}\quad\text{and}\quad \overline{O(m_2\cap m_{110111})}_{1+2+3}=\tfrac32&
\end{flalign*}

\subsection{(In)equalities induced by dynamical conditioning}
{\bf  Inequalities induced by $G_{3,\epsilon}(P_1\cap A_{000101})\subset P_2$:}
Together with Corollary \ref{IFFINCLU}, the expression of $G_{3,\epsilon}$ and the value $B_{3,000101}=\tfrac14\left(\begin{array}{c}2\\ 1\\ 0\end{array}\right)$ imply that the dynamics condition $G_{3,\epsilon}(P_1\cap A_{000101})\subset P_{2}$ expresses on constraint matrices as
\begin{equation}
2(1-\epsilon) O(m_1\cap m_{000101})+\tfrac{\epsilon}2\left(\begin{array}{c}2\\ 1\\ 0\end{array}\right)\subset m_2
\label{COND000101}
\end{equation}
In particular, $\overline{O(m_1\cap m_{000101})}_1=\tfrac12$ implies $1\leq \overline{m_2}_1$ and therefore we must have (hence the anticipated expression of $m_2\cap m_{110111}$ above)
\[
\overline{m_2}_1=1\quad\text{and}\quad \underline{m_2}_3=0.
\]
\medskip

\noindent
{\bf Equalities induced by $G_{3,\epsilon}(P_2\cap A_{110111})=P_{1}$:} 
The expression of $G_{3,\epsilon}$ and the value $B_{3,110111}=\tfrac14\left(\begin{array}{c}2\\ 3\\ 0\end{array}\right)$ imply that the dynamics condition $G_{3,\epsilon}(P_2\cap A_{110111})=P_{1}$ expresses on constraint matrices as
\begin{equation}
2(1-\epsilon) O(m_2\cap m_{110111})+\tfrac{\epsilon}2\left(\begin{array}{c}2\\ 3\\ 0\end{array}\right)-\left(\begin{array}{c}1\\ 1\\ 0\end{array}\right)=m_1
\label{COND110111}
\end{equation}
For those entries of $O(m_2\cap m_{110111})$ that are trivial, this equation immediately gives some entries of $m_1$. In this way, one gets (hence the anticipated expression of $m_1\cap m_{000101}$ above)
\[
\overline{m_1}_1=1-\epsilon,\quad \underline{m_1}_2=\underline{m_1}_{2+3}=\tfrac{\epsilon}2,\quad \underline{m_1}_3=0\quad\text{and}\quad \overline{m_1}_{1+2}=\overline{m_1}_{1+2+3}=1-\tfrac{\epsilon}2.
\]
Moreover, Claim \ref{SYM-M2} implies that we have $\underline{O(m_2\cap m_{110111})}_{1+2}+\overline{O(m_2\cap m_{110111})}_{2+3}=2$, which in turn implies, using the corresponding entries of the equation
\[
\underline{m_1}_{1+2}+\overline{m_1}_{2+3}=1.
\]

\subsection{Computing the entries of $m_1$ and $m_2$}
Here, the entries of $m_1$ and $m_2$ that remain to be determined are obtained by considering the remaining (in)equalities from \eqref{COND000101} and \eqref{COND110111} above.
\medskip

\noindent
{\bf Computing $(\underline{m_1}_1,\underline{m_2}_1,\overline{m_2}_3)$:} The inequality that implies the entry $\underline{O(m_1\cap m_{000101})}_1$ in \eqref{COND000101} and the equality with $\underline{O(m_2\cap m_{110111})}_1$ in \eqref{COND110111}  constitute the following independent two-dimensional system, whose solution also provides an expression for $\overline{m_2}_3$ when using the property in Claim \ref{SYM-M2}. We have
\[
\left\{\begin{array}{l}
2(1-\epsilon)\underline{m_1}_1-\underline{m_2}_1\geq -\epsilon\\
2(1-\epsilon)\underline{m_2}_1-\underline{m_1}_1=1-\epsilon
\end{array}\right.
\quad\Longleftrightarrow\quad\left\{\begin{array}{l}
\underline{m_1}_{1}=1-2p^\ast+2(1-\epsilon)\delta_1\\
\underline{m_2}_{1}=2p^\ast+\delta_1\quad \Longleftrightarrow\quad \overline{m_2}_3=1-2p^\ast-\delta_1
\end{array}\right.
\]
for arbitrary $\delta_1\geq 0$,  where $p^\ast=\tfrac{2-\epsilon}{2(3-2\epsilon)}$ (defined before Statement \ref{SOLP4} above) solves the equation \cite{S18}
\begin{equation*}
p^\ast=1-\tfrac{\epsilon}2-2(1-\epsilon)p^\ast.
\end{equation*}

The following statement is proved independently, see Section 2.2 of Supplementary Material.
\begin{Claim}
We must have $\underline{m_2}_1+\tfrac12\leq\underline{m_2}_{1+2}$, or equivalently $\overline{m_2}_{2+3}\leq \overline{m_2}_3+\tfrac12$.
\label{m1p2}
\end{Claim} 
\medskip

\noindent
{\bf Computing $(\overline{m_1}_2,\underline{m_2}_2,\overline{m_2}_2)$:} From Claim \ref{m1p2}, we have
\[
\overline{m_2}_{2}\leq \overline{m_2}_{2+3}-\underline{m_2}_3=\overline{m_2}_{2+3}\leq \overline{m_2}_3+\tfrac12=\tfrac32-\underline{m_2}_1.
\]
Therefore, the equality that implies $\overline{O(m_2\cap m_{110111})}_{2}$ in \eqref{COND110111} and the inequality with $\overline{O(m_1\cap m_{000101})}_{2}$ in \eqref{COND000101} constitute the following independent two-dimensional system, whose solution also provides an expression for $\underline{m_2}_2$ when using the property in Claim \ref{SYM-M2}. We have
\[
\left\{\begin{array}{l}
2(1-\epsilon)\overline{m_1}_2-\overline{m_2}_2\leq -\tfrac{\epsilon}2\\
2(1-\epsilon)\overline{m_2}_2-\overline{m_1}_2=1-\tfrac{3\epsilon}2
\end{array}\right.
\quad\Longleftrightarrow\quad\left\{\begin{array}{l}
\overline{m_1}_{2}=p^\ast-2(1-\epsilon)\delta_2\\
\overline{m_2}_{2}=1-p^\ast-\delta_2\quad \Longleftrightarrow\quad \underline{m_2}_2=p^\ast+\delta_2
\end{array}\right.
\]
for arbitrary $\delta_2\geq 0$.
\medskip

\noindent
{\bf Computing $(\underline{m_1}_{1+2},\overline{m_1}_{2+3},\underline{m_2}_{1+2},\overline{m_2}_{2+3})$:} From Claim \ref{m1p2}, the equality that implies $\underline{O(m_2\cap m_{110111})}_{1+2}$ in \eqref{COND110111} and the inequality with $\underline{O(m_1\cap m_{000101})}_{1+2}$ in \eqref{COND000101} constitute the following independent two-dimensional system, 
\[
\left\{\begin{array}{l}
2(1-\epsilon)\underline{m_1}_{1+2}-\underline{m_2}_{1+2}\geq -\tfrac{3\epsilon}2\\
2(1-\epsilon)\underline{m_2}_{1+2}-\underline{m_1}_{1+2}=2-\tfrac{5\epsilon}2
\end{array}\right.
\]
whose solution also provides the expression of $(\overline{m_1}_{2+3},\overline{m_2}_{2+3})$, when using symmetries relationships
\[
\left\{\begin{array}{l}
\underline{m_1}_{1+2}=1-p^\ast+2(1-\epsilon)\delta_3\\
\underline{m_2}_{1+2}=1+p^\ast+\delta_3
\end{array}\right.
\quad\Longleftrightarrow\quad\left\{\begin{array}{l}
\overline{m_1}_{2+3}=p^\ast-2(1-\epsilon)\delta_3\\
\overline{m_2}_{2+3}=1-p^\ast-\delta_3
\end{array}\right.
\]
for arbitrary $\delta_3\geq 0$.
\medskip

\noindent
{\bf Computing $\overline{m_1}_3$:} From Claim \ref{m1p2}, the equality that implies $\overline{O(m_2\cap m_{110111})}_{3}$ together with the expression of $\overline{m_2}_{2+3}$ yields the following expression
\[
\overline{m_1}_3=p^\ast-\tfrac{\epsilon}2-2(1-\epsilon)\delta_3.
\]
\medskip

\noindent
{\bf Computing $(\underline{m_1}_{1+2+3},\underline{m_2}_{1+2+3},\overline{m_2}_{1+2+3})$:} Claim \ref{m1p2} implies $\underline{m_2}_1+\tfrac12\leq\underline{m_2}_{1+2+3}$; hence the above equality that implies $\underline{O(m_2\cap m_{110111})}_{1+2+3}$ in \eqref{COND000101} and the above inequality with $\underline{O(m_1\cap m_{000101})}_{1+2+3}$ in \eqref{COND110111} constitute the same independent two-dimensional system as the one for $(\underline{m_1}_{1+2},\underline{m_2}_{1+2})$. Its solution also provides an expression for $\overline{m_2}_{1+2+3}$ when using the property in Claim \ref{SYM-M2}. We have
\[
\left\{\begin{array}{l}
\underline{m_1}_{1+2+3}=1-p^\ast+2(1-\epsilon)\delta_4\\
\underline{m_2}_{1+2+3}=1+p^\ast+\delta_4\quad \Longleftrightarrow\quad \overline{m_2}_{1+2+3}=2-p^\ast-\delta_4
\end{array}\right.
\]
for arbitrary $\delta_4\geq 0$.
\medskip

\noindent
{\bf Computing $(\overline{m_2}_{1+2},\underline{m_2}_{2+3})$:} At this stage, all entries of $m_1$ and $m_2$ have been determined, excepted $\overline{m_2}_{1+2}$ and  $\underline{m_2}_{2+3}$. All information \eqref{COND110111} has been used above; hence we use instead that we must have $\overline{m_2}_{1+2}=\overline{O(m_2)}_{1+2}$, {\sl ie.}\ 
\begin{align*}
\overline{m_2}_{1+2}&\leq\min\left\{\overline{m_2}_1+\overline{m_2}_2,\overline{m_2}_{1+2+3}-\underline{m_2}_3, \overline{m_2}_1+\overline{m_2}_{2+3}-\underline{m_2}_3,\overline{m_2}_2+\overline{m_2}_{1+2+3}-\underline{m_2}_{2+3}\right\}\\
&=\min\left\{2-p^\ast-\max\{\delta_2,\delta_3,\delta_4\},1-2p^\ast-\delta_2-\delta_4+\overline{m_2}_{1+2}\right\}
\end{align*}
which imposes
\[
\delta_2+\delta_4\leq 1-2p^\ast\quad\text{and}\quad
\left\{\begin{array}{l}
\overline{m_2}_{1+2}=2-p^\ast-\delta_5\\
\underline{m_2}_{2+3}=p^\ast+\delta_5
\end{array}\right.
\]
for arbitrary $\delta_5\geq \max\{\delta_2,\delta_3,\delta_4\}$.

In summary, we have
\begin{align}
m_1=&\left(
\begin{array}{cccc}
1-\epsilon&p^\ast-2(1-\epsilon)\delta_2&p^\ast-\tfrac{\epsilon}2-2(1-\epsilon)\delta_3&1-\tfrac{\epsilon}2\\
1-2p^\ast+2(1-\epsilon)\delta_1&\tfrac{\epsilon}2&0&1-p^\ast+2(1-\epsilon)\delta_3
\end{array}\right.\nonumber\\
&\left.
\begin{array}{cc}
p^\ast-2(1-\epsilon)\delta_3&1-\tfrac{\epsilon}2\\
\tfrac{\epsilon}2&1-p^\ast+2(1-\epsilon)\delta_4
\end{array}\right)^T
\label{m1}
\end{align}
and
\begin{equation}
m_2=\left(
\begin{array}{cccccc}
1&1-p^\ast-\delta_2&1-2p^\ast-\delta_1&2-p^\ast-\delta_5&1-p^\ast-\delta_3&2-p^\ast-\delta_4\\
2p^\ast+\delta_1&p^\ast+\delta_2&0&1+p^\ast+\delta_3&p^\ast+\delta_5&1+p^\ast+\delta_4
\end{array}\right)^T
\label{m2}
\end{equation}

\subsection{Proof of Statement \ref{SOLP3}}
To prove this statement, one needs to check that, in addition to $m_1=O(m_1)$ and $m_2=O(m_2)$, the entries of $m_1$ and $m_2$ satisfy all the required conditions in the previous sections. This verification uses elementary algebra and assumes $\epsilon\in (0,\tfrac12)$ which in particular implies
\[
\tfrac12< 2p^\ast<1,
\]
and $\min_{i\in\{1,\cdots ,4\}}\delta_i\geq 0$ and $\delta_5\geq\max\{\delta_2,\delta_3,\delta_4\}$.
\medskip

\noindent 
{\bf Checking the basic inequalities $0\leq \underline{m_k}_{i+\cdots +j}\leq \overline{m_k}_{i+\cdots +j}\leq j-i+1$:} For $k=1$ and $k=2$ we respectively get the following additional restrictions
\[
\left\{\begin{array}{l}
2(1-\epsilon)\delta_1\leq 2p^\ast-\epsilon\\
2(1-\epsilon)\max\{\delta_2,\delta_3,\delta_4\}\leq p^\ast-\tfrac{\epsilon}2
\end{array}\right.
\quad\text{and}\quad
\left\{\begin{array}{l}
\max\{\delta_1,\delta_3+\delta_5\}\leq 1-2p^\ast\\
\max\{\delta_2,\delta_4\}\leq\tfrac{1}2-p^\ast
\end{array}\right.
\]
\medskip

\noindent 
{\bf Checking the inequalities in Claim \ref{BASICm1m2}:} These inequalities yield the following additional restrictions
\[
2(1-\epsilon)\delta_1\leq 2p^\ast-\tfrac12\quad\text{and}\quad
\delta_5\leq \tfrac12-p^\ast.
\]
At this stage, all restrictions imposed by the localisation conditions in Problem \ref{CONDPRO3} have been included.
\medskip

\noindent 
{\bf Checking the inequalities in Claim \ref{m1p2}:} The two inequalities are equivalent and yield the following additional restriction
\[
\delta_1\leq \tfrac12-p^\ast+\delta_3
\]
\medskip

\noindent 
{\bf Checking the remaining inequalities in the condition $G_{3,\epsilon}(P_1\cap A_{000101})\subset P_{2}$:} The inequalities in \eqref{COND000101} that involve the entries
\[
\underline{O(m_1\cap m_{000101})}_2,\ \underline{O(m_1\cap m_{000101})}_3,\ \overline{O(m_1\cap m_{000101})}_3,\ \overline{O(m_1\cap m_{000101})}_{1+2},
\] 
and
\[
\underline{O(m_1\cap m_{000101})}_{2+3}\ \text{and}\ \overline{O(m_1\cap m_{000101})}_{1+2+3}
\]
remain to be considered. Together with the inequality $2(1-\epsilon)\max\{\delta_2,\delta_3,\delta_4\}\leq p^\ast-(1-\epsilon)^2$ (which will be obtained independently below) they yield the following additional constraints
\begin{flalign*}
&2(1-\epsilon)\min\left\{p^\ast-\tfrac{\epsilon}2-2(1-\epsilon)\delta_3,2p^\ast-\tfrac12-4(1-\epsilon)\delta_3\right\}\leq 1-2p^\ast-\delta_1&\\
&\delta_5\leq \epsilon(\tfrac32-\epsilon)-p^\ast\quad \left(\text{and note that}\ \epsilon(\tfrac32-\epsilon)\leq \tfrac12\ \text{for}\ \epsilon\in (0,\tfrac12)\right)&
\end{flalign*} 
\medskip

\noindent 
{\bf Checking the inequalities resulting from the condition $G_{3,\epsilon}(P_1\cap A_{100101})\subset P_{1}$:} Using Claim \ref{BASICm1m2} and some of the entries computed above, we have
\[
m_1\cap m_{100101}=\left(
\begin{array}{cccccc}
1-\epsilon&\overline{m_1}_2&\overline{m_1}_3&1-\tfrac{\epsilon}2&\overline{m_1}_{2+3}&1-\tfrac{\epsilon}2\\
\tfrac12&\tfrac{\epsilon}2&0&\underline{m_1}_{1+2}&\tfrac{\epsilon}2&\underline{m_1}_{1+2+3}
\end{array}\right)^T.
\]
Using also the explicit entries of $m_1$ above, the expression of $O(m_1\cap m_{100101})$ then simplifies as follows (see Section 2.1 of Supplementary Material)
\begin{align*}
&\underline{O(m_1)}_1=\max\left\{\tfrac12,1-2p^\ast+2(1-\epsilon)(\max\{\delta_2,\delta_4\}+\delta_3)\right\}\quad\text{and}\quad\overline{O(m_1)}_1=1-\epsilon\\
&\underline{O(m_1)}_2=\tfrac{\epsilon}2\quad\text{and}\quad\overline{O(m_1)}_2=\min\left\{p^\ast-2(1-\epsilon)\delta_2,\tfrac{1-\epsilon}2\right\}\\
&\underline{O(m_1)}_3=0\quad\text{and}\quad\overline{O(m_1)}_3=\min\left\{p^\ast-\tfrac{\epsilon}2-2(1-\epsilon)\delta_3,\ \tfrac12-\epsilon\right\}\\
&\underline{O(m_1)}_{1+2}=\max\left\{1-p^\ast+2(1-\epsilon)\delta_3,\tfrac{1+\epsilon}2\right\}\quad\text{and}\quad\overline{O(m_1)}_{1+2}=1-\tfrac{\epsilon}2\\
&\underline{O(m_1)}_{2+3}=\tfrac{\epsilon}2\quad\text{and}\quad\overline{O(m_1)}_{2+3}=\min\left\{p^\ast-2(1-\epsilon)\delta_3,\ \tfrac{1-\epsilon}2\right\}\\
&\underline{O(m_1)}_{1+2+3}=\max\left\{1-p^\ast+2(1-\epsilon)\delta_4,\tfrac{1+\epsilon}2\right\}\quad\text{and}\quad\overline{O(m_1)}_{1+2+3}=1-\tfrac{\epsilon}2
\end{align*}
The expression of $G_{3,\epsilon}$ and the value $B_{3,100101}=\left(\begin{array}{c}1\\ 0\\ 0\end{array}\right)$ imply that the dynamics condition $G_{3,\epsilon}(P_1\cap A_{100101})\subset P_{1}$ expresses on constraint matrices as
\[
2(1-\epsilon) O(m_1\cap m_{100101})-(1-2\epsilon)\left(\begin{array}{c}1\\ 0\\ 0\end{array}\right)\subset m_1
\]
Using the expression of $O(m_1\cap m_{100101})$ above, one obtains using elementary algebra that this condition yields the following additional conditions
\begin{align*}
&2(1-\epsilon)\delta_1\leq 2p^\ast-1+\epsilon\\
&2(1-\epsilon)\max\{\delta_2,\delta_3,\delta_4\}\leq p^\ast-(1-\epsilon)^2
\end{align*}
\medskip

\noindent 
{\bf Checking the inequalities resulting from the condition $G_{3,\epsilon}(P_2\cap A_{110212})\subset P_2$:}
Claim \ref{BASICm1m2}, the values above of $\overline{m_2}_1$ and $\underline{m_2}_3$ and also the inequalities $\overline{m_2}_{2+3}\leq \overline{m_2}_{2}+\overline{m_2}_{3}\leq \tfrac32$ and $\overline{m_2}_{1+2+3}\leq \overline{m_2}_{1}+\overline{m_2}_{2+3}\leq \tfrac52$, we have
\[
m_2\cap m_{110212}=\left(
\begin{array}{cccccc}
1&\overline{m_2}_2&\overline{m_2}_3&\overline{m_2}_{1+2}&\overline{m_2}_{2+3}&\overline{m_2}_{1+2+3}\\
\underline{m_2}_1&\tfrac12&0&\tfrac32&\tfrac12&\tfrac32
\end{array}\right)^T,
\]
The expression of the corresponding optimized matrix simplifies as follows (see Section 2.1 of Supplementary Material)
\begin{flalign*}
&\underline{O(m_2\cap m_{110212})}_1=\max\left\{2p^\ast+\delta_1,\tfrac12+p^\ast+\delta_2\right\}\quad\text{and}\quad \overline{O(m_2\cap m_{110212})}_1=1&\\
&\underline{O(m_2\cap m_{110212})}_2=\tfrac12\quad\text{and}\quad \overline{O(m_2\cap m_{110212})}_2=1-p^\ast-\delta_2&\\
&\underline{O(m_2\cap m_{110212})}_3=0\quad\text{and}\quad \overline{O(m_2\cap m_{110212})}_3=\min\left\{1-2p^\ast-\delta_1,\tfrac12-p^\ast-\delta_4\right\}&\\
&\underline{O(m_2\cap m_{110212})}_{1+2}=\tfrac32\quad\text{and}\quad \overline{O(m_2\cap m_{110212})}_{1+2}=2-p^\ast-\delta_5&\\
&\underline{O(m_2\cap m_{110212})}_{2+3}=\tfrac12\quad\text{and}\quad \overline{O(m_2\cap m_{110212})}_{2+3}=\min\left\{1-p^\ast-\delta_3,\tfrac32-2p^\ast-(\delta_2+\delta_4)\right\}&\\
&\underline{O(m_2\cap m_{110212})}_{1+2+3}=\tfrac32\quad\text{and}\quad \overline{O(m_2\cap m_{110212})}_{1+2+3}=2-p^\ast-\delta_4&
\end{flalign*}
The expression of $G_{3,\epsilon}$ and the value $B_{3,110212}=\left(\begin{array}{c}1\\ 1\\ 0\end{array}\right)$ imply that the dynamics condition $G_{3,\epsilon}(P_2\cap A_{110212})\subset P_{2}$ expresses on constraint matrices as
\[
2(1-\epsilon) O(m_2\cap m_{110212})-(1-2\epsilon)\left(\begin{array}{c}1\\ 1\\ 0\end{array}\right)\subset m_2
\]
Using the expression of $O(m_2\cap m_{110212})$ above, one obtains after elementary algebra that this condition yields the following additional conditions
\begin{flalign*}
&\delta_1\leq \tfrac12-p^\ast+\min\{\delta_2,\delta_4\}\quad\text{and}\quad \delta_1\leq \epsilon(1-2p^\ast)+2(1-\epsilon)\min\{\delta_2,\delta_4\}\quad\text{and}\quad \delta_5\leq \epsilon-p^\ast&
\end{flalign*}
\medskip

\noindent 
{\bf Checking the inequalities resulting from the condition $G_{3,\epsilon}(P_2\cap A_{110112})\subset \sigma_{4231}(P_{1})$:} Using Claim \ref{BASICm1m2} and some of the entries computed above (including that $\overline{m_2}_{1+2+3}\leq 2$, we have
\[
m_2\cap m_{110212}=\left(
\begin{array}{cccccc}
1&\overline{m_2}_2&\overline{m_2}_3&\tfrac32&\overline{m_2}_{2+3}&\overline{m_2}_{1+2+3}\\
\underline{m_2}_1&\tfrac12&0&\underline{m_2}_{1+2}&\tfrac12&\tfrac32
\end{array}\right)^T,
\]
Using also the explicit entries of $m_2$ above, the expression of $O(m_2\cap m_{110112})$ then simplifies as follows (see Section 2.1 of Supplementary Material)
\begin{align*}
&\underline{O(m_2\cap m_{110112})}_1=\max\left\{2p^\ast+\delta_1,\tfrac12+p^\ast+\delta_3\right\}\quad\text{and}\quad \overline{O(m_2\cap m_{110112})}_1=1\\
&\underline{O(m_2\cap m_{110112})}_2=\tfrac12\quad\text{and}\quad\overline{O(m_2\cap m_{110112})}_2=\min\left\{1-p^\ast-\delta_2,\ \tfrac32-2p^\ast-\delta_1\right\}\\
&\underline{O(m_2\cap m_{110112})}_3=0\quad\text{and}\quad \overline{O(m_2\cap m_{110112})}_3=\min\left\{1-2p^\ast-\delta_1,\ \tfrac12-p^\ast-\delta_3\right\}\\
&\underline{O(m_2\cap m_{110112})}_{1+2}=\max\left\{1+p^\ast+\delta_3,\ 2p^\ast+\delta_1+\tfrac12\right\}\quad\text{and}\quad\overline{O(m_2\cap m_{110112})}_{1+2}=\tfrac32\\
&\underline{O(m_2\cap m_{110112})}_{2+3}=\tfrac12\quad\text{and}\quad\overline{O(m_2\cap m_{110112})}_{2+3}=\min\left\{1-p^\ast-\delta_3, \tfrac52-4p^\ast-2\delta_1\right\}\\
&\underline{O(m_2\cap m_{110112})}_{1+2+3}=\tfrac32\quad\text{and}\quad\overline{O(m_2\cap m_{110112})}_{1+2+3}=\min\left\{2-p^\ast-\max\{\delta_3,\delta_4\},\ \tfrac52-2p^\ast-\delta_1\right\}
\end{align*}
In addition, the expression of $m_1$ implies that $0\leq \underline{m_1}_{1+2}<\overline{m_1}_{1+2}$ and $0\leq \underline{m_1}_{2+3}<\overline{m_1}_{2+32}$. Therefore, we have
\[
\sigma_{4231}|_{P_{m_1}} (x)=(1-x_2-x_3,x_2,1-x_1-x_2),
\]
which satisfies \eqref{SPETRANS} (but not the assumptions of the second Claim \ref{CLAIMSPETRANS}), and then 
\[
\sigma_{4231}(m_1)=\left(
\begin{array}{cccccc}
1-\underline{m_1}_{2+3}&\overline{m_1}_2&1-\underline{m_1}_{1+2}&1-\underline{m_1}_{3}&1-\underline{m_1}_{1}&2-\underline{m_1}_{1+2+3}\\
1-\overline{m_1}_{2+3}&\underline{m_1}_2&1-\overline{m_1}_{1+2}&1-\overline{m_1}_{3}&1-\overline{m_1}_{1}&2-\overline{m_1}_{1+2+3}
\end{array}\right)^T
\]

The expression of $G_{3,\epsilon}$ and the value $B_{3,110112}=\left(\begin{array}{c}1\\ 1\\ 0\end{array}\right)$ imply that the dynamics condition $G_{3,\epsilon}(P_{m_2}\cap A_{110112})\subset \sigma_{4231}(P_{m_1})$ requires that the following inequality on constraints matrices hold 
\[
2(1-\epsilon) O(m_2\cap m_{110112})+\tfrac{\epsilon}2\left(\begin{array}{c}3\\ 3\\ 1\end{array}\right)-\left(\begin{array}{c}1\\ 1\\ 0\end{array}\right)\subset \sigma_{4231}(m_{1})
\]
Using the expression of $O(m_2\cap m_{110112})$ and $m_1$ above, one checks, after elementary algebra, that this condition does not bring additional restriction on the parameters.
\medskip

\noindent 
{\bf Checking optimality:} Using the explicit values of $m_1$ and $m_2$ in the conditions $m_1=O(m_1)$ and $m_2=O(m_2)$ respectively yields the following additional restrictions (see Section 2.3 of Supplementary material)
\[
\delta_3+\max\left\{\delta_2,\delta_4\right\}\leq \delta_1,\quad \delta_3\leq \min\{\delta_2,\delta_4\}\quad\text{and}\quad
\delta_5+\delta_1-\min\{\delta_2,\delta_4\}\leq 1-2p^\ast
\]
\medskip

\noindent 
{\bf Resulting conditions on parameters:} Collecting all conditions on parameters above, one finally obtains the following inequalities
\begin{equation}
\left\{\begin{array}{l}
0\leq 2(1-\epsilon)\delta_1\leq 2p^\ast-(1-\epsilon)\\
\delta_1\leq \epsilon(1-2p^\ast)+2(1-\epsilon)\min\{\delta_2,\delta_4\}\\
0\leq 2(1-\epsilon)\max\{\delta_2,\delta_4\}\leq p^\ast-(1-\epsilon)^2
\end{array}\right.
\quad\text{and}\quad
\left\{\begin{array}{l}
0\leq \delta_3\leq \min\{\delta_2,\delta_4\} \\
\delta_3+\max\left\{\delta_2,\delta_4\right\}\leq \delta_1\\
\max\{\delta_2,\delta_4\}\leq \delta_5\leq\epsilon-p^\ast\\
\delta_5+\delta_1\leq 1-2p^\ast+\min\{\delta_2,\delta_4\}
\end{array}\right.
\label{DELTAEPS}
\end{equation}
Now, the condition $(1-\epsilon)^2\leq p^\ast$ is equivalent to $\epsilon\geq \epsilon_3$ and this inequality implies $p^\ast <\epsilon$. Together with $\tfrac12<2p^\ast<1$, we conclude that when $\epsilon\geq \epsilon_3$, there exists $\Delta_\epsilon\subset \R^5$, which contains the origin,\footnote{Of note, $\Delta_{\epsilon_3}=\{0\}$.} such that every quintuple $\delta=(\delta_i)_{i=1}^5\in\Delta_\epsilon$ satisfies the inequalities above. In particular, every quintuple as follows
\[
\left\{\begin{array}{l}
0\leq \delta_1\leq \min\left\{\tfrac{2p^\ast-(1-\epsilon)}{2(1-\epsilon)},\epsilon(1-2p^\ast)\right\}\\
0\leq \min\{\delta_2,\delta_4\}\leq\max\{\delta_2,\delta_4\}\leq \min\left\{\tfrac{p^\ast-(1-\epsilon)^2}{2(1-\epsilon)},\delta_1,1-2p^\ast-\delta_1,\epsilon-p^\ast\right\}\\
0\leq \delta_3\leq \min\left\{\delta_2,\delta_4,\delta_1-\max\{\delta_2,\delta_4\}\right\}\\
\max\{\delta_2,\delta_4\}\leq \delta_5\leq \min\left\{1-2p^\ast-\delta_1,\epsilon-p^\ast\right\}
\end{array}\right.
\]
belong to $\Delta_\epsilon$ (so that $\Delta_\epsilon$ must have positive volume when $\epsilon>\epsilon_3$).
\medskip

\noindent 
{\bf Proof that $\text{Orb}_{\langle \sigma_{4231},\sigma_{1324}\rangle}(P_1,P_2)$ is an AsIUP with respect to $\Sigma=1-\text{Id}$.} This immediately follows from the facts that, when the explicit dependence on the parameters is included, we have for $i=1,2$
\begin{itemize}
\item[] $P_{m_i(\delta)}\subset P_{m_i(0)}$ $\forall \delta\in\Delta_\epsilon$,
\item[] $\Sigma(P_{m_i(0)})\cap \text{Orb}_{\langle \sigma_{4231},\sigma_{1324}\rangle}(P_{m_1(0)},P_{m_2(0)})$ \cite{S18}.
\end{itemize}
Statement \ref{SOLP3} is proved. \hfill $\Box$

\section{Analysis of the Conditioning Problem \ref{CONDPRO4}}
\subsection{Proof of Statement \ref{NOSOLUTIONP4}}\label{A-NOSOLUTIONP4}
The constraint matrices under consideration here are those related to the canonical coefficient matrix; hence we use the same notations as those in Appendix \ref{A-P3}.

The proof proceeds by contradiction and aims to show that if $m=O(m)$ is an optimized constraint matrix of a polytope $P_{m}$ which satisfies Conditioning Problem \ref{CONDPRO4}, then its entries must satisfy
\[
2(1-\epsilon)\overline{m}_1+\underline{m}_1\leq 1-\epsilon\leq 2(1-\epsilon)\underline{m}_1+ \overline{m}_1
\]
which clearly implies $\overline{m}_1\leq \underline{m}_1$; hence $P_m=\emptyset$.
\medskip

\noindent
{\bf Optimized constraint matrices associated with atoms involved in} Conditioning Problem \ref{CONDPRO4}: In addition to $m_{000101}$, the optimized constraint matrices $m_\omega$ associated with the involved atoms $P_{m_{\omega}}$ are
\[
m_{000000}=\left(
\begin{array}{cccccc}
\tfrac12&\tfrac12&\tfrac12&\tfrac12&\tfrac12&\tfrac12\\
0&0&0&0&0&0
\end{array}\right)^T,
\quad
m_{000001}=\left(
\begin{array}{cccccc}
\tfrac12&\tfrac12&\tfrac12&\tfrac12&\tfrac12&1\\
0&0&0&0&0&\tfrac12
\end{array}\right)^T,
\]

From now on, we assume that $m=O(m)$ is an optimzed constraint matrix of a polytope $P=P_m$ which satisfies Conditioning Problem \ref{CONDPRO4}.
\medskip

\noindent
 {\bf Restricting the range of constraint matrices entries:} Similar arguments to those in the proof of Claim \ref{BASICm21} can be used to show that the following statement holds.
\begin{Claim}
We must have
\[
\overline{m}_1,\overline{m}_2,\overline{m}_3,\underline{m}_{1+2},\overline{m}_{2+3},\underline{m}_{1+2+3}<\tfrac12<\overline{m}_{1+2},\overline{m}_{1+2+3}<1.
\]
\label{D1}
\end{Claim}

To proceed, it suffices to consider those entries of the atomic restrictions that are involved in the determination of the coordinates $\underline{m}_1$ and $\overline{m}_1$. 
\medskip

\noindent
 {\bf Expression of $O(m\cap m_{000000})$:} According to the previous Claim and the definition of $m_{000000}$, we have
\[
m\cap m_{000000}=\left(
\begin{array}{cccccc}
\overline{m}_1&\overline{m}_2&\overline{m}_3&\tfrac12&\overline{m}_{2+3}&\tfrac12\\
\underline{m}_1&\underline{m}_2&\underline{m}_3&\underline{m}_{1+2}&\underline{m}_{2+3}&\underline{m}_{1+2+3}
\end{array}\right)^T
\]
which implies 
\begin{align*}
\overline{O(m\cap m_{000000})}_2&=\min\left\{\overline{m}_2,\ \tfrac12-\underline{m}_1,\ \overline{m}_{2+3}-\underline{m}_3,\ \tfrac12-\underline{m}_1-\underline{m}_3,\ \tfrac12+\overline{m}_{2+3}-\underline{m}_{1+2+3}\right\}\\
&=\min\left\{\overline{m}_2,\ \tfrac12-\underline{m}_1-\underline{m}_3\right\}
\end{align*}
following similar arguments as before.

That the condition $G_{3,\epsilon}(P\cap A_{000000})\subset P\cap A_{000101}$ expresses on constraints matrix (thanks to Corollary \ref{IFFINCLU}) as 
\begin{equation}
2(1-\epsilon)O(m\cap m_{000000})\subset m\cap m_{000101}
\label{INCL000}
\end{equation}
implies in particular, and together with the expression of $m\cap m_{000101}$ below, that we must have $\tfrac12-\underline{m}_1-\underline{m}_3<\overline{m}_2$. Anticipating that we are going to show that $\underline{m}_3=0$, we must in fact have
\begin{equation}
\tfrac12-\underline{m}_1<\overline{m}_2
\label{LOWER}
\end{equation}
\medskip

\noindent
{\bf Expression of $O(m\cap m_{000001})$:} According to Claim \ref{D1} and the definition of $m_{000001}$, we have
\[
m\cap m_{000001}=\left(
\begin{array}{cccccc}
\overline{m}_1&\overline{m}_2&\overline{m}_3&\tfrac12&\overline{m}_{2+3}&\overline{m}_{1+2+3}\\
\underline{m}_1&\underline{m}_2&\underline{m}_3&\underline{m}_{1+2}&\underline{m}_{2+3}&\tfrac12
\end{array}\right)^T
\]
which implies 
\begin{align*}
\overline{O(m\cap m_{000001})}_2&=\min\left\{\overline{m}_2,\tfrac12-\underline{m}_1,\ \overline{m}_{2+3}-\underline{m}_3,\ \overline{m}_{1+2+3}-\underline{m}_1-\underline{m}_3,\overline{m}_{2+3}\right\}\\
&=\min\left\{\overline{m}_2, \tfrac12-\underline{m}_1\right\}=\tfrac12-\underline{m}_1
\end{align*}
where we used Claim \ref{D1} and the inequality \eqref{LOWER}. Moreover, Claim \ref{D1} implies $\overline{m}_{1+2}<1$ so that we have
\[
\sigma_{3124}|_{P_{m}} (x)=(1-x_1-x_2,x_1,x_2+x_3),
\]
which satisfies \eqref{SPETRANS}. Therefore, the condition $G_{3,\epsilon}(P\cap A_{000001})\subset \sigma_{3124}(P)$ requires that the following inequality on constraints matrices holds
\begin{equation}
2(1-\epsilon)O(m\cap m_{000001})+\tfrac{\epsilon}2\left(\begin{array}{c}1\\ 0\\ 1\end{array}\right)\subset \sigma_{3124}(m)
\label{LAST001}
\end{equation}
In particular, the inequality $2(1-\epsilon) \overline{O(m\cap m_{000001})}_2\leq \overline{\sigma_{3124}(m)}_2=\overline{m}_1$ then yields $1-\epsilon\leq 2(1-\epsilon)\underline{m}_1+ \overline{m}_1$, which is nothing but the second of the desired inequalities.

In order to show $\underline{m}_3=0$, we first compute
\[
\underline{O(m\cap m_{000001})}_{1+2+3}=\max\left\{\tfrac12,\ \underline{m}_1+\underline{m}_{2+3},\ \underline{m}_3+\underline{m}_{1+2},\ \underline{m}_1+\underline{m}_2+\underline{m}_3,\ \underline{m}_{1+2}+\underline{m}_{2+3}-\overline{m}_2\right\}=\tfrac12
\]
because $\underline{O(m)}_{1+2+3}=\underline{m}_{1+2+3}<\tfrac12$. 
The inequality 
\[
1+\underline{m}_3=\underline{\sigma_{3124}(m)}_{1+2+3}\leq 2(1-\epsilon)\underline{O(m\cap m_{000001})}_{1+2+3}+\epsilon,
\]
then yields $\underline{m}_3\leq 0$. Therefore, we must have $\underline{m}_3=0$.
\medskip

\noindent
{\bf Expression of $O(m\cap m_{000101})$:} According to Claim \ref{D1} and the definition of $m_{000101}$, we have
\[
m\cap m_{000101}=\left(
\begin{array}{cccccc}
\overline{m}_1&\overline{m}_2&\overline{m}_3&\overline{m}_{1+2}&\overline{m}_{2+3}&\overline{m}_{1+2+3}\\
\underline{m}_1&\underline{m}_2&\underline{m}_3&\tfrac12&\underline{m}_{2+3}&\tfrac12
\end{array}\right)^T
\]
which implies 
\begin{align*}
\overline{O(m\cap m_{000101})}_1&=\min\left\{\overline{m}_1,\ \overline{m}_{1+2}-\underline{m}_2,\ \overline{m}_{1+2+3}-\underline{m}_{2+3},\ \overline{m}_{1+2+3}-\underline{m}_2-\underline{m}_3,\ \overline{m}_3+\overline{m}_{1+2}-\underline{m}_{2+3}\right\}\\
&=\overline{m}_1
\end{align*}
Moreover, Claim \ref{D1} implies that we have
\[
\sigma_{2134}|_{P^\alpha_m}(x)=\left(1-x_1,x_1+x_2,x_3\right)
\]
which satisfies \eqref{SPETRANS}. Hence, the condition $G_{3,\epsilon}(P\cap A_{000101})\subset \sigma_{2134}(P)$ requires that the following inequality on constraints matrices hold
\begin{equation}
2(1-\epsilon)O(m\cap m_{000101})+\tfrac{\epsilon}2\left(\begin{array}{c}2\\ 1\\ 0\end{array}\right)\subset \sigma_{2134}(m)
\label{LASTINC}
\end{equation}
and we have $\overline{\sigma_{2134}(m)}_1=1-\underline{m}_1$. The first of the desired inequality is a direct consequence of $2(1-\epsilon)\overline{O(m\cap m_{000101})}_1+\epsilon\leq \overline{\sigma_{2134}(m)}_1$.\hfill $\Box$

\subsection{Proof of Statement \ref{SOLP4}}\label{A-SOLP4}
This appendix contains the structure of the proof of Statement \ref{SOLP4}. Computational details can be found in the Mathematica notebook "AsIUPG3SeondBif.nb" in Supplementary Material. 
As announced, we use constraint matrices in order to characterize polytopes. Given the half-spaces that are involved in the definitions of both the solution polytope $P$ and the atoms involved in Problem \ref{CONDPRO4}, and given that the symmetries $\sigma_{2134}$ and $\sigma_{3124}$ should satisfy the condition \eqref{SPETRANS}, consider the matrix $\alpha$ such that the polytopes $P^{\alpha}_m$ are defined as follows 
\[
P^{\alpha}_m=\left\{x\in\R^3\ :\ \begin{array}{c}
\underline{m}_1<x_1<\overline{m}_1\\
\underline{m}_2<x_2<\overline{m}_2\\
\underline{m}_3<x_3<\overline{m}_3\\
\underline{m}_{1+2}<x_1+x_2<\overline{m}_{1+2}\\
\underline{m}_{2+3}<x_2+x_3<\overline{m}_{2+3}\\
\underline{m}_{1+2+3}<x_1+x_2+x_3<\overline{m}_{1+2+3}\\
\underline{m}_{1+2\cdot 2+3\cdot 3}<x_1+2x_2+3x_3<\overline{m}_{1+2\cdot 2+3\cdot 3}\\
\underline{m}_{p_1^\ast\cdot 1+p_2^\ast\cdot 2+p_3^\ast\cdot 3}<p_1^\ast x_1+p_2^\ast x_2+p_3^\ast x_3<\overline{m}_{p_1^\ast\cdot 1+p_2^\ast\cdot 2+p_3^\ast\cdot 3}\\
\underline{m}_{-p^\ast\cdot 1+p_2^\ast\cdot 2+p_3^\ast\cdot 3}<-p^\ast x_1+p_2^\ast x_2+p_3^\ast x_3<\overline{m}_{-p^\ast\cdot 1+p_2^\ast\cdot 2+p_3^\ast\cdot 3}\\
\underline{m}_{-p^\ast\cdot 1-2p^\ast\cdot 2+p_3^\ast\cdot 3}<-p^\ast x_1-2p^\ast x_2+p_3^\ast x_3<\overline{m}_{-p^\ast\cdot 1-2p^\ast\cdot 2+p_3^\ast\cdot 3}
\end{array}\right\}
\]
where
{\small\[
m=\left(
\begin{array}{cccccccccc}
\overline{m}_1&\overline{m}_2&\overline{m}_3&\overline{m}_{1+2}&\overline{m}_{2+3}&\overline{m}_{1+ 2+ 3}&\overline{m}_{1+2\cdot 2+3\cdot 3}&\overline{m}_{p_1^\ast\cdot 1+p_2^\ast\cdot 2+p_3^\ast\cdot 3}&\overline{m}_{-p^\ast\cdot 1+p_2^\ast\cdot 2+p_3^\ast\cdot 3}&\overline{m}_{-p^\ast\cdot 1-2p^\ast\cdot 2+p_3^\ast\cdot 3}\\
\underline{m}_1&\underline{m}_2&\underline{m}_3&\underline{m}_{1+2}&\underline{m}_{2+3}&\underline{m}_{1+ 2+3}&\underline{m}_{1+2\cdot 2+3\cdot 3}&\underline{m}_{p_1^\ast\cdot 1+p_2^\ast\cdot 2+p_3^\ast\cdot 3}&\underline{m}_{-p^\ast\cdot 1+p_2^\ast\cdot 2+p_3^\ast\cdot 3}&\underline{m}_{-p^\ast\cdot 1-2p^\ast\cdot 2+p_3^\ast\cdot 3}
\end{array}\right)^T
\]}
Each entry of the corresponding optimization procedure $m\mapsto O(m)$ involves 60 combinations of entries.
\medskip

\noindent 
{\bf Checking that $P$ is not empty.} The polytope $P$ defined by \eqref{DEFP} can be regarded as a polytope $P^\alpha_m$ where the entries of $m$ are given by \eqref{DEFP} together with default values for a polytope inside $[0,\tfrac12]^3$, {\sl eg.}
\[
m=\left(
\begin{array}{cccccccccc}
\tfrac12&\tfrac12&\tfrac12&1&\tfrac12&\tfrac32&1&p_2^\ast&0&0\\
0&\tfrac{\epsilon}2&0&\tfrac{\epsilon}2(3-2\epsilon)&0&0&0&\tfrac{p_3^\ast}2&-\tfrac{p^\ast}2+\tfrac{p_3^\ast}2&-\tfrac{3p^\ast}2+\tfrac{p_3^\ast}2
\end{array}\right)^T
\]
As indicated in Lemma \ref{OPTIVEC}, in order to check that $P^\alpha_m$ is not empty, it suffices to verify that $\underline{O(m)}_i<\overline{O(m)}_i$ for all $i$. 
\medskip

\noindent 
{\bf Checking that $P$ is localized as required.} This amounts to verifying that the entries of $m$ satisfy the condition in Claim \ref{D1}. 
\medskip

\noindent 
{\bf Checking the dynamical conditions.} 

\noindent 
{\sl Checking $G_{3,\epsilon}(P\cap A_{000000})\subset P\cap A_{000101}$.} A preliminary step is to compute the optimized constraint matrices associated with the atoms $A_{000000}$ and $A_{000101}$. This is done by first considering constraint matrices $m_{000000}$ and $m_{000101}$ in Appendix \ref{A-P3} (obtained for the matrix associated with the symbolic partition) and related elementary bounds, and then by applying the optimization procedure.

Once the atoms are characterized, atomic restrictions $P\cap A_{000000}$ and $P\cap A_{000101}$ are computed using intersections and again the optimization procedure. Then, explicit computations show that, in the current context, the inequalities \eqref{INCL000} (which are equivalent to the condition $G_{3,\epsilon}(P\cap A_{000000})\subset P\cap A_{000101}$) hold iff $\epsilon \in \left[\tfrac{5-\sqrt{17}}2,\tfrac12\right)$.
\medskip

\noindent
{\sl Checking $G_{3,\epsilon}(P\cap A_{000001})\subset \sigma_{3124}(P)$.}
The expression of $\sigma_{3124}|_{P^\alpha_m}$ in the previous section also satisfies \eqref{SPETRANS} in the current context and we have 
{\small \begin{align*}
\sigma_{3124}(m)=&\left(
\begin{array}{cccccccc}
1-\underline{m}_{1+2}&\overline{m}_1&\overline{m}_{2+3}&1-\underline{m}_2&\overline{m}_{1+2+3}&1+\overline{m}_3&1+\overline{m}_{1+2\cdot 2+3\cdot 3}&p_1^\ast+\overline{m}_{-p^\ast\cdot 1-2p^\ast\cdot 2+p_3^\ast\cdot 3}\\
1-\overline{m}_{1+2}&\underline{m}_1&\underline{m}_{2+3}&1-\overline{m}_2&\underline{m}_{1+2+3}&1+\underline{m}_3&1+\underline{m}_{1+2\cdot 2+3\cdot 3}&p_1^\ast+\underline{m}_{-p^\ast\cdot 1-2p^\ast\cdot 2+p_3^\ast\cdot 3}
\end{array}\right.\\
&\left.
\begin{array}{cc}
\overline{m}_{p_1^\ast\cdot 1+p_2^\ast\cdot 2+p_3^\ast\cdot 3}-p^\ast&\overline{m}_{-p^\ast\cdot 1+p_2^\ast\cdot 2+p_3^\ast\cdot 3}-p^\ast\\
\underline{m}_{p_1^\ast\cdot 1+p_2^\ast\cdot 2+p_3^\ast\cdot 3}-p^\ast&\underline{m}_{-p^\ast\cdot 1+p_2^\ast\cdot 2+p_3^\ast\cdot 3}-p^\ast
\end{array}
\right)^T
\end{align*}}
Then, explicit computations show that, in the current context, the inequalities \eqref{LAST001} (which are equivalent to the condition $G_{3,\epsilon}(P\cap A_{000001})\subset \sigma_{3124}(P)$) hold for all $\epsilon \in \left[\tfrac{5-\sqrt{17}}2,\tfrac12\right)$.
\medskip

\noindent
{\sl Checking $G_{3,\epsilon}(P\cap A_{000101})\subset \sigma_{2134}(P)$.}
The expression of $\sigma_{2134}|_{P^\alpha_m}$ in the previous section also satisfies \eqref{SPETRANS} in the current context and we have 
{\small \begin{align*}
\sigma_{2134}(m)=&\left(
\begin{array}{ccccccccc}
1-\underline{m}_1&\overline{m}_{1+2}&\overline{m}_3&1+\overline{m}_2&\overline{m}_{1+2+3}&1+\overline{m}_{2+3}&1+\overline{m}_{1+2\cdot 2+3\cdot 3}&p_1^\ast+\overline{m}_{-p^\ast\cdot 1+p_2^\ast\cdot 2+p_3^\ast\cdot 3}\\
1-\overline{m}_1&\underline{m}_{1+2}&\underline{m}_3&1+\underline{m}_2&\underline{m}_{1+2+3}&1+\underline{m}_{2+3}&1+\underline{m}_{1+2\cdot 2+3\cdot 3}&p_1^\ast+\underline{m}_{-p^\ast\cdot 1+p_2^\ast\cdot 2+p_3^\ast\cdot 3}
\end{array}\right.\\
&\left.
\begin{array}{cc}
\overline{m}_{p_1^\ast\cdot 1+p_2^\ast\cdot 2+p_3^\ast\cdot 3}-p^\ast&\overline{m}_{-p^\ast\cdot 1-2p^\ast\cdot 2+p_3^\ast\cdot 3}-p^\ast\\
\underline{m}_{p_1^\ast\cdot 1+p_2^\ast\cdot 2+p_3^\ast\cdot 3}-p^\ast&\underline{m}_{-p^\ast\cdot 1-2p^\ast\cdot 2+p_3^\ast\cdot 3}-p^\ast
\end{array}
\right)^T
\end{align*}}
Then, explicit computations show that, in the current context, the inequalities \eqref{LASTINC} (which are equivalent to the condition $G_{3,\epsilon}(P\cap A_{000101})\subset \sigma_{2134}(P)$) hold for all $\epsilon \in \left[\tfrac{5\sqrt{17}}2,\tfrac12\right)$.
\medskip

\noindent
{\sl Proof that $\text{Orb}_{\langle \sigma_{2134},\sigma_{1324}\rangle}(P)$ is an AsIUP with respect to $\Sigma=1-\text{Id}$ .} Using the expressions of the optimized matrices and of the symmetry-induced transformations on these matrices (and that these transformations all satisfy \eqref{SPETRANS}), one easily checks that the following relations hold 
\begin{align*}
&\sup\left\{x_3\ :\ x\in P\right\}=\sup\left\{x_3\ :\ x\in \sigma_{2134}(P)\right\}=\tfrac{(1-\epsilon)^2}3\\
&\sup\left\{x_3\ :\ x\in \sigma_{1324}(P)\right\}=\sup\left\{x_3\ :\ x\in \sigma_{3124}(P)\right\}=p^\ast\\
&\sup\left\{x_3\ :\ x\in \sigma_{3214}(P)\right\}=\sup\left\{x_3\ :\ x\in \sigma_{2314}(P)\right\}=1-p^\ast
\end{align*}
and $\tfrac{(1-\epsilon)^2}3<p^\ast<1-p^\ast$ for all $\epsilon\in (0,\tfrac12)$. It follows that
\[
\sup\left\{x_3\ :\ x\in \text{Orb}_{\langle \sigma_{2134},\sigma_{1324}\rangle}(P)\right\}=1-p^\ast<1-\tfrac{(1-\epsilon)^2}3 \leq \inf\left\{x_3\ :\ x\in \Sigma\left(\text{Orb}_{\langle \sigma_{2134},\sigma_{1324}\rangle}(P)\right)\right\}
\]
from where the desired conclusion easily follows from the commutation of $\Sigma$ with every transformation in the subgroup $\langle \sigma_{2134},\sigma_{1324}\rangle$. \hfill $\Box$
\medskip

\noindent
{\bf Addendum: sketch of proof of the equality $G_{3,\epsilon}(P\cap A_{000101})= \sigma_{2134}(P)$.}
\begin{itemize}
\item Consider the coefficient matrix obtained from the canonical one above by removing the first row, ie. ignore the constraint on $x_1$.
\item Show that $\sigma_{2134}$ and this new coefficient matrix satisfy the conditions of statement {\sl (ii)} in Claim \ref{CLAIMSPETRANS}.
\item Check finally that, when $\epsilon \in \left[\tfrac{5-\sqrt{17}}2,\tfrac12\right)$, the optimized matrix $O(m)$ associated with $P$ satisfies the equalities obtained from the inequalities \eqref{LASTINC} by replacing $\subset$ by $=$ and $\sigma_{2134}(m)$ by $\sigma_{2134}\circ O(m)$.
\end{itemize}
\end{document}


\title{Supplementary material of the paper "Conditioning problems for invariant sets of expanding piecewise affine mappings: Application to loss of ergodicity in globally coupled maps"}
\author{Bastien Fernandez and Fanni M.\ S\'elley}
\date{}   
\maketitle


\leftline{\small\today.}

\section{Details of the analysis of Conditioning Problem \ref{CONDPRO1}}
This section provides additional details of the analysis of the conditioning problem \ref{CONDPRO1} presented in Appendix \ref{A-P1} of the main text. All notations and labels refer to that appendix uniquely.

\subsection{Simplifying the entries of the optimized matrices $O(m\cap m_{001})$ and $O(m\cap m_{011})$}
This section presents the arguments of the simplification of the expressions of the entries $O(m\cap m_{001})$ and $O(m\cap m_{011})$ obtained from those of $m\cap m_{001}$ and $m\cap m_{011}$, via the expression of the function $O$ in this context, and by using in particular, the restrictions on the entries of $m$ obtained in Claims \ref{BASICm21} and \ref{ma12}.

\subsubsection{Simplifying the expression of $O(m\cap m_{001})$.}\label{SM-001}
\begin{flalign*}
\bullet\ \underline{O(m\cap m_{001})}_1&=\max\left\{\underline{m}_1, \tfrac{\underline{m}_{a\cdot 1+2}-\tfrac12}{a}, \underline{m}_{1+a\cdot 2}-\tfrac{a}2, \tfrac{a\underline{m}_{a\cdot 1+2}-\min\{\overline{m}_{1+a\cdot 2},\tfrac{1+a}2\}}{a^2-1}\right\}&\\
&=\max\left\{\underline{m}_1, \tfrac{\underline{m}_{a\cdot 1+2}-\tfrac12}{a}, \underline{m}_{1+a\cdot 2}-\tfrac{a}2\right\}&
\end{flalign*}
where the second inequality follows from the fact that $\underline{m}_1=\underline{O(m)}_1\geq \tfrac{a\underline{m}_{a\cdot 1+2}-\overline{m}_{1+a\cdot 2}}{a^2-1}$ and $\underline{m}_{a\cdot 1+2}\leq\tfrac{1+a}2$. 
\begin{flalign*}
\bullet\ \underline{O(m\cap m_{001})}_2&=\underline{m}_2&
\end{flalign*}
because the expression of $\underline{O(m\cap m_{001})}_{2}$ only involves components for which  $\overline{(m\cap m_{001})}_i=\overline{m}_i$ and $\underline{(m\cap m_{001})}_i=\underline{m}_i$. (In particular, the list does not involve $\overline{(m\cap m_{001})}_2$ nor $\overline{(m\cap m_{001})}_{1+a\cdot 2}$.)
\begin{flalign*}
\bullet\ \underline{O(m\cap m_{001})}_{a\cdot 1+2}&=\max\left\{\underline{m}_{a\cdot 1+2}, a\underline{m}_1+\underline{m}_2, \tfrac{(a^2-1)\underline{m}_1+\underline{m}_{1+a\cdot 2}}{a}, a\underline{m}_{1+a\cdot 2}-\tfrac{a^2-1}2\right\}&\\
&=\max\left\{\underline{m}_{a\cdot 1+2}, a\underline{m}_{1+a\cdot 2}-\tfrac{a^2-1}2\right\}&
\end{flalign*}
and this maximum remains to be computed.
\begin{flalign*}
\bullet\ \underline{O(m\cap m_{001})}_{1+a\cdot 2}&=\underline{m}_{1+a\cdot 2}&
\end{flalign*}
for the same reason as the one above for $\underline{O(m\cap m_{001})}_2$.
\begin{flalign*}
\bullet\ \overline{O(m\cap m_{001})}_1&=\min\left\{\overline{m}_1, \tfrac{\min\{\overline{m}_{a\cdot 1+2},\tfrac{1+a}2\}-\underline{m}_2}{a}, \min\{\overline{m}_{1+a\cdot 2},\tfrac{1+a}2\}-a\underline{m}_2, \tfrac{a\min\{\overline{m}_{a\cdot 1+2},\tfrac{1+a}2\}-\underline{m}_{1+a\cdot 2}}{a^2-1}\right\}&\\
&=\min\left\{\overline{m}_1, \tfrac{\tfrac{1+a}2-\underline{m}_2}{a}, \tfrac{1+a}2-a\underline{m}_2, \tfrac{\tfrac{a(1+a)}2-\underline{m}_{1+a\cdot 2}}{a^2-1}\right\}
=\overline{m}_1&
\end{flalign*}
where the last equality follows from the inequalities $\underline{m}_2\leq \tfrac12$ and $a\geq 1$, which imply $\tfrac12\leq \tfrac{\tfrac{1+a}2-\underline{m}_2}{a}\leq \tfrac{1+a}2-a\underline{m}_2$ and $\tfrac12\leq \tfrac{\tfrac{a(1+a)}2-\underline{m}_{1+a\cdot 2}}{a^2-1}$, which follows from $\underline{m}_{1+a\cdot 2}\leq\tfrac{1+a}2$. 
\begin{flalign*}
\bullet\ \overline{O(m\cap m_{001})}_2&=\min\left\{\tfrac12, \tfrac{\min\{\overline{m}_{1+a\cdot 2},\tfrac{1+a}2\}-\underline{m}_1}{a}, \overline{m}_{a\cdot 1+2}-a\underline{m}_1, \tfrac{a\min\{\overline{m}_{1+a\cdot 2},\tfrac{1+a}2\}- \underline{m}_{a\cdot 1+2}}{a^2-1}\right\}&\\
&=\min\left\{\tfrac12, \tfrac{\tfrac{1+a}2-\underline{m}_1}{a},\tfrac{\tfrac{a(1+a)}2- \underline{m}_{a\cdot 1+2}}{a^2-1}\right\}=\tfrac12&
\end{flalign*}
\begin{flalign*}
\bullet\ \overline{O(m\cap m_{001})}_{a\cdot 1+2}&=\min\left\{\overline{m}_{a\cdot 1+2}, a\overline{m}_1+\tfrac12, \tfrac{(a^2-1)\overline{m}_1+\min\{\overline{m}_{1+a\cdot 2},\tfrac{1+a}2\}}{a}, a\min\{\overline{m}_{1+a\cdot 2},\tfrac{1+a}2\}-(a^2-1)\underline{m}_2\right\}&\\
&=\min\left\{\overline{m}_{a\cdot 1+2},a\overline{m}_1+\tfrac12, \tfrac{(a^2-1)\overline{m}_1+\tfrac{1+a}2}{a}, \tfrac{a(1+a)}2-(a^2-1)\underline{m}_2\right\}&\\
&=\min\left\{\overline{m}_{a\cdot 1+2},\tfrac{1+a}2\right\}=\overline{m}_{a\cdot 1+2}&
\end{flalign*}
where the third equality uses $\overline{m}_1\leq \tfrac12$ and $\underline{m}_2\leq \tfrac12$, $a\geq 1$ and Claim \ref{ma12} of the main text, and the last one follows from $\overline{m}_1\leq \tfrac12$. 
\begin{flalign*}
\bullet\ \overline{O(m\cap m_{001})}_{1+a\cdot 2}&=\min\left\{\overline{m}_{1+a\cdot 2}, \tfrac{1+a}2,\tfrac{a}2+\overline{m}_1, \tfrac{\tfrac{a^2-1}2+\overline{m}_{a\cdot 1+2}}{a}\right\}=\min\left\{\overline{m}_{1+a\cdot 2}, \tfrac{a}2+\overline{m}_1, \tfrac{\tfrac{a^2-1}2+\overline{m}_{a\cdot 1+2}}{a}\right\}&
\end{flalign*}

\subsubsection{Simplifying the expression of $O(m\cap m_{011})$.}\label{SM-011}
\begin{flalign*}
\bullet\ \underline{O(m\cap m_{011})}_1&=\max\left\{\underline{m}_1, \tfrac{\max\{\underline{m}_{a\cdot 1+2},\tfrac12\}-\overline{m}_2}{a}, \max\{\underline{m}_{1+a\cdot 2},\tfrac{a}2\}-a\overline{m}_2\right\}&\\
&=\max\left\{\underline{m}_1, \tfrac{\tfrac12-\overline{m}_2}{a}, \tfrac{a}2-a\overline{m}_2\right\}=\underline{m}_1&
\end{flalign*}
where we used $\tfrac12\leq \overline{m}_2$ and $\underline{m}_1\geq 0$. 
\begin{flalign*}
\bullet\ \underline{O(m\cap m_{011})}_2&=\max\left\{\tfrac12, \tfrac{\max\{\underline{m}_{1+a\cdot 2},\tfrac{a}2\}-\overline{m}_1}{a}, \underline{m}_{a\cdot 1+2}-a\overline{m}_1, \tfrac{a\max\{\underline{m}_{1+a\cdot 2},\tfrac{a}2\}-\overline{m}_{a\cdot 1+2}}{a^2-1}\right\}&\\
&=\max\left\{\tfrac12,\tfrac12-\tfrac{\overline{m}_1}{a}, \tfrac{\frac{a^2}2-\overline{m}_{a\cdot 1+2}}{a^2-1}\right\}=\tfrac12&
\end{flalign*}
where the second equality relies on $\underline{O(m)}_2=\underline{m}_2<\tfrac12$ and the third one uses $\overline{m}_1\geq 0$ and $\tfrac12\leq \overline{m}_{a\cdot 1+2}$. 
\begin{flalign*}
\bullet\ \underline{O(m\cap m_{011})}_{a\cdot 1+2}&=\max\left\{\underline{m}_{a\cdot 1+2}, a\underline{m}_1+\tfrac12, \tfrac{(a^2-1)\underline{m}_1+\max\{\underline{m}_{1+a\cdot 2},\tfrac{a}2\}}{a}, a\max\{\underline{m}_{1+a\cdot 2},\tfrac{a}2\}-(a^2-1)\overline{m}_2\right\}&\\
&=\max\left\{\underline{m}_{a\cdot 1+2}, a\underline{m}_1+\tfrac12, \tfrac{(a^2-1)\underline{m}_1+\tfrac{a}2}{a}, 2-(a^2-1)\overline{m}_2\right\}&\\
&=\max\left\{\underline{m}_{a\cdot 1+2}, a\underline{m}_1+\tfrac12 \right\}&
\end{flalign*}
where the last reduction uses $2-(a^2-1)\overline{m}_2\leq \tfrac12$ which follows from $\tfrac12\leq \overline{m}_2$.
\begin{flalign*}
\bullet\ \underline{O(m\cap m_{011})}_{1+a\cdot 2}&=\max\left\{\max\{\underline{m}_{1+a\cdot 2},\tfrac{a}2\}, \tfrac{a}2+\underline{m}_1, \tfrac{\tfrac{a^2-1}2+\underline{m}_{a\cdot 1+2}}{a},a\underline{m}_{a\cdot 1+2}-(a^2-1)\overline{m}_1\right\}&\\
&=\max\left\{\underline{m}_{1+a\cdot 2},\tfrac{a}2+\underline{m}_1, \tfrac{\tfrac{a^2-1}2+\underline{m}_{a\cdot 1+2}}{a}\right\}&
\end{flalign*}
\begin{flalign*}
\bullet\ \overline{O(m\cap m_{011})}_1&=\min\left\{\overline{m}_1, \tfrac{\overline{m}_{a\cdot 1+2}-\tfrac12}{a}, \overline{m}_{1+a\cdot 2}-\tfrac{a}2, \tfrac{a\overline{m}_{a\cdot 1+2}-\max\{\underline{m}_{1+a\cdot 2},\tfrac{a}2\}}{a^2-1}\right\}&\\
&=\min\left\{\overline{m}_1, \tfrac{\overline{m}_{a\cdot 1+2}-\tfrac12}{a}, \overline{m}_{1+a\cdot 2}-\tfrac{a}2\right\}&
\end{flalign*}
where we used $\tfrac12\leq \overline{m}_{a\cdot 1+2}$.
\begin{flalign*}
\bullet\ \overline{O(m\cap m_{011})}_2&=\overline{m}_2&
\end{flalign*}
\begin{flalign*}
\bullet\ \overline{O(m\cap m_{011})}_{a\cdot 1+2}&=\min\left\{\overline{m}_{a\cdot 1+2}, a\overline{m}_1+\overline{m}_2, \tfrac{(a^2-1)\overline{m}_1+\overline{m}_{1+a\cdot 2}}{a}, a\overline{m}_{1+a\cdot 2}-\tfrac{a^2-1}2\right\}&\\
&=\min\left\{\overline{m}_{a\cdot 1+2},a\overline{m}_{1+a\cdot 2}-\tfrac{a^2-1}2\right\}&
\end{flalign*}
\begin{flalign*}
\bullet\ \overline{O(m\cap m_{011})}_{1+a\cdot 2}&=\overline{m}_{1+a\cdot 2}&
\end{flalign*}

\subsection{Proof of Claim \ref{ma12Bis}}\label{SM-ma12Bis}
By contradiction, assume that $\underline{m}_{a\cdot 1+2}\leq a(1+a) r_\varrho-\tfrac{a^2-1}2$. Then the second equation before Claim \ref{ma12Bis} in the main text implies
\[
\overline{m}_{1+a\cdot 2}=(1+a)\left(a(1-\epsilon)(1-2r_\varrho)+\epsilon(1-2\varrho)\right).
\]
We have $a(1-\epsilon)(1-2r_\varrho)+\epsilon(1-2\varrho)>\tfrac{r_\varrho}{a}+\tfrac{a-1}{2a}$ for all $a>1$ and $\epsilon,\varrho\in (0,\tfrac12)$; hence the last equation before Claim \ref{ma12Bis} implies (as in the main text)
\[
\underline{m}_{a\cdot 1+2}=(1+a)\left(\tfrac{1-\epsilon}{a}(1-2r_\varrho)+\epsilon(1-2\varrho)\right)
\]
and then, the inequality $\underline{m}_{a\cdot 1+2}\leq a(1+a) r_\varrho-\tfrac{a^2-1}2$ is equivalent to $\epsilon\leq 1-\tfrac{a}2$. Under this condition, using also $\underline{m}_1=\epsilon(1-2\varrho)$ and $\underline{m}_{1+a\cdot 2}=(1+a)r_\varrho$ (and the simplified expression of $\underline{O(m\cap m_{011})}_{1+a\cdot 2}$ in the previous section), the following equality holds 
\[
\underline{O(m\cap m_{011})}_{1+a\cdot 2}=\max\left\{\underline{m}_{1+a\cdot 2},\tfrac{a}2+\underline{m}_1, \tfrac{\tfrac{a^2-1}2+\underline{m}_{a\cdot 1+2}}{a}\right\}=(1+a)r_\varrho.
\]
When evaluated at the entry $\underline{O(m\cap m_{011})}_{1+a\cdot 2}$, the condition $2(1-\epsilon)O(m\cap m_{011})+(2\epsilon-1){0\choose 1}\subset m$ yields the inequality $\underline{m}_{1+a\cdot 2}\geq 2a$, which is equivalent to $a\leq \tfrac{1-2\varrho\epsilon}{5-2(2-\varrho)\epsilon}$. However, one checks that $\tfrac{1-2\varrho\epsilon}{5-2(2-\varrho)\epsilon}\leq \tfrac13$ for all $\epsilon,\varrho\in (0,\tfrac12)$; hence the previous inequality cannot hold. \hfill $\Box$

\subsection{Checking $m_a=O(m_a)$}\label{SM-OPTI}
The definition of the optimization function $O$ in the main text indicates that a constraint matrix $m=O(m)$ is indeed an optimized matrix iff we have 
\begin{align*}
\underline{m}_1&\geq \max\left\{\tfrac{\underline{m}_{a\cdot 1+2}-\overline{m}_2}{a}, \underline{m}_{1+a\cdot 2}-a\overline{m}_2, \tfrac{a\underline{m}_{a\cdot 1+2}-\overline{m}_{1+a\cdot 2}}{a^2-1}\right\}\\
\underline{m}_2&\geq\max\left\{\tfrac{\underline{m}_{1+a\cdot 2}-\overline{m}_1}{a}, \underline{m}_{a\cdot 1+2}-a\overline{m}_1, \tfrac{a\underline{m}_{1+a\cdot 2}-\overline{m}_{a\cdot 1+2}}{a^2-1}\right\}\\
\underline{m}_{a\cdot 1+2}&\geq\max\left\{ a\underline{m}_1+\underline{m}_2, \tfrac{(a^2-1)\underline{m}_1+\underline{m}_{1+a\cdot 2}}{a}, a\underline{m}_{1+a\cdot 2}-(a^2-1)\overline{m}_2\right\}\\
\underline{m}_{1+a\cdot 2}&\geq\max\left\{a\underline{m}_2+\underline{m}_1, \tfrac{(a^2-1)\underline{m}_2+\underline{m}_{a\cdot 1+2}}{a},a\underline{m}_{a\cdot 1+2}-(a^2-1)\overline{m}_1\right\}\\
\overline{m}_1&\leq\min\left\{\tfrac{\overline{m}_{a\cdot 1+2}-\underline{m}_2}{a}, \overline{m}_{1+a\cdot 2}-a\underline{m}_2, \tfrac{a\overline{m}_{a\cdot 1+2}-\underline{m}_{1+a\cdot 2}}{a^2-1}\right\}\\
\overline{m}_2&\leq\min\left\{\tfrac{\overline{m}_{1+a\cdot 2}-\underline{m}_1}{a}, \overline{m}_{a\cdot 1+2}-a\underline{m}_1, \tfrac{a\overline{m}_{1+a\cdot 2}-\underline{m}_{a\cdot 1+2}}{a^2-1}\right\}\\
\overline{m}_{a\cdot 1+2}&\leq\min\left\{ a\overline{m}_1+\overline{m}_2, \tfrac{(a^2-1)\overline{m}_1+\overline{m}_{1+a\cdot 2}}{a}, a\overline{m}_{1+a\cdot 2}-(a^2-1)\underline{m}_2\right\}\\
\overline{m}_{1+a\cdot 2}&\leq \min\left\{a\overline{m}_2+\overline{m}_1, \tfrac{(a^2-1)\overline{m}_2+\overline{m}_{a\cdot 1+2}}{a},a\overline{m}_{a\cdot 1+2}-(a^2-1)\underline{m}_1\right\}
\end{align*}
For the matrix $m_a$ whose entries are recalled here, 
{\small \begin{equation*}
m_a=\left(
\begin{array}{cccc}
r_{\varrho}&1-2\epsilon\left(1-\varrho-\epsilon(1-2\varrho)\right)&(1+a)r_{\varrho}&\tfrac{1+a}{a}\left(a r_\varrho+2(a-1)(1-\epsilon)^2(1-2r_\varrho)\right)\\
\epsilon(1-2\varrho)&r_{\varrho}&\tfrac{1+a}{a}\left(r_\varrho+(a-1)\epsilon(1-2\varrho)\right)&(1+a)r_{\varrho}
\end{array}\right)^T
\end{equation*}}
these inequalities simplify as follows after elementary algebra, using $r_\varrho-\epsilon(1-2\varrho)=(1-\epsilon)(1-2r_\varrho)$ 
and the conditions on the parameters
\begin{align*}
0&\geq \max\left\{\tfrac{1-2a(1-\epsilon))}{a^2}, 1-2a(1-\epsilon), \tfrac{-2(1-\epsilon)}{a}\right\}\\
0&\geq\max\left\{0, -1,0\right\}\\
\tfrac{1}{a}&\geq\max\left\{ 0,\tfrac{1}{a},a-2 (a^2-1)(1-\epsilon)\right\}\\
0&\geq\max\left\{-1, -\tfrac{a^2-1}{a^2},-(a^2-1)\right\}\\
0&\leq\min\left\{0, 1, 0\right\}\\
0&\leq\min\left\{\tfrac{a-2(1-\epsilon)}{a^2}, a-2(1-\epsilon),\tfrac{1}{a}\right\}\\
0&\leq\min\left\{ 1, \tfrac{a^2-1}{a^2}, (a^2-1)\right\}\\
-\tfrac{2(1-\epsilon)}{a}&\leq \min\left\{\tfrac{-2(1-\epsilon)}{1+a}, \tfrac{-2(1-\epsilon)}{a},-2(1-\epsilon)+a-1\right\}
\end{align*}
These inequalities evidently hold for $\epsilon\in [\epsilon_a,\tfrac12)$ since we then have $\tfrac1{2a}<1-\epsilon< \tfrac{a}2$. Moreover, according to statement {\sl (iii)} of Lemma 2.1 in the main text, those inequalities that are strict, namely the first, fourth, sixth and seventh ones, correspond to edges of $P^{\alpha_a}_{m_a}$, namely these edges are
\[
x_1=\underline{m_a}_1,\ x_1+ax_2=\underline{m_a}_{1+a\cdot 2},\ x_2=\overline{m_a}_2,\ \text{and}\quad ax_1+x_2=\overline{m_a}_{a\cdot 1+2}.
\]

\section{Details of the analysis of Conditioning Problem \ref{CONDPRO3}}
This section provides additional details of the analysis of Conditioning Problem \ref{CONDPRO3} presented in Appendix \ref{A-P2} of the main text. All notations and labels refer to that appendix uniquely.

\subsection{Simplifying the entries of the optimized constraint matrices associated with atomic restrictions}
This section presents the arguments of the simplification of the expressions of the entries  of the optimized constraint matrices associated with atomic restrictions, via the expression of the function $O$ in this context, and by using in particular, the restrictions on the entries of $m_1$ and $m_2$ obtained in Claim \ref{BASICm1m2}.

\subsubsection{Simplifying the expression of $O(m_1\cap m_{000101})$.}
The definition of the intersection of two constraint matrices together with the inequalities in Claim \ref{BASICm1m2} of the main text imply 
\[
m_1\cap m_{000101}=\left(
\begin{array}{cccccc}
\tfrac12&\overline{m_1}_2&\overline{m_1}_3&\min\{\overline{m_1}_{1+2},1\}&\overline{m_1}_{2+3}&\min\{\overline{m_1}_{1+2+3},1\}\\
\underline{m_1}_1&\underline{m_1}_2&\underline{m_1}_3&\underline{m_1}_{1+2}&\underline{m_1}_{2+3}&\underline{m_1}_{1+2+3}
\end{array}\right)^T.
\]
The computation of $O(m_1\cap m_{000101})$ then proceeds as follows.
\begin{flalign*}
\bullet\ \underline{O(m_1\cap m_{000101})}_1&=\underline{m_1}_1&
\end{flalign*}
because the expression of $\underline{O(m_1\cap m_{000101})}_{1}$ only involves components for which  $\overline{(m_1\cap m_{000101})}_{i,j}=\overline{m_1}_{i,j}$ and $\underline{(m_1\cap m_{000101})}_{i,j}=\underline{m_1}_{i,j}$.
\begin{flalign*}
\bullet\ \underline{O(m_1\cap m_{000101})}_2&=\max\left\{\underline{m_1}_2,\underline{m_1}_{1+2}-\tfrac12,\underline{m_1}_{1+2}+\underline{m_1}_{2+3}-\min\{\overline{m_1}_{1+2+3},1\}\right\}&\\
&=\max\left\{\underline{m_1}_2,\underline{m_1}_{1+2}-\tfrac12,\underline{m_1}_{1+2}+\underline{m_1}_{2+3}-1\}\right\}=\max\left\{\underline{m_1}_2,\underline{m_1}_{1+2}-\tfrac12\right\}&
\end{flalign*}
where the second equality is a consequence of $\underline{m_1}_2=\underline{O(m_1)}_2$ and the last one follows from $\underline{m_1}_{2+3}\leq \tfrac12$, which is a consequence of Claim \ref{SYM-M2}.
\begin{flalign*}
\bullet\ \underline{O(m_1\cap m_{000101})}_3&=\max\left\{\underline{m_1}_3,\underline{m_1}_{1+2+3}-\min\{\overline{m_1}_{1+2},1\},\underline{m_1}_{1+2+3}-\tfrac12-\overline{m_1}_2\right\}&\\
&=\max\left\{\underline{m_1}_3,\underline{m_1}_{1+2+3}-1,\underline{m_1}_{1+2+3}-\tfrac12-\overline{m_1}_2\right\}=\max\left\{\underline{m_1}_3,\underline{m_1}_{1+2+3}-\tfrac12-\overline{m_1}_2\right\}
\end{flalign*}
where the last equality follows from $\underline{m_1}_{1+2+3}<1$.
\begin{flalign*}
\bullet\ \underline{O(m_1\cap m_{000101})}_{1+2}&=\underline{m_1}_{1+2}&
\end{flalign*}
\begin{flalign*}
\bullet\ \underline{O(m_1\cap m_{000101})}_{2+3}&=\max\left\{\underline{m_1}_{2+3},\underline{m_1}_{1+2+3}-\tfrac12,\underline{m_1}_2+\underline{m_1}_{1+2+3}-\min\{\overline{m_1}_{1+2},1\}\right\}&\\
&=\max\left\{\underline{m_1}_{2+3},\underline{m_1}_{1+2+3}-\tfrac12,\underline{m_1}_2+\underline{m_1}_{1+2+3}-1\right\}
=\max\left\{\underline{m_1}_{2+3},\underline{m_1}_{1+2+3}-\tfrac12\right\}&
\end{flalign*} 
using $\underline{m_1}_2\leq\tfrac12$.
\begin{flalign*}
\bullet\ \underline{O(m_1\cap m_{000101})}_{1+2+3}&=\underline{m_1}_{1+2+3}&
\end{flalign*}
\begin{align*}
\bullet\ \overline{O(m_1\cap m_{000101})}_1&=\min\left\{\tfrac12,\min\{\overline{m_1}_{1+2},1\}-\underline{m_1}_2,\min\{\overline{m_1}_{1+2+3},1\}-\underline{m_1}_{2+3},\underline{m_1}_{1+2+3}-\overline{m_1}_2-\overline{m_1}_3\right.\\
&\left.\overline{m_1}_3+\min\{\overline{m_1}_{1+2},1\}-\underline{m_1}_{2+3}\right\}=\min\left\{\tfrac12,1-\underline{m_1}_2,1-\underline{m_1}_{2+3},\overline{m_1}_3+1-\underline{m_1}_{2+3}\right\}=\tfrac12&
\end{align*} 
using $\tfrac12<\overline{m_1}_1=\overline{O(m_1)}_1$ and $\underline{m_1}_2,\underline{m_1}_{2+3}<\tfrac12$ and $0<\overline{m_1}_3$.
\begin{flalign*}
\bullet\ \overline{O(m_1\cap m_{000101})}_2&=\min\left\{\overline{m_1}_2,\min\{\overline{m_1}_{1+2},1\}-\underline{m_1}_1,\min\{\overline{m_1}_{1+2+3},1\}-\underline{m_1}_1-\underline{m_1}_3,\right.&\\
&\left. \min\{\overline{m_1}_{1+2},1\}+\overline{m_1}_{2+3}-\underline{m_1}_{1+2+3}\right\}&\\
&=\min\left\{\overline{m_1}_2,1-\underline{m_1}_1,1-\underline{m_1}_1-\underline{m_1}_3,1+\overline{m_1}_{2+3}-\underline{m_1}_{1+2+3}\right\}&\\
&=\min\left\{\overline{m_1}_2,1-\underline{m_1}_1-\underline{m_1}_3\right\}&
\end{flalign*}
\begin{flalign*}
\bullet\ \overline{O(m_1\cap m_{000101})}_3&=\min\left\{\overline{m_1}_3,\min\{\overline{m_1}_{1+2+3},1\}-\underline{m_1}_{1+2},\tfrac12+\overline{m_1}_{2+3}-\underline{m_1}_{1+2}\right\}&\\
&=\min\left\{\overline{m_1}_3,1-\underline{m_1}_{1+2},\tfrac12+\overline{m_1}_{2+3}-\underline{m_1}_{1+2}\right\}=\min\left\{\overline{m_1}_3,\tfrac12+\overline{m_1}_{2+3}-\underline{m_1}_{1+2}\right\}&
\end{flalign*}
where the last equality is a consequence of $\overline{m_1}_3,\overline{m_1}_{2+3}\leq \tfrac12$.
\begin{flalign*}
\bullet\ \overline{O(m_1\cap m_{000101})}_{1+2}&=\min\left\{\overline{m_1}_{1+2},1,\tfrac12+\overline{m_1}_2,\min\{\overline{m_1}_{1+2+3},1\}-\underline{m_1}_3,\overline{m_1}_2+\min\{\overline{m_1}_{1+2+3},1\}-\underline{m_1}_{2+3}\right\}&\\
&=\min\left\{\overline{m_1}_{1+2},\tfrac12+\overline{m_1}_2,\overline{m_1}_2+1-\underline{m_1}_{2+3}\right\}=\min\left\{\overline{m_1}_{1+2},\tfrac12+\overline{m_1}_2\right\}&
\end{flalign*}
using $\overline{m_1}_2,\underline{m_1}_{2+3}<\tfrac12$ and $0\leq  \underline{m_1}_3$.
\begin{flalign*}
\bullet\ \overline{O(m_1\cap m_{000101})}_{2+3}&=\min\left\{\overline{m_1}_{2+3},\min\{\overline{m_1}_{1+2+3},1\}-\underline{m_1}_{1},\overline{m_1}_3+\min\{\overline{m_1}_{1+2},1\}-\underline{m_1}_{1},\right.&\\
&\left. \overline{m_1}_2+\min\{\overline{m_1}_{1+2+3},1\}-\underline{m_1}_{1+2}\right\}&\\
&=\min\left\{\overline{m_1}_{2+3},1-\underline{m_1}_{1},\overline{m_1}_3+1-\underline{m_1}_{1},\overline{m_1}_2+1-\underline{m_1}_{1+2}\right\}&\\
&=\min\left\{\overline{m_1}_{2+3},1-\underline{m_1}_{1}\right\}=\overline{m_1}_{2+3}
\end{flalign*}
using $\underline{m_1}_{1},\overline{m_1}_{2+3}<\tfrac12$.
\begin{flalign*}
\bullet\ \overline{O(m_1\cap m_{000101})}_{1+2+3}&=\min\left\{\overline{m_1}_{1+2+3},1,\tfrac12+\overline{m_1}_{2+3},\overline{m_1}_3+\min\{\overline{m_1}_{1+2},1\},\min\{\overline{m_1}_{1+2},1\}+\overline{m_1}_{2+3}-\underline{m_1}_2\right\}&\\
&=\min\left\{\overline{m_1}_{1+2+3},1,\tfrac12+\overline{m_1}_{2+3},\overline{m_1}_3+1,1+\overline{m_1}_{2+3}-\underline{m_1}_2\right\}
=\min\left\{\overline{m_1}_{1+2+3},\tfrac12+\overline{m_1}_{2+3}\right\}
\end{flalign*} 
using $\underline{m_1}_2,\overline{m_1}_{2+3}<\tfrac12$ and $0<\overline{m_1}_3$.

\subsubsection{Simplifying the expression of $O(m_1\cap m_{100101})$.}
The definition of the intersection of two constraint matrices together with the inequalities in Claim \ref{BASICm1m2} of the main text imply 
\[
m_1\cap m_{100101}=\left(
\begin{array}{cccccc}
\overline{m_1}_1&\overline{m_1}_2&\overline{m_1}_3&\overline{m_1}_{1+2}&\overline{m_1}_{2+3}&\overline{m_1}_{1+2+3}\\
\tfrac12&\underline{m_1}_2&\underline{m_1}_3&\underline{m_1}_{1+2}&\underline{m_1}_{2+3}&\underline{m_1}_{1+2+3}
\end{array}\right)^T.
\]
The computation of $O(m_1\cap m_{100101})$ is rather straightforward. It suffices to identify, those combinations in the definition of $O$ that involve $\underline{m}_1$ and use that $m_1=O(m_1)$. We then have
\begin{align*}
&\underline{O(m_1)}_1=\max\left\{\tfrac12,\ \underline{m}_{1+2}-\overline{m}_2,\ \underline{m}_{1+2+3}-\overline{m}_{2+3}\right\}\quad\text{and}\quad\overline{O(m_1)}_1=\overline{m}_1\\
&\underline{O(m_1)}_2=\underline{m_1}_2\quad\text{and}\quad\overline{O(m_1)}_2=\min\left\{\overline{m}_2,\ \overline{m}_{1+2}-\tfrac12\right\}\\
&\underline{O(m_1)}_3=\max\left\{\underline{m}_3,\ \tfrac12+\underline{m}_{2+3}-\overline{m}_{1+2}\right\}\quad\text{and}\quad\overline{O(m_1)}_3=\min\left\{\overline{m}_3,\ \overline{m}_{1+2+3}-\tfrac12-\underline{m}_2\right\}\\
&\underline{O(m_1)}_{1+2}=\max\left\{\underline{m}_{1+2},\ \tfrac12+\underline{m}_2\right\}\quad\text{and}\quad\overline{O(m_1)}_{1+2}=\overline{m}_{1+2}\\
&\underline{O(m_1)}_{2+3}=\underline{m}_{2+3}\quad\text{and}\quad\overline{O(m_1)}_{2+3}=\min\left\{\overline{m}_{2+3},\ \overline{m}_{1+2+3}-\tfrac12\right\}\\
&\underline{O(m_1)}_{1+2+3}=\max\left\{\underline{m}_{1+2+3},\ \tfrac12+\underline{m}_{2+3}\right\}\quad\text{and}\quad\overline{O(m_1)}_{1+2+3}=\overline{m}_{1+2+3}
\end{align*}


\subsubsection{Simplifying the expression of $O(m_2\cap m_{110111})$.}
The definition of the intersection of two constraint matrices together with the inequalities in Claim \ref{BASICm1m2} of the main text imply 
\[
m_2\cap m_{110111}=\left(
\begin{array}{cccccc}
\overline{m_2}_1&\overline{m_2}_2&\overline{m_2}_3&\tfrac32&\min\{\overline{m_2}_{2+3},1\}&\tfrac32\\
\underline{m_2}_1&\tfrac12&\underline{m_2}_3&\max\{\underline{m_2}_{1+2},1\}&\tfrac12&\max\{\underline{m_2}_{1+2+3},1\}
\end{array}\right)^T,
\]
The computation of $O(m_2\cap m_{110111})$ then proceeds as follows.
\begin{flalign*}
\bullet\ \underline{O(m_2\cap m_{110111})}_1&=\max\left\{\underline{m_2}_1, \max\{\underline{m_2}_{1+2},1\}-\overline{m_2}_2,\max\{\underline{m_2}_{1+2+3},1\}-\overline{m_2}_{2+3}, \right.&\\
&\left. \max\{\underline{m_2}_{1+2+3},1\}-\overline{m_2}_{2}-\overline{m_2}_{3},\underline{m_2}_3+\max\{\underline{m_2}_{1+2},1\}-\overline{m_2}_{2+3}\right\}&\\
&=\max\left\{\underline{m_2}_1, 1-\overline{m_2}_2, 1-\overline{m_2}_{2+3},1-\overline{m_2}_{2}-\overline{m_2}_{3}, \underline{m_2}_3+1-\overline{m_2}_{2+3}\right\}&\\
&=\max\left\{\underline{m_2}_1, 1-\overline{m_2}_2, \underline{m_2}_3+1-\overline{m_2}_{2+3}\right\}=\max\left\{\underline{m_2}_1, 1-\overline{m_2}_2\right\}=\underline{m_2}_1&
\end{flalign*}
where the second equality is a consequence of $\underline{m_2}_1=\underline{O(m_2)}_1$ and the other ones use $0<\overline{m_2}_{3}$, $\overline{m_2}_2=\overline{O(m_2)}_2$ and $1-\overline{m_2}_2\leq\tfrac12\leq\underline{m_2}_1$.
\begin{flalign*}
\bullet\ \underline{O(m_2\cap m_{110111})}_2&=\max\left\{\tfrac12, \max\{\underline{m_2}_{1+2},1\}-\overline{m_2}_1, \tfrac12-\overline{m_2}_{3}, \max\{\underline{m_2}_{1+2+3},1\}-\overline{m_2}_{1}-\overline{m_2}_{3},\right.&\\
&\left. \max\{\underline{m_2}_{1+2},1\}-1\right\}=\max\left\{\tfrac12, 1-\overline{m_2}_1, \tfrac12-\overline{m_2}_{3},1-\overline{m_2}_{1}-\overline{m_2}_{3}, \underline{m_2}_{1+2}-1\right\}=\tfrac12.
\end{flalign*}
using $\underline{O(m_2)}_2=\underline{m_2}_2<\tfrac12$, $\tfrac12<\overline{m_2}_1$, $0<\overline{m_2}_{3}$ and $\underline{m_2}_{1+2}\leq\tfrac32$.
\begin{flalign*}
\bullet\ \underline{O(m_2\cap m_{110111})}_3&=\max\left\{\underline{m_2}_3, \tfrac12-\overline{m_2}_2, \max\{\underline{m_2}_{1+2+3},1\}-\tfrac32, \max\{\underline{m_2}_{1+2+3},1\}-\overline{m_2}_{1}-\overline{m_2}_{2}, \underline{m_2}_1-1\right\}&\\
&=\max\left\{\underline{m_2}_3, \tfrac12-\overline{m_2}_2, \underline{m_2}_{1+2+3}-\tfrac32, \tfrac12,1-\overline{m_2}_{1}-\overline{m_2}_{2}, \underline{m_2}_1-1\right\}=\underline{m_2}_3&
\end{flalign*}
using $\tfrac12<\overline{m_2}_{1}$, $\underline{m_2}_{1+2+3}\leq \tfrac32$, $\tfrac12\leq \overline{m_2}_2$ and $\underline{m_2}_1<1$.
\begin{flalign*}
\bullet\ \underline{O(m_2\cap m_{110111})}_{1+2}&=\max\left\{\underline{m_2}_{1+2},1, \underline{m_2}_1+\tfrac12, \max\{\underline{m_2}_{1+2+3},1\}-\overline{m_2}_3, \underline{m_2}_1+\tfrac12-\overline{m_2}_{3},\right.&\\
&\left. \tfrac12+\max\{\underline{m_2}_{1+2+3},1\}-\overline{m_2}_{2+3}\right\}&\\
&=\max\left\{\underline{m_2}_{1+2},\underline{m_2}_1+\tfrac12, \tfrac12+\underline{m_2}_{1+2+3},-\overline{m_2}_{2+3},\tfrac32-\overline{m_2}_{2+3}\right\}=\max\left\{\underline{m_2}_{1+2}, \underline{m_2}_1+\tfrac12\right\}&
\end{flalign*}
using $\tfrac12\leq \underline{m_2}_1,\overline{m_2}_{2+3}$ and $\underline{O(m_2)}_1=\underline{m_2}_1$.
\begin{flalign*}
\bullet\ \underline{O(m_2\cap m_{110111})}_{2+3}&=\max\left\{\tfrac12, \tfrac12+\underline{m_2}_3, \max\{\underline{m_2}_{1+2+3},1\}-\overline{m_2}_1, \underline{m_2}_3+\max\{\underline{m_2}_{1+2},1\}-\overline{m_2}_{1}, \right.&\\
&\left.\max\{\underline{m_2}_{1+2+3},1\}-1\right\}&\\
&=\max\left\{\tfrac12+\underline{m_2}_3, 1-\overline{m_2}_1, \underline{m_2}_3+1-\overline{m_2}_{1}, \underline{m_2}_{1+2+3}-1,0\right\}=\tfrac12+\underline{m_2}_3&
\end{flalign*}
using $\tfrac12<\overline{m_2}_1$ and $\underline{m_2}_{1+2+3}\leq\tfrac32$.
\begin{flalign*}
\bullet\ \underline{O(m_2\cap m_{110111})}_{1+2+3}&=\max\left\{\underline{m_2}_{1+2+3},1, \underline{m_2}_1+\tfrac12, \underline{m_2}_3+\max\{\underline{m_2}_{1+2},1\}, \underline{m_2}_{1}+\tfrac12+\underline{m_2}_3, \right.&\\
&\left. \max\{\underline{m_2}_{1+2},1\}+\tfrac12-\overline{m_2}_{2}\right\}&\\
&=\max\left\{\underline{m_2}_{1+2+3},1, \underline{m_2}_{1}+\tfrac12, \underline{m_2}_3+1,\underline{m_2}_{1}+\tfrac12+ \underline{m_2}_3,\underline{m_2}_{1+2}+\tfrac12-\overline{m_2}_{2},\tfrac32-\overline{m_2}_{2}\right\}&\\
&=\max\left\{\underline{m_2}_{1+2+3}, \underline{m_2}_{1}+\tfrac12+ \underline{m_2}_3\right\}&
\end{flalign*}
using $\tfrac12\leq \underline{m_2}_{1}$ and $\underline{O(m_2)}_1=\underline{m_2}_1$ and $\tfrac12\leq \overline{m_2}_{2}$.
\begin{flalign*}
\bullet\ \overline{O(m_2\cap m_{110111})}_1&=\min\left\{\overline{m_2}_1,1,1, 1-\underline{m_2}_{3}, \overline{m_2}_3+1\right\}=\overline{m_2}_1&
\end{flalign*}
using the symmetry in Claim \ref{SYM-M2} of the main text.
\begin{flalign*}
\bullet\ \overline{O(m_2\cap m_{110111})}_2&=\min\left\{\overline{m_2}_2, \tfrac32-\underline{m_2}_1, \tfrac32-\underline{m_2}_{1}-\underline{m_2}_{3}, \tfrac32+\overline{m_2}_{2+3}-\max\{\underline{m_2}_{1+2+3},1\}\right\}&\\
&=\min\left\{\overline{m_2}_2, \tfrac32-\underline{m_2}_{1}-\underline{m_2}_{3}\right\}&
\end{flalign*}
using $\underline{O(m_2)}_1=\underline{m_2}_1$ and $0\leq \underline{m_2}_{3}$.
\begin{flalign*}
\bullet\ \overline{O(m_2\cap m_{110111})}_3&=\min\left\{\overline{m_2}_3, \overline{m_2}_{2+3}-\tfrac12, \tfrac32-\max\{\underline{m_2}_{1+2},1\}, 1-\underline{m_2}_{1}, \overline{m_2}_1+\overline{m_2}_{2+3}-\max\{\underline{m_2}_{1+2},1\}\right\}&\\
&=\min\left\{\overline{m_2}_3, \overline{m_2}_{2+3}-\tfrac12, \tfrac32-\underline{m_2}_{1+2}, \tfrac12,1-\underline{m_2}_{1}, \overline{m_2}_1+\overline{m_2}_{2+3}-1\right\}&\\
&=\min\left\{\overline{m_2}_3, \overline{m_2}_{2+3}-\tfrac12\right\}
\end{flalign*}
using $\overline{m_2}_3\leq \tfrac12\leq \underline{m_2}_{1}$ and symmetries in Claim \ref{SYM-M2}.
\begin{flalign*}
\bullet\ \overline{O(m_2\cap m_{110111})}_{1+2}&=\min\left\{\tfrac32, \overline{m_2}_1+\overline{m_2}_2, \tfrac32-\underline{m_2}_3, \overline{m_2}_1+\overline{m_2}_{2+3}-\underline{m_2}_{3}, \overline{m_2}_2+1\right\}=\tfrac32-\underline{m_2}_3&
\end{flalign*}
using $\tfrac32\leq \overline{m_2}_{1+2}=\overline{O(m_2)}_{1+2}$ and $\tfrac12\leq\overline{m_2}_2$.
\begin{flalign*}
\bullet\ \overline{O(m_2\cap m_{110111})}_{2+3}&=\min\left\{\overline{m_2}_{2+3}, \tfrac32-\underline{m_2}_1, \overline{m_2}_3+\tfrac32-\underline{m_2}_{1}, \overline{m_2}_2+\tfrac32-\max\{\underline{m_2}_{1+2},1\}\right\}&\\
&=\min\left\{\overline{m_2}_{2+3}, \tfrac32-\underline{m_2}_1, \overline{m_2}_2+\tfrac12,\overline{m_2}_2+\tfrac32-\underline{m_2}_{1+2}\right\}&\\
&=\min\left\{\overline{m_2}_{2+3}, \tfrac32-\underline{m_2}_1, \overline{m_2}_2+\tfrac12\right\}=\min\left\{\overline{m_2}_{2+3}, \tfrac32-\underline{m_2}_1\right\}&
\end{flalign*}
using $\underline{O(m_2)}_1=\underline{m_2}_1$ and $\overline{m_2}_3\leq\tfrac12$ and $\underline{O(m_2)}_2=\underline{m_2}_2$.
\begin{flalign*}
\bullet\ \overline{O(m_2\cap m_{110111})}_{1+2+3}&=\min\left\{\tfrac32, \overline{m_2}_1+\overline{m_2}_{2+3}, \overline{m_2}_3+\tfrac32, \overline{m_2}_{1}+\overline{m_2}_{2}+\overline{m_2}_3, \tfrac32+\overline{m_2}_{2+3}-\underline{m_2}_{2}\right\}&\\
&=\min\left\{\tfrac32,\overline{m_2}_3+\tfrac32,\tfrac32+\overline{m_2}_{2+3}-\underline{m_2}_{2}\right\}=\min\left\{\tfrac32,\overline{m_2}_3+\tfrac32\right\}=\tfrac32&
\end{flalign*}
using $\tfrac32\leq \overline{m_2}_{1+2+3}=\overline{O(m_2)}_{1+2+3}$ and $\overline{m_2}_{3}=\overline{O(m_2)}_{3}$. 

\subsubsection{Simplifying the expression of $O(m_2\cap m_{110212})$.}
The definition of the intersection of two constraint matrices together with the inequalities in Claim \ref{BASICm1m2} of the main text imply 
\[
m_2\cap m_{110111}=\left(
\begin{array}{cccccc}
\overline{m_2}_1&\overline{m_2}_2&\overline{m_2}_3&\overline{m_2}_{1+2}&\overline{m_2}_{2+3}&\overline{m_2}_{1+2+3}\\
\underline{m_2}_1&\tfrac12&\underline{m_2}_3&\tfrac32&\tfrac12&\tfrac32
\end{array}\right)^T,
\]
The computation of $O(m_2\cap m_{110212})$ then proceeds as follows.
\begin{flalign*}
\bullet\ \underline{O(m_2\cap m_{110212})}_1&=\max\left\{\underline{m_2}_1,\tfrac32-\overline{m_2}_{2},\tfrac32-\overline{m_2}_{2+3}\right\}=\max\left\{\underline{m_2}_1,\tfrac32-\overline{m_2}_{2}\right\}&
\end{flalign*}
because $\overline{m_2}_{2}\leq \overline{m_2}_{2+3}$.
\begin{flalign*}
\bullet\ \underline{O(m_2\cap m_{110212})}_2&=\max\left\{\tfrac12,\tfrac32-\overline{m_2}_1,\tfrac12-\overline{m_2}_{3},\tfrac32-\overline{m_2}_{1}-\overline{m_2}_{3},2-\overline{m_2}_{1+2+3}\right\}=\tfrac32-\overline{m_2}_1&
\end{flalign*}
using $\tfrac12<\overline{m_2}_1$, $0<\overline{m_2}_{3}$ and $\tfrac32\leq \overline{m_2}_{1+2+3}$.
\begin{flalign*}
\bullet\ \underline{O(m_2\cap m_{110212})}_3&=\max\left\{\underline{m_2}_3,\tfrac12-\overline{m_2}_{1},\tfrac32-\overline{m_2}_{1+2}\right\}=\underline{m_2}_3&
\end{flalign*}
using  $\tfrac12<\overline{m_2}_{1}$ and $\tfrac32\leq \overline{m_2}_{1+2}$.
\begin{flalign*}
\bullet\ \underline{O(m_2\cap m_{110212})}_{1+2}&=\max\left\{\tfrac32,\underline{m_2}_1+\tfrac12,\tfrac32-\overline{m_2}_3,2-\overline{m_2}_{2+3}\right\}=\tfrac32&
\end{flalign*}
using $0<\overline{m_2}_3$ and $\tfrac12\leq \overline{m_2}_{2+3} $.
\begin{flalign*}
\bullet\ \underline{O(m_2\cap m_{110212})}_{2+3}&=\max\left\{\tfrac12,\tfrac12+\underline{m_2}_3,\tfrac32-\overline{m_2}_1,\underline{m_2}_3+\tfrac32-\overline{m_2}_{1},2-\overline{m_2}_{1+2}\right\}=\underline{m_2}_3+\tfrac32-\overline{m_2}_{1}&
\end{flalign*}
using $\tfrac32\leq \overline{m_2}_{1+2}$.
\begin{flalign*}
\bullet\ \underline{O(m_2\cap m_{110212})}_{1+2+3}&=\max\left\{\tfrac32,\underline{m_2}_1+\tfrac12,\underline{m_2}_3+\tfrac32,\underline{m_2}_{1}+\tfrac12+\underline{m_2}_3,2-\overline{m_2}_{2}\right\}=\underline{m_2}_3+\tfrac32&
\end{flalign*}
using $\underline{m_2}_1\leq 1$, $0\leq \underline{m_2}_3$ and $\tfrac12\leq \overline{m_2}_{2}$.
\begin{flalign*}
\bullet\ \overline{O(m_2\cap m_{110212})}_1&=\min\left\{\overline{m_2}_1,\overline{m_2}_{1+2}-\tfrac12\right\}=\min\left\{\overline{m_2}_1,\overline{m_2}_{1+2}-\tfrac12\right\}=\overline{m_2}_1&
\end{flalign*}
using $\overline{m_2}_1\leq 1$ and $\tfrac32\leq \overline{m_2}_{1+2}$.
\begin{flalign*}
\bullet\ \overline{O(m_2\cap m_{110212})}_2&=\min\left\{\overline{m_2}_2,\overline{m_2}_{1+2}+\overline{m_2}_{2+3}-\tfrac32\right\}=\overline{m_2}_2&
\end{flalign*}
using $\overline{m_2}_2\leq \overline{m_2}_{2+3}$ and $\tfrac32\leq \overline{m_2}_{1+2}$.
\begin{flalign*}
\bullet\ \overline{O(m_2\cap m_{110212})}_3&=\min\left\{\overline{m_2}_3,\overline{m_2}_{2+3}-\tfrac12,\overline{m_2}_{1+2+3}-\tfrac32\right\}=\min\left\{\overline{m_2}_3,\overline{m_2}_{1+2+3}-\tfrac32\right\}&
\end{flalign*}
because $\overline{m_2}_{2+3}\leq \overline{m_2}_{1+2+3}$.
\begin{flalign*}
\bullet\ \overline{O(m_2\cap m_{110212})}_{1+2}&=\min\left\{\overline{m_2}_{1+2},\overline{m_2}_2+\overline{m_2}_{1+2+3}-\tfrac12\right\}=\overline{m_2}_{1+2}&
\end{flalign*}
because $\overline{m_2}_2\geq \tfrac12$ and $\overline{m_2}_{1+2+3}\geq \overline{m_2}_{1+2}$. 
\begin{flalign*}
\bullet\ \overline{O(m_2\cap m_{110212})}_{2+3}&=\min\left\{\overline{m_2}_{2+3},\overline{m_2}_2+\overline{m_2}_{1+2+3}-\tfrac32\right\}&
\end{flalign*}
\begin{flalign*}
\bullet\ \overline{O(m_2\cap m_{110212})}_{1+2+3}&=\min\left\{\overline{m_2}_{1+2+3},\overline{m_2}_{1+2}+\overline{m_2}_{2+3}-\tfrac12\right\}=\overline{m_2}_{1+2+3}&
\end{flalign*}
because  $\overline{m_2}_{1+2+3}\leq \overline{m_2}_{2+3}+\overline{m_2}_1\leq \overline{m_2}_{2+3}+1\leq \overline{m_2}_{1+2}+\overline{m_2}_{2+3}-\tfrac12$.

\subsubsection{Simplifying the expression of $O(m_2\cap m_{110112})$.}
The definition of the intersection of two constraint matrices together with the inequalities in Claim \ref{BASICm1m2} of the main text imply 
\[
m_2\cap m_{110212}=\left(
\begin{array}{cccccc}
1&\overline{m_2}_2&\overline{m_2}_3&\tfrac32&\overline{m_2}_{2+3}&\overline{m_2}_{1+2+3}\\
\underline{m_2}_1&\tfrac12&0&\underline{m_2}_{1+2}&\tfrac12&\tfrac32
\end{array}\right)^T,
\]
The computation of $O(m_2\cap m_{110112})$ then proceeds as follows (using also the relations in Claim \ref{SYM-M2}).
\begin{align*}
&\underline{O(m_2\cap m_{110112})}_1=\max\left\{\underline{m_2}_1,\tfrac32-\overline{m_2}_{2+3}\right\}\\
&\underline{O(m_2\cap m_{110112})}_2=\max\left\{\tfrac12,\ \underline{m_2}_{1+2}-1\right\}=\tfrac12\\
&\underline{O(m_2\cap m_{110112})}_3=\max\left\{0,\ \tfrac12-\overline{m_2}_2\right\}=0\\
&\underline{O(m_2\cap m_{110112})}_{1+2}=\max\left\{\underline{m_2}_{1+2},\ \underline{m_2}_1+\tfrac12,\ \tfrac32-\overline{m_2}_3,\ 2-\overline{m_2}_{2+3}\right\}=\max\left\{\underline{m_2}_{1+2},\ \underline{m_2}_1+\tfrac12\right\}\\
&\underline{O(m_2\cap m_{110112})}_{2+3}=\max\left\{\tfrac12,\ \underline{m_2}_{1+2}-1\right\}=\tfrac12\\
&\underline{O(m_2\cap m_{110112})}_{1+2+3}=\max\left\{\tfrac32,\ \underline{m_2}_1+\tfrac12,\ \underline{m_2}_{1+2}\right\}=\tfrac32
\end{align*}
and
\begin{align*}
&\overline{O(m_2\cap m_{110112})}_1=\min\left\{1,\ \overline{m_2}_{1+2+3}-\tfrac12\right\}=1\\
&\overline{O(m_2\cap m_{110112})}_2=\min\left\{\overline{m_2}_2,\ \tfrac32-\underline{m_2}_1\right\}\\
&\overline{O(m_2\cap m_{110112})}_3=\min\left\{\overline{m_2}_3,\ \overline{m_2}_{2+3}-\tfrac12,\ 1+\overline{m_2}_{2+3}-\underline{m_2}_{1+2}\right\}=\min\left\{\overline{m_2}_3,\ \overline{m_2}_{2+3}-\tfrac12\right\}\\
&\overline{O(m_2\cap m_{110112})}_{1+2}=\min\left\{\tfrac32,\ 1+\overline{m_2}_2,\ \overline{m_2}_{1+2+3},\ \overline{m_2}_2+\overline{m_2}_{1+2+3}-\tfrac12\right\}=\tfrac32\\
&\overline{O(m_2\cap m_{110112})}_{2+3}=\min\left\{\overline{m_2}_{2+3}, \overline{m_2}_3+\tfrac32-\underline{m_2}_1\right\}\\
&\overline{O(m_2\cap m_{110112})}_{1+2+3}=\min\left\{\overline{m_2}_{1+2+3},\ 1+\overline{m_2}_{2+3},\ \overline{m_2}_3+\tfrac32\right\}
\end{align*}

\subsection{Proof of Claim \ref{m1p2}}
By contradiction, assume that we had $\underline{m_2}_{1+2}<\underline{m_2}_1+\tfrac12$ and then $\overline{m_2}_3+\tfrac12<\overline{m_2}_{2+3}$. Then, the equalities induced on the entries $\underline{O(m_2\cap m_{110111})}_{1+2}$ and $\overline{O(m_2\cap m_{110111})}_{2+3}$ by the condition $G_{3,\epsilon}(P_{2}\cap A_{110111})=P_{1}$ would yield
\[
\left\{\begin{array}{l}
\underline{m_1}_{1+2}=2(1-\epsilon)\underline{m_2}_1+\tfrac{3\epsilon}2-1\\
\overline{m_1}_{2+3}=2(1-\epsilon)(1-\underline{m_2}_1)+\tfrac{\epsilon}2
\end{array}\right.
\]
Moreover, the equality induced on the entry $\overline{O(m_2\cap m_{110111})}_2$ would read $2(1-\epsilon)\overline{m_2}_3=\overline{m_1}_3$. Since $2(1-\epsilon)>1$, the inequality on $\overline{O(m_1\cap m_{000101})}_3$ induced $G_{3,\epsilon}(P_1\cap A_{000101})\subset P_2$ would imply that we ought to have $\tfrac12+\overline{m_1}_{2+3}-\underline{m_1}_{1+2}<\overline{m_1}_3$ and then this inequality would read 
\[
2(1-\epsilon)(\tfrac12+\overline{m_1}_{2+3}-\underline{m_1}_{1+2})\leq \overline{m_2}_3=1-\underline{m_2}_1
\]
Using the expressions above, simple algebra show that we would necessarily have
\[
\underline{m_2}_1\geq \tfrac{6\epsilon^2-13\epsilon+6}{8\epsilon^2-16\epsilon+7}.
\]
However, this inequality is incompatible with the inequality $\underline{m_2}_1\leq 2(1-\epsilon)\underline{m_1}_1+\epsilon\leq \epsilon (3-2\epsilon)$, hence the contradiction. The latter inequality follows from using $\underline{m_1}_1\leq \epsilon$ as the inequality induced at the entry $\underline{O(m_1\cap m_{100101})}_{1}$ from the condition $G_{3,\epsilon}(P_{1}\cap A_{100101})=P_{1}$.

\subsection{Checking $m_1=O(m_1)$ and $m_2=O(m_2)$}
The definition of the optimization function $O$ in the main text indicates that a constraint matrix $m=O(m)$ is indeed an optimized matrix iff we have\footnote{in particular, those combinations of three elements in the definition of $O(m)$ can be removed because they are dominated by combinations of two elements, when the inequalities are considered altogether.}
\begin{align*}
\underline{m}_1&\geq \max\left\{\underline{m}_{1+2}-\overline{m}_2,\ \underline{m}_{1+2+3}-\overline{m}_{2+3}\right\}\\
\underline{m}_2&\geq \max\left\{\underline{m}_{1+2}-\overline{m}_1,\ \underline{m}_{2+3}-\overline{m}_3\right\}\\
\underline{m}_3&\geq \max\left\{\underline{m}_{2+3}-\overline{m}_2,\ \underline{m}_{1+2+3}-\overline{m}_{1+2}\right\}\\
\underline{m}_{1+2}&\geq \max\left\{\underline{m}_1+\underline{m}_2,\ \underline{m}_{1+2+3}-\overline{m}_3\right\}\\
\underline{m}_{2+3}&\geq \max\left\{\underline{m}_2+\underline{m}_3,\ \underline{m}_{1+2+3}-\overline{m}_1\right\}\\
\underline{m}_{1+2+3}&\geq \max\left\{\underline{m}_1+\underline{m}_{2+3},\ \underline{m}_3+\underline{m}_{1+2}\right\}\\
\overline{m}_1&\leq \min\left\{\overline{m}_{1+2}-\underline{m}_2,\ \overline{m}_{1+2+3}-\underline{m}_{2+3}\right\}\\
\overline{m}_2&\leq \min\left\{\overline{m}_{1+2}-\underline{m}_1,\ \overline{m}_{2+3}-\underline{m}_3\right\}\\
\overline{m}_3&\leq \min\left\{\overline{m}_{2+3}-\underline{m}_2,\ \overline{m}_{1+2+3}-\underline{m}_{1+2}\right\}\\
\overline{m}_{1+2}&\leq \min\left\{\overline{m}_1+\overline{m}_2,\ \overline{m}_{1+2+3}-\underline{m}_3\right\}\\
\overline{m}_{2+3}&\leq \min\left\{\overline{m}_2+\overline{m}_3,\ \overline{m}_{1+2+3}-\underline{m}_1\right\}\\
\overline{m}_{1+2+3}&\leq \min\left\{\overline{m}_1+\overline{m}_{2+3},\ \overline{m}_3+\overline{m}_{1+2}\right\}
\end{align*}
For the matrix $m_1$ whose entries are recalled here 
\begin{align*}
m_1=&\left(
\begin{array}{cccc}
1-\epsilon&p^\ast-2(1-\epsilon)\delta_2&p^\ast-\tfrac{\epsilon}2-2(1-\epsilon)\delta_3&1-\tfrac{\epsilon}2\\
1-2p^\ast+2(1-\epsilon)\delta_1&\tfrac{\epsilon}2&0&1-p^\ast+2(1-\epsilon)\delta_3
\end{array}\right.\nonumber\\
&\left.
\begin{array}{cc}
p^\ast-2(1-\epsilon)\delta_3&1-\tfrac{\epsilon}2\\
\tfrac{\epsilon}2&1-p^\ast+2(1-\epsilon)\delta_4
\end{array}\right)^T
\end{align*}
these inequalities bring, after elementary algebra, the following additional limitations on the parameters,
\[
\delta_3+\max\left\{\delta_2,\delta_4\right\}\leq \delta_1\quad\text{and}\quad \delta_3\leq \min\{\delta_2,\delta_4\} 
\]
For the matrix $m_2$ whose entries are recalled here
\[
m_2=\left(
\begin{array}{cccccc}
1&1-p^\ast-\delta_2&1-2p^\ast-\delta_1&2-p^\ast-\delta_5&1-p^\ast-\delta_3&2-p^\ast-\delta_4\\
2p^\ast+\delta_1&p^\ast+\delta_2&0&1+p^\ast+\delta_3&p^\ast+\delta_5&1+p^\ast+\delta_4
\end{array}\right)^T
\]
these inequalities bring, after elementary algebra, the following additional limitations on the parameters,
\[
\delta_5+\delta_1-\min\{\delta_2,\delta_4\}\leq 1-2p^\ast
\]

\section{Details of the computations related to the proof of Statement \ref{SOLP4}}
see Mahematica notebook "AsIUPG3SecondBif.nb".